\documentclass[12pt,twoside]{memoireuqam1.3}

\usepackage{graphicx}
\usepackage[french]{babel}
\usepackage[utf8]{inputenc} 
\usepackage[T1]{fontenc}
\usepackage{url}
\usepackage{amsmath, amsfonts, amssymb, amsthm}
\usepackage{tikz, tkz-graph, tikz-cd}
\usetikzlibrary{arrows, automata, backgrounds, calc, shapes, patterns, fit, shadows.blur}
\usepackage[textwidth=2cm]{todonotes}
\usepackage{enumitem}  
\usepackage{mathtools}  
\usepackage{genyoungtabtikz} 
\usepackage{aliascnt}  
\usepackage[colorlinks = true, linkcolor = blue, urlcolor  = blue, citecolor = blue, anchorcolor = black]{hyperref}
\usepackage{cleveref}
\usepackage{sectsty}  
\usepackage{shuffle}  
\usepackage{float}  
\usepackage[normalem]{ulem}  
\usepackage{cases}  
\usepackage{stackengine}  
\usepackage[most]{tcolorbox}
\usepackage[font=footnotesize]{caption}
\usepackage{arydshln}  
\usepackage{adjustbox}  
\usepackage{afterpage}  
\usepackage{etoolbox}  
                    
\NoAutoSpaceBeforeFDP

\usepackage{imakeidx}
\makeindex[intoc]
\patchcmd{\tableofcontents}{\indexname}{INDEX DES NOTATIONS}{}{}

\usepackage{regexpatch}
\makeatletter
\xpatchcmd{\@todo}{\setkeys{todonotes}{#1}}{\setkeys{todonotes}{inline,#1}}{}{}
\makeatother

\allsectionsfont{\normalfont\bfseries}
\chapterfont{\normalfont\bfseries\centering\raggedbottom\textsc}

\makeatletter
\def\th@remark{%
	\thm@headfont{\bfseries}%
	\normalfont 
	\thm@preskip\topsep \divide\thm@preskip\tw@
	\thm@postskip\thm@preskip
}
\makeatother

\newtheorem{thm}{Théorème}
\newaliascnt{lem}{thm}
\newaliascnt{ex}{thm}
\newaliascnt{prop}{thm}
\newaliascnt{corl}{thm}
\newaliascnt{rmq}{thm}
\newaliascnt{defn}{thm}
\newaliascnt{cj}{thm}
\newaliascnt{pb}{thm}
\newtheorem{corl}[corl]{Corollaire}
\newtheorem{lem}[lem]{Lemme}
\newtheorem{prop}[prop]{Proposition}
\newtheorem{cj}[cj]{Conjecture}
\theoremstyle{remark}
\newtheorem{defn}[defn]{Définition}
\newtheorem{rmq}[rmq]{Remarque}
\newtheorem{ex}[ex]{Exemple}
\newtheorem{pb}[pb]{Problème}
\aliascntresetthe{lem}
\aliascntresetthe{prop}
\aliascntresetthe{corl}
\aliascntresetthe{rmq}
\aliascntresetthe{defn}
\aliascntresetthe{ex}
\aliascntresetthe{cj}
\aliascntresetthe{pb}

\addto\extrasfrench{}
\DeclareCaptionLabelSeparator{colon}{~: }

\newcommand{\C}{\mathbb{C}}
\newcommand{\N}{\mathbb{N}}

\newcommand{\CSn}{\C S_n}
\DeclareMathOperator{\sh}{\mathsf{sh}}
\DeclareMathOperator{\del}{\partial}

\DeclareMathOperator{\isoproj}{isoproj}
\DeclareMathOperator{\projlift}{projlift}
\DeclareMathOperator{\proj}{proj}
\DeclareMathOperator{\eig}{eig}
\DeclareMathOperator{\diag}{diag}

\DeclareMathOperator{\type}{type}

\newcommand{\tabw}[1]{t_{#1}}
\DeclareMathOperator{\Id}{Id}
\DeclareMathOperator{\tr}{tr}
\DeclareMathOperator{\im}{im}

\DeclareMathOperator{\noninv}{noninv}

\newcommand{\bhr}{\textsf{BHR}}
\newcommand{\A}{\C\langle A\rangle}  

\newcommand{\dia}{\YFrench \yng}
\newcommand{\smalldia}{\YFrench \Yboxdim{5pt}\yng}
\newcommand{\tab}{\YFrench \young}

\newcommand{\Nu}[1]{\nu_{#1}}
\newcommand{\allnuk}{\{\Nu{k}\}_{k \in \mathbb{N}}}

\newcommand{\gammak}[1]{\gamma_{#1}}
\newcommand{\allgammak}{\{\gammak{k}\}_{k \in \mathbb{N}}}

\newcommand{\prodnuk}{\prod_{j\in [n]} \#\{l \in [k] \mid a_l= j\}!}
\newcommand{\shdelUnak}{\sh_{a_1}\circ\ldots\circ\sh_{a_k}\circ\del_{a_1}\circ\ldots\circ\del_{a_k}}

\newcommand{\eqexp}[2]{\overset{\mathclap{\normalfont\mbox{\tiny #1}}}{=\qquad}#2}
\newcommand{\bbox}[1]{\adjustbox{cfbox=blue}{$#1$}}

\newcommand{\createcard}[1]{%
	\begin{tikzpicture}
	\node [ inner sep=0pt, outer sep=0pt,
	clip, rounded corners=3pt,
	blur shadow={shadow blur steps=5} ]
	(image) at (0,0) {
		\includegraphics[scale=0.225, trim={1pt, 1.5pt, 1pt, 1.5pt}, clip]{#1}
	};
	\node [draw=black, thin, fit=(image), rounded corners=3pt, inner sep=0pt] {};
	\end{tikzpicture}%
	\hspace{5pt}%
}
\def\carda{\createcard{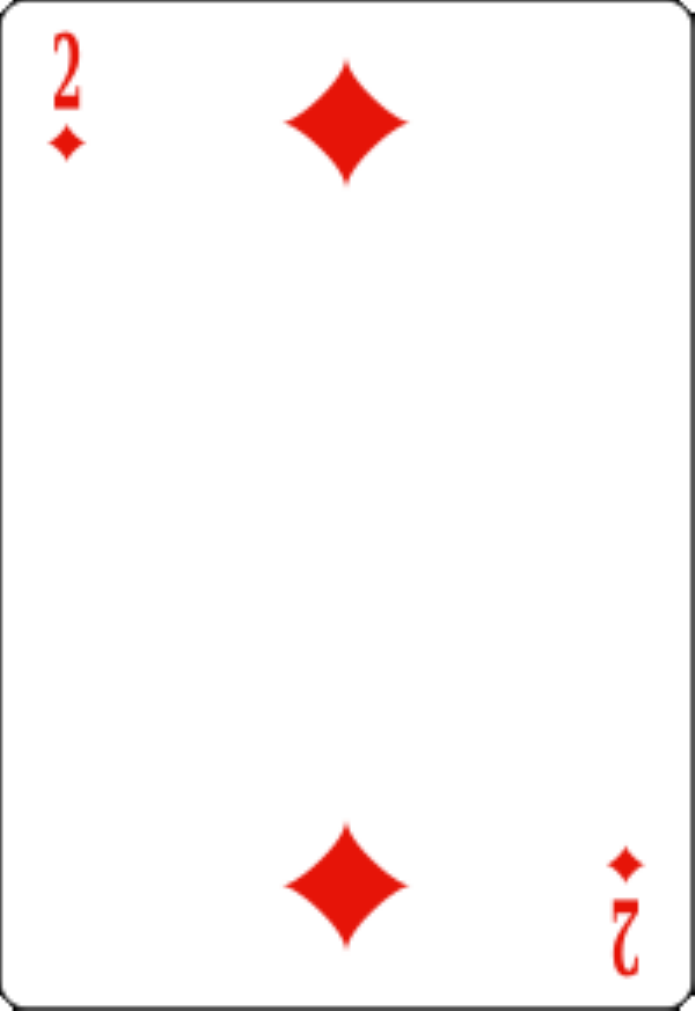}}
\def\cardb{\createcard{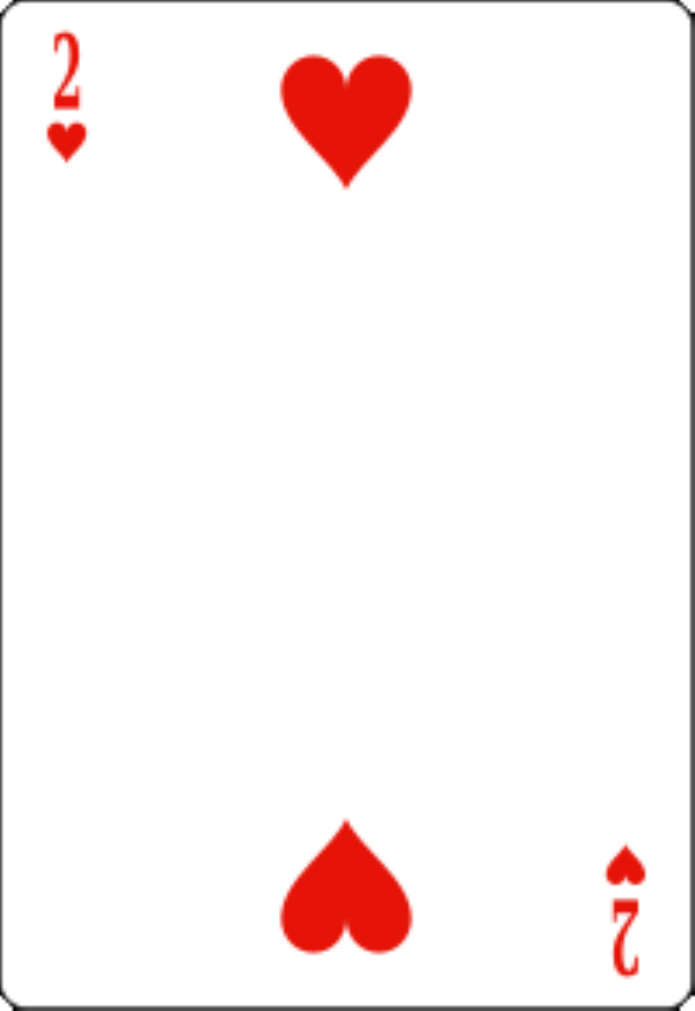}}
\def\cardc{\createcard{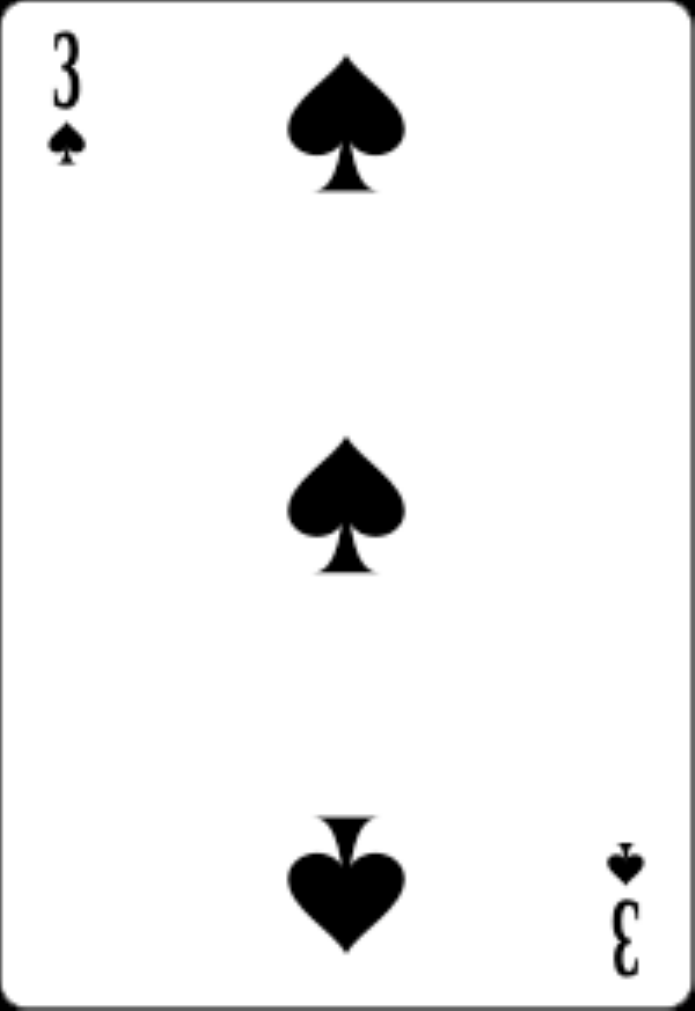}}
\def\cardd{\createcard{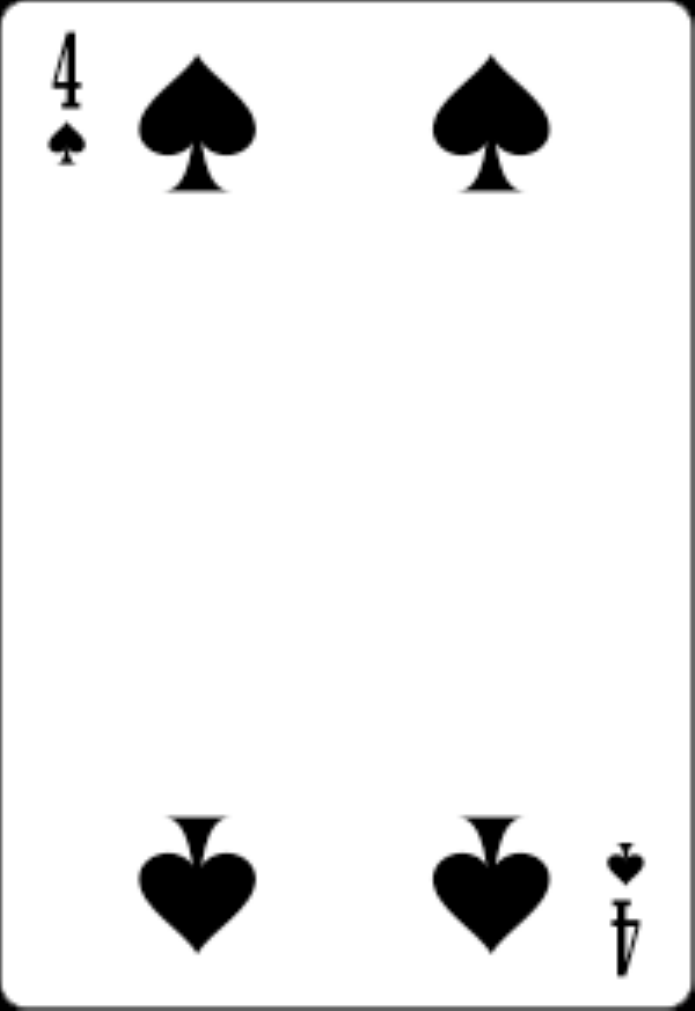}}
\def\carde{\createcard{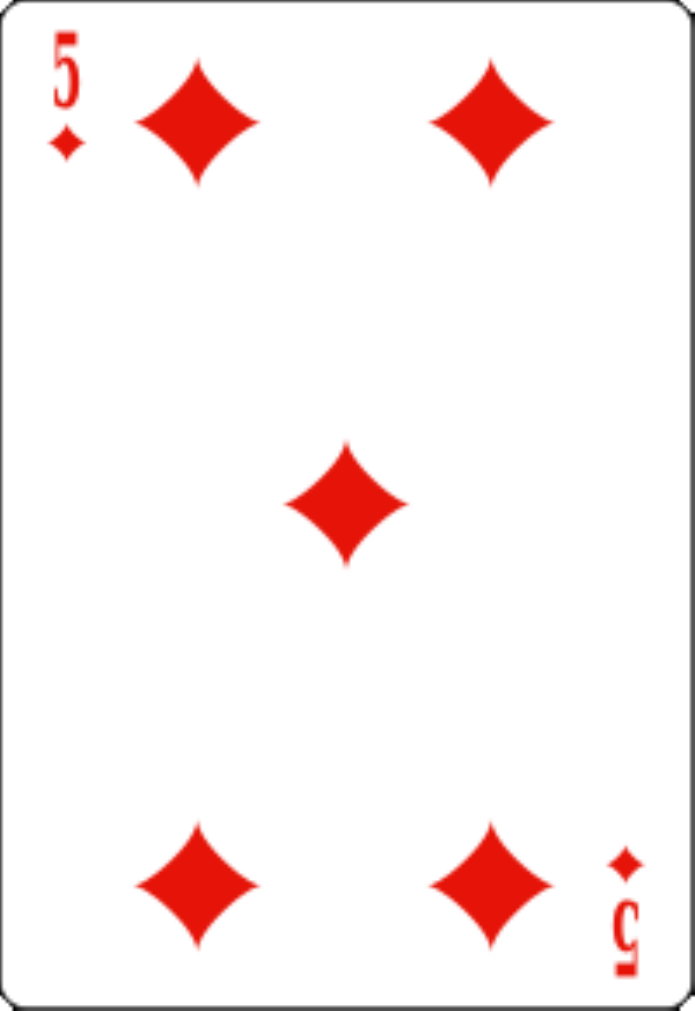}}
\def\cardf{\createcard{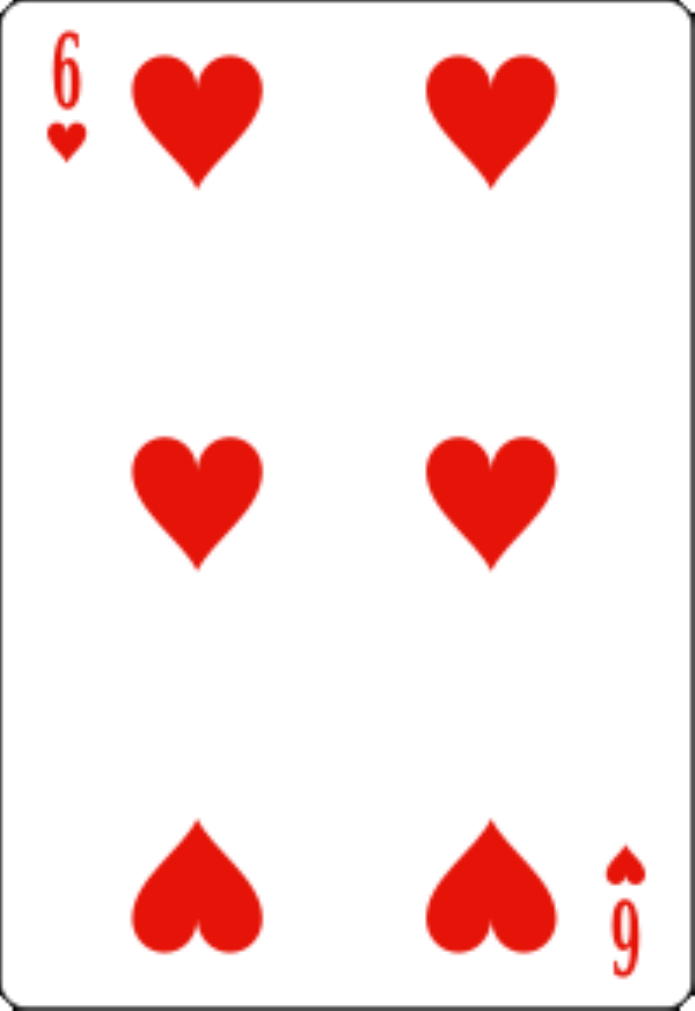}}
\def\cardg{\createcard{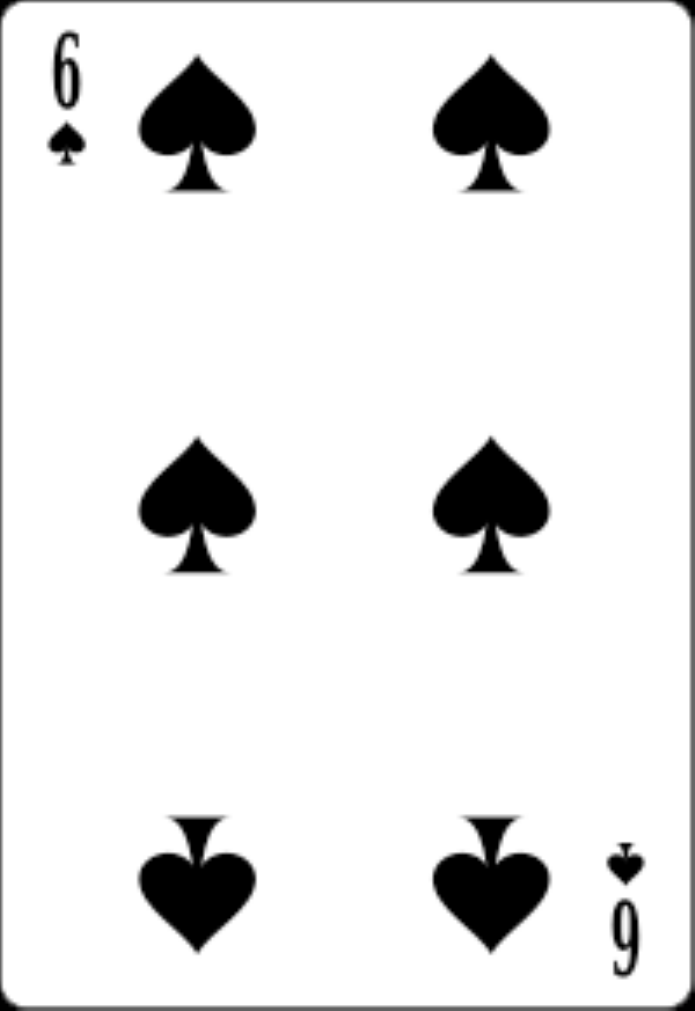}}

\def\versocarte{\createcard{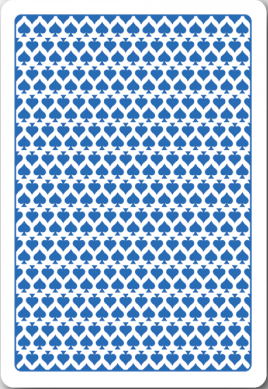}}
\def\tuileScrabble{
\begin{tikzpicture}
	\node [clip] (image) at (0,0) {\includegraphics[scale=0.2]{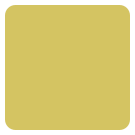}};
\end{tikzpicture}%
\hspace{5pt}}
\newcommand{\Scrabble}[2]{
	\begin{tikzpicture}
	\node (image) at (0,0) {\tuileScrabble};
	\node (lettre) at (0.1,0)  {{\Large \textsf{#1$_{_{#2}}$}}};
	\end{tikzpicture}
}

\allowdisplaybreaks
\begin{document}
\sloppy
\title{Valeurs propres des opérateurs de mélange symétrisés}
\author{Nadia Lafrenière}
\uqamthese 
\matiere{mathématiques}

\thispagestyle{empty}        
\maketitle

\renewcommand \bibname{R\'ef\'erences}

\renewcommand \listfigurename{Liste des figures}
\renewcommand \appendixname{Annexe}
\renewcommand \figurename{Figure}
\renewcommand \tablename{Tableau}
\def\tableautorefname{tableau}
\def\appendixautorefname{annexe}

\pagenumbering{roman} 
\addtocounter{page}{1} 
\chapter*{Remerciements}
Les premiers remerciements vont sans conteste à mon directeur, Franco Saliola, qui a cru en moi tout au long de ce projet et qui m'a encouragé dans toutes mes (bonnes et moins bonnes) initiatives. Tu es un excellent motivateur, toujours de bonne humeur et avec un rire contagieux. Merci pour le sujet passionnant, tes patientes explications et ta disponibilité. Tes intérêts diversifiés en mathématiques et ton souci de bien communiquer les mathématiques ont su nourrir ma passion des mathématiques durant les quatre dernières années (et plus), et je t'en suis très reconnaissante. Merci aussi pour tes nombreux conseils et pour le code qui a permis l'élaboration des conjectures.

Merci à ma grande amie Stéphanie Schanck avec qui une partie du travail de cette thèse a été réalisée. Ta curiosité, ta passion des maths, ton énergie et ta bonne humeur sont si inspirantes! Nos rencontres de recherche avec Franco me manquent déjà.
	
Merci aux membres du jury d'avoir pris le temps de  lire ma thèse. Votre intérêt et vos commentaires ont été éclairants. Plus personnellement, merci à Greg Warrington, François Bergeron et Christophe Reutenauer.

Merci à Megan Bernstein, sans qui je n'aurais pas pu écrire la section sur le temps de mélange.

Les Fonds de recherche Québec et l'Institut des sciences mathématiques m'ont apporté le support financier essentiel pour réaliser mon projet de doctorat. Merci de supporter la recherche fondamentale et l'avancement des connaissances.

La très belle communauté du LaCIM ne serait pas la même sans les professeurs, les étudiantes et les étudiants. Merci pour les mots croisés, les discussions de mathématiques et les bons moments. Des remerciements particuliers vont à celles qui partagent mon bureau, pour tout le support durant la rédaction de la thèse~: Pauline Hubert, Mélodie Lapointe et Florence Maas-Gariépy. Merci également à Aram Dermenjian, Fanny Desjardins, Herman Goulet-Ouellet et Émile Nadeau. 

Un grand merci à Laura Colmenarejo, pour ses encouragements soutenus au cours de mon doctorat.

Enfin, merci à ma famille et à mes amies et amis d'être des personnes si extraordinaires!
	
\tableofcontents 
\listoftables 
\listoffigures 
\begin{abstract}
	L'opérateur de mélange doublement aléatoire explique, par exemple, l'évolution d'un paquet de cartes à jouer au fil de l'exécution de la procédure consistant à piger une carte aléatoirement et à la remettre à une position aléatoire. 
	Si on pige plusieurs cartes avant de toutes les replacer, on obtient plutôt une famille d'opérateurs de mélange symétrisés, tels que définis par Victor Reiner, Franco Saliola et Volkmar Welker.
	
	Cette thèse décrit comment obtenir les valeurs propres de ces opérateurs de mélange. Pour y arriver, on s'appuie notamment sur les travaux d'Anton Dieker et de Franco Saliola, qui ont calculé les valeurs propres du mélange doublement aléatoire (qui est parmi les mélanges étudiés). Ici, on calcule les valeurs propres pour tous les opérateurs de la famille. On procède avec l'aide de la théorie de la représentation du groupe symétrique~: on décompose l'espace vectoriel sur lequel les mélanges agissent (dont une base est formée des permutations des cartes) en modules simples. Ceux-ci sont désignés par les tableaux standards, et l'algorithme de calcul des valeurs propres est entièrement combinatoire~: elles sont obtenues à partir de calculs sur les tableaux standards.
	
	Les techniques utilisées ici permettent de démontrer d'une nouvelle façon que ces opérateurs de mélange commutent entre eux et de démontrer que les valeurs propres sont entières et positives, résolvant une conjecture laissée ouverte par Reiner, Saliola et Welker. De plus, connaître les valeurs propres pourrait notamment permettre de déterminer le nombre nécessaire d'itérations d'un mélange pour qu'un paquet de cartes soit bien mélangé. 
	
	Un second ensemble de mélanges aussi introduit par Reiner, Saliola et Welker est étudié. Nous présentons plusieurs conjectures concernant leurs valeurs propres.

\textbf{Mots-clefs~:} Combinatoire algébrique, théorie de la représentation du groupe symétrique, chaînes de Markov, tableaux standards, mélanges de cartes, valeurs propres, opérateurs de mélange symétrisés
	
\end{abstract}

\renewcommand\abstractname{English abstract}
\begin{abstract}
\textbf{English title:} Eigenvalues of Symmetrized Shuffling Operators

The random-to-random shuffling operator explains, for example, the evolution of a deck of cards subject to the following random process: draw a card randomly from the deck and reinsert it at a random position.
If one instead draws more than one card at a time before reinserting, then the resulting operator is an example of a family of symmetrized shuffling operators studied by Victor Reiner, Franco Saliola and Volkmar Welker.

This thesis describes a way to obtain the eigenvalues of these operators.
We build on the work of Anton Dieker and Franco Saliola, who computed the eigenvalues of the random-to-random shuffle.
Here, we compute the eigenvalues for all the operators of the family.
We proceed with the help of the representation theory of the symmetric group.
We decompose the vector space on which the shuffles act into simple modules for the symmetric group.
These modules correspond to standard Young tableaux, and the algorithm to compute the eigenvalues is combinatorial because it computes the eigenvalues directly from the standard Young tableaux.

As a corollary of our main result, we solve several conjectures of Reiner, Saliola and Welker, including showing that the eigenvalues are all nonnegative integers.
Furthermore, the techniques used here allow us to give a new proof of their result that these symmetrized shuffling operators commute.
Knowing the eigenvalues is the key step in one method of computing the number of shuffles one needs to execute to get a perfectly shuffled deck, which is briefly explored.

We also study a second family of shuffles introduced by Reiner, Saliola and Welker. We present many conjectures about their eigenvalues.

\textbf{Keywords:} Algebraic combinatorics, representation theory of the symmetric group, Markov chains, standard tableaux, card shuffling, eigenvalues, symmetrized shuffling operators
\end{abstract}



\begin{introduction}
	Supposons que vous êtes en famille, que vous jouez à un jeu de cartes, et que vous venez de terminer une manche. Pour que les prochains tours soient équitables, il faut que les cartes soient bien mélangées. Après tout, peut-être que votre soeur a une mémoire phénoménale, qui lui permettrait d'avoir mémorisé l'ordre dans lequel les cartes ont été rangées; elle serait alors largement avantagée! 
	
	Alors, comment s'assurer que les cartes soient bien mélangées? C'est une question qui taraude la communauté mathématique depuis qu'elle a été posée par Henri Poincaré \cite{poincare} et étudiée vers la moitié du \mbox{XX$^e$ siècle} \cite{borel,gilbert}. Si les développements ont été sommaires au début de l'étude des mélanges de cartes, ils ont explosé dans les années 1980, notamment avec les travaux de Persi Diaconis. Encore aujourd'hui, c'est un sujet d'actualité; de nombreux articles sont écrits sur le sujet.
	
	Cette thèse porte sur certains mélanges introduits par Victor Reiner, Franco Saliola et Volkmar Welker qu'on appelle les mélanges symétrisés. Nous étudions deux de ces types de mélange~:
	\begin{enumerate}
		\item le mélange doublement aléatoire s'exécute en retirant une carte, puis en la réinsérant n'importe où. Au lieu de retirer une seule carte, on peut en piger plusieurs (un nombre déterminé d'avance), les mélanger, puis les réinsérer.
		\item on pourrait aussi diviser tout le paquet en un nombre fixé de paquets de $1$ et de paquets de $2$ cartes, puis battre les paquets ensemble pour n'en obtenir qu'un.
	\end{enumerate}
	Bien que ces mélanges semblent très distincts, Reiner, Saliola et Welker ont illustré qu'ils ont en commun d'être associés à des matrices dont les valeurs propres sont \emph{jolies} (c'est-à-dire entières et qui s'expriment potentiellement de façon combinatoire). De plus, les opérateurs de la première famille commutent deux à deux, tout comme ceux de la deuxième famille, une propriété rare parmi les mélanges symétrisés, dont font partie ceux présentés plus haut.
	
	Le calcul des valeurs propres est laborieux, notamment parce que les matrices associées aux mélanges sont très grosses. On peut cependant y parvenir en utilisant la théorie de la représentation du groupe symétrique. C'est une technique souvent utilisée pour étudier les mélanges et leurs valeurs propres. Les éléments importants de théorie de la représentation pour comprendre la thèse sont présentés au \autoref{chap:th_rep}. Au \autoref{chap:1efamille}, on présente les résultats pour la première famille, mais on délègue les preuves au \autoref{chap:preuves}. Pour démontrer les résultats sur la première famille, on utilise largement l'approche adoptée par Anton Dieker et Franco Saliola, qui ont trouvé toutes les valeurs propres des matrices associées au mélange doublement aléatoire \cite{DS}. Enfin, on présente au \autoref{chap:2efamille} la deuxième famille et ce qu'on sait sur ce mélange.\\
%
\end{introduction}
\chapter{Préliminaires}\label{chap:chap0}
\section{Premières notions}
\paragraph{Permutations}
Un objet d'études central dans cette thèse est l'ensemble des permutations. En effet, nous étudions des mélanges de collections finies d'objets. Ces mélanges ne sont rien d'autre qu'une succession de permutations, choisies en fonction d'une distribution de probabilités fixée.

Une \textit{permutation} est une bijection d'un ensemble vers lui-même. Plus intuitivement, c'est une réorganisation des objets dans l'espace, et la permutation correspond à la façon dont on déplace les objets.

Pour une permutation $\sigma$\index{$\sigma$} d'un ensemble de $n$ éléments, on utilise principalement trois notations~:
\begin{itemize}
	\item La notation sur deux lignes, dans laquelle on place les nombres de $1$ à $n$ sur la première ligne, alignés sur les nombres $\sigma(1), \sigma(2), \ldots, \sigma(n)$. Elle n'est pas très utilisée dans cette thèse.
	\item La notation sur une ligne est définie à partir de la notation sur deux lignes. On retire toutefois la première ligne, puisqu'il s'agit toujours des nombres de $1$ à $n$. À moins d'indication contraire, c'est la notation qui est utilisée dans ce document.
	\item La notation en cycles consiste à découper la permutation en cycles disjoints, puis à inscrire les cycles entre parenthèses. Notons que, comme ces cycles sont disjoints, on peut les écrire dans n'importe quel ordre. Les points fixes, qui sont des cycles de longueur $1$, sont parfois omis.
\end{itemize}

\begin{ex}
	La permutation
	$\left(\begin{smallmatrix}
	1 & 2 & 3 & 4 & 5 & 6 & 7\\
	1 & 4 & 7 & 5 & 2 & 6 & 3
	\end{smallmatrix}\right)$
	peut aussi s'écrire $1475263$ et $(37)(245)$.
\end{ex}

Le nombre de permutations de $n$ est noté $n!$\index{$n$ factorielle@$n"!$}, prononcé \textit{$n$ factorielle} et correspond au nombre \mbox{$n \cdot (n-1) \cdots 3 \cdot 2 \cdot 1$}.

L'ensemble des entiers de $1$ à $n$ est noté $[n] = \{1,2, \ldots, n\}$. \index{$[n]$}

\paragraph{Groupe symétrique}
Un ensemble muni d'une loi de composition qui satisfait certaines propriétés algébriques (associativité de la loi de composition, existence d'un élément neutre et existence d'un inverse pour chaque élément) est un groupe, et un groupe est abélien si, pour toute paire d'éléments $g$ et $h$ du groupe, $gh=hg$; on dit alors que $g$ et $h$ commutent. Les permutations de $n$ forment un des exemples les plus connus de groupe non-abélien. Ce groupe est appelé le \textit{groupe symétrique}.
\begin{ex}
	Les permutations ne commutent pas, comme on peut le voir avec les permutations de 3 éléments $(12) = \left(\begin{smallmatrix}
	1 & 2 & 3\\ 2 & 1 & 3
	\end{smallmatrix}\right)$ et $(23)  = \left(\begin{smallmatrix}
	1 & 2 & 3\\ 1 & 3 & 2
	\end{smallmatrix}\right)$.
	Ainsi,
	\[ (12) \circ (23) =  (123) \text{ et }  (23) \circ (12) = (132). \]
\end{ex}

\paragraph{Vocabulaire de la combinatoire des mots}
Les permutations sont très utiles pour exprimer le déplacement d'objets, mais aussi pour lister les objets d'une collection eux-mêmes lorsqu'ils sont tous distincts. Cependant, il arrive qu'on veuille mélanger des éléments qui ne sont pas forcément distincts. On cherche alors une certaine forme de permutation qui admettrait des répétitions.

Un outil qui permet de faire ça est les mots. Un \textit{mot} est une suite ordonnée de lettres prises dans un alphabet fixé (et fini). Nous ne nous intéressons qu'aux mots finis dans le cadre de cette thèse.

Un mot $u$ est un \textit{sous-mot} de $w$ si toutes les lettres de $u$ apparaissent aussi dans $w$, dans le même ordre et de façon pas forcément contiguë.

\begin{ex}
	Le mot de la langue française \textit{symétrique} contient comme sous-mot le mot \textit{yéti}.
\end{ex}

\paragraph{Opérations binaires sur les mots~: concaténation et battage}
Nous utiliserons deux opérations binaires sur les mots. La première est la concaténation de deux mots $u$ et $v$, notée $u\cdot v$, \index{$u\cdot v$} et son résultat est le mot formé de longueur $|u|+|v|$ dont $u$ est le préfixe et $v$, le suffixe.

\begin{ex}
	La concaténation des mots de la langue française \textit{pré} et \textit{fixe} donne le mot \textit{préfixe}.
	
	Le nom \textit{parcelle} est la concaténation du nom \textit{parc} et du pronom \textit{elle}.
\end{ex}

L'autre opération, un peu plus complexe, est celle de \textit{battage}. Le battage de deux mots $u$ et $v$, noté $u \shuffle v$,\index{Opérateurs!$\shuffle$} ressemble, intuitivement, au résultat obtenu en intercalant les lettres de $u$ et $v$, tout en gardant celles de $u$ (respectivement de $v$) dans le même ordre que dans le mot original. Le résultat est la combinaison linéaire de tous les résultats possibles, et les coefficients sont le nombre de façons d'obtenir chacun des résultats.
\begin{ex}
	L'adjectif \textit{versatile} peut être obtenu par battage des mots \textit{veste}  et \textit{rail}. On peut le voir à la \autoref{fig:versatile}, où les lettres de \textit{veste} sont en rouge et soulignées, alors que celle de \textit{rail} sont en bleu.
	\begin{figure}
		\begin{center}
		{\Large {\color{red}\underline{v} \underline{e}} {\color{blue} r}{\color{red}\underline{ s}} {\color{blue} a }{\color{red}\underline{t}} {\color{blue}  i l }{\color{red} \underline{e}} }
	\end{center}
	\caption{Le mot \textit{versatile} peut être vu comme un résultat du battage des mots \textit{veste} et \textit{rail}.}\label{fig:versatile}
	\end{figure}
	Ainsi,
	\[\text{veste} \shuffle \text{rail} = \ldots + \text{versatile} + \ldots\]
\end{ex}
\begin{ex}
	Parmi les battages des mots \textit{tain} et \textit{rente}, on retrouve le nom \textit{trentaine} et l'adjectif \textit{tarentine}. Ainsi,
	\[ \text{tain} \shuffle \text{rente} = \ldots + \text{tarentine} + \ldots + \text{trentaine} + \ldots \]
\end{ex}
\begin{ex}
	Le battage des suites de lettre \textit{bo} et \textit{om} donne la combinaison linéaire suivante, qui comprend le mot \textit{boom}~:
	\[\text{bo} \shuffle \text{om} = 2 \cdot \text{boom} + \text{bomo} +\text{obom} + \text{obmo} + \text{ombo}.\]
\end{ex}

Le battage de mots ne semble peut-être pas si intuitif. Cependant, il peut sembler naturel aux joueurs de cartes à jouer. Il s'agit en effet d'une étape cruciale de plusieurs mélanges habituels, notamment le mélange américain (connu en anglais sous le nom de \textit{riffle shuffle}), dans lequel on sépare le paquet en deux plus petits paquets, que l'on bat pour obtenir un nouveau paquet. Le tout est illustré à la \autoref{fig:melange_americain}.
\afterpage{
	\begin{figure}
		\centering
		\includegraphics[height=3cm]{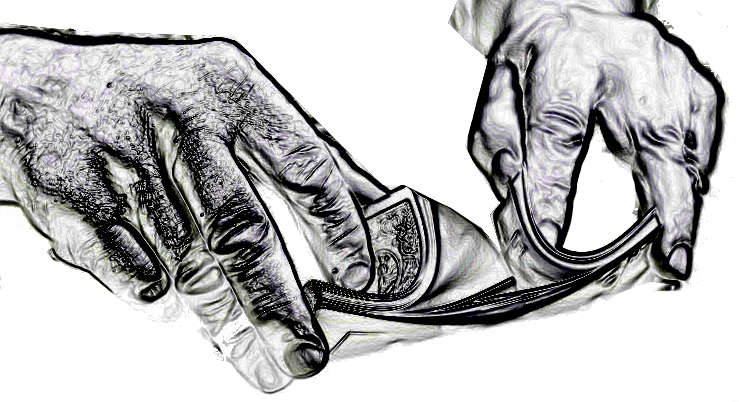}
		\caption[Battage de cartes dans le mélange américain.]{Battage de cartes dans le mélange américain.\footnotemark}\label{fig:melange_americain}
	\end{figure}
	\footnotetext{Source de l'image~: Todd Klassy, sur Wikimedia Commons, modifiée.}
}

\paragraph{Retrait d'un objet dans une liste}
On note le retrait d'un objet dans une liste par l'accent circonflexe lorsque le contexte est clair.\index{$\hat{i}$} Plus explicitement, \[ (1,2,\ldots, i-1, i+1, \ldots, n) = (1, 2, \ldots, \widehat{i}, \ldots, n). \]

\section{Opérateurs de mélange}

Cette thèse s'intéresse aux mélanges d'une collection d'objets. Il peut s'agir d'une grande diversité d'objets, à commencer par les paquets de cartes à jouer, qui sont souvent cités en exemple. Les cartes à jouer ont la propriété particulièrement intéressante que tous les objets soient distincts (dans un paquet de cartes qui n'est pas truqué, les cartes ne se répètent pas), mais qu'ils se ressemblent tous (les cartes ont toutes la même forme, la même taille, le même poids; et la chance qu'une carte soit pigée ne dépend d'aucune de ces caractéristiques). Toutefois, d'autres collections d'objets peuvent être étudiées, et les propriétés sont alors différentes. Les séquences d'ADN sont des collections avec des objets (les nucléotides) répétés, mais similaires \cite{durrett}. En revanche, des livres dans une bibliothèque \cite{phatarfod} et des fichiers dans un ordinateur \cite{donnelly} sont des collections d'objets généralement distincts, mais qui ont aussi des poids distincts. Nous ne nous intéressons pas au dernier cas, mais nous permettons, dans cette thèse, que des objets puissent être répétés.

La littérature sur les mélanges est abondante. En particulier, les mélanges de cartes ont attiré l'attention de nombreux chercheurs, qui ont en revanche trouvé beaucoup de façons de mélanger un paquet de cartes. Parmi les mélanges les plus étudiés, on retrouve~:
\begin{itemize}
	\item le \textit{mélange américain}, qui est l'un des plus utilisés pour les jeux de cartes. L'affirmation voulant que sept itérations du mélange soient nécessaires et suffisantes pour obtenir un paquet bien mélangé s'est même taillé une place dans la presse généraliste, par exemple dans le New York Times. Pour en savoir plus, voir \cite{AD, BD, Diaconis_dev_RS}.
	\item le \textit{mélange par-dessus la main} est aussi utilisé par les joueurs de carte, et particulièrement par les joueurs de cartes occasionnels. Robin Pemantle a confirmé une impression populaire, soit que ce mélange est particulièrement inefficace, puisqu'il faut mélanger un paquet de $n$ cartes un nombre de fois de l'ordre de $n^2 \log(n)$ avant d'avoir atteint un niveau satisfaisant de hasard \cite{pemantle, jonasson}.
	\item même le \textit{mélange spatial}, dans lequel on dispose les cartes sur une table et on les déplace approximativement a été modélisé et étudié  récemment \cite{DPal}.
	\item le \textit{mélange par transpositions} est parmi les premiers à avoir été étudié \cite{DiSh}. C'est aussi le premier mélange pour lequel l'étude a nécessité d'utiliser la théorie de la représentation, une approche reprise ici.
\end{itemize}
Afin de définir les opérateurs de mélange symétrisés, nous aurons besoin de certains mélanges moins habituels. Ceux-ci sont listés à la \autoref{ssec:bhr}. Parmi eux, les plus importants sont les trois mélanges décrits ici.
\paragraph{De l'aléatoire au dessus~: la bibliothèque de Tsetlin}
Un des mélanges les plus étudiés est celui de la bibliothèque de Tsetlin. On peut s'imaginer qu'on dispose d'une collection de livres ordonnée d'une certaine façon; contrairement aux cartes à jouer, les livres représentent généralement le cas où les objets sont sélectionnés selon des probabilités différentes. On s'imagine plutôt facilement que, lorsque vient le temps de sélectionner un livre dans une vraie bibliothèque, la popularité des différents livres n'est pas égale. Dans le mélange de la bibliothèque de Tsetlin, les objets sont numérotés de $1$ à $n$, et on sélectionne le livre $i$ avec probabilité $p_i$. Avant de pouvoir piger un autre livre, on doit rendre le livre $i$. On le remet alors sur \textit{le dessus} de la collection; on définit une telle position s'il n'y a pas d'endroit évident pour le dessus. Ce mélange a été étudié par de nombreuses personnes, dont plusieurs qui étaient motivées par les similitudes entre ce mélange et la gestion de la mémoire cache dans les ordinateurs. Pour en savoir plus, on peut notamment lire \cite{fill}, qui présente un historique élaboré de l'étude de ce mélange.

Un cas particulier de la bibliothèque de Tsetlin est lorsque tous les livres ont la même chance d'être sélectionnés. Ça peut paraître curieux pour une bibliothèque, et c'est pourquoi on privilégie alors l'exemple du mélange de cartes. Dans ce cas, on appelle le mélange celui \textit{de l'aléatoire au dessus}, ou \textit{random-to-top} en anglais. Pour plus d'informations sur ce mélange, particulièrement lorsqu'il est vu comme un opérateur de \bhr, voir l'\autoref{ex:tsetlin}.
\paragraph{L'opérateur du dessus à l'aléatoire}
Un mélange plutôt similaire à la bibliothèque de Tsetlin est celui où, plutôt que de retirer un objet et de le placer sur le dessus, on place l'objet du dessus n'importe où. Ces mélanges se ressemblent étant donné qu'ils correspondent à l'opération inverse, dans le sens suivant~: si on pouvait obtenir une collection d'objets placés dans l'ordre donné par la permutation $\sigma$ à partir de la permutation $\tau$ avec la bibliothèque de Tsetlin avec une certaine probabilité, la probabilité est la même d'obtenir une collection placée selon $\tau$ à partir de $\sigma$ avec le mélange du dessus à l'aléatoire. Nous introduisons les matrices associées aux opérateurs de mélange à la prochaine sous-section; une fois qu'on aura les outils nécessaires, on pourra tout simplement dire que les matrices des opérateurs de l'aléatoire au dessus et du dessus à l'aléatoire sont mutuellement transposées.

Le mélange du dessus à l'aléatoire (\textit{top-to-random}, en anglais) a par exemple été étudié par \cite{AD}, pour la version expliquée ici et dans le cas d'une collection d'objets tous de même poids. Une version plus générale a également été définie, dans laquelle on retire un nombre fixé d'objets sur le dessus de la collection, puis on les réinsère un à un, sans forcément préserver l'ordre \cite{DFP}.
\paragraph{Le mélange doublement aléatoire}\label{par:r2r_def}
Enfin, le mélange doublement aléatoire, ou \textit{random-to-random} en anglais, consiste à retirer n'importe quel objet de façon aléatoire dans notre collection, puis à le replacer n'importe où. C'est en quelque sorte la composition des mélanges de la bibliothèque de Tsetlin et de celui du dessus à l'aléatoire.

Imaginé par Persi Diaconis et Laurent Saloff-Coste \cite{DSC}, il a attiré beaucoup d'attention depuis \cite{UR, SCZ,Subag,QM,DS}, notamment parce que son temps de mélange, c'est-à-dire le nombre d'itérations requises pour bien mélanger les objets, s'est avéré très difficile à calculer. Il a finalement été calculé récemment par Megan Bernstein et Evita Nestoridi \cite{BN}. L'ingrédient-clef pour le calcul était une description combinatoire des valeurs propres de l'opérateur, trouvée par Anton Dieker et Franco Saliola \cite{DS}. Leur approche est décrite à la \autoref{ssec:r2r}.

Certains opérateurs de mélanges symétrisés correspondent à une généralisation du mélange doublement aléatoire, dans laquelle nous déplaçons un plus grand nombre d'objets; on discute de ces mélanges au \autoref{chap:1efamille}. On dira plus tard que le mélange doublement aléatoire est la version symétrisée de la bibliothèque de Tsetlin.

\subsection{Matrices, valeurs propres, vecteurs propres}

Les opérateurs de mélange sont des opérateurs linéaires, ce qui fait qu'on peut définir des matrices leur étant associées. Les lignes et les colonnes de ces matrices carrées sont indexées par les dispositions possibles de notre collection d'objets. En particulier, lorsque les objets sont tous distincts, ce sont les permutations qui indexent les lignes et les colonnes. L'entrée de la matrice de l'opérateur $\varphi$ à la position $(c_1, c_2)$ est le nombre de façons d'obtenir des objets disposés selon $c_2$ à partir de $c_1$ en utilisant une fois l'opérateur $\varphi$.

\begin{ex}
	Considérons une collection de trois objets distincts, \textsf{1}, \textsf{2} et \textsf{3}, de poids respectivement $p_1$, $p_2$ et $p_3$. Si l'on place le dessus de la liste à la fin du mot, on obtient la matrice du mélange de la bibliothèque de Tsetlin~:
	\[\bordermatrix{
		&\textsf{123}&\textsf{132}&\textsf{213}&\textsf{231}&\textsf{312}&\textsf{321}\cr
		\textsf{123}& p_3 & p_2 & 0 & p_1 & 0 & 0 \cr
		\textsf{132}& p_3 & p_2 & 0 & 0 & 0 & p_1 \cr
		\textsf{213}& 0 & p_2 & p_3 & p_1 & 0 & 0 \cr
		\textsf{231}& 0 & 0 & p_3 & p_1 & p_2 & 0 \cr
		\textsf{312}& p_3 & 0 & 0 & 0 & p_2 & p_1 \cr
		\textsf{321}& 0 & 0 & p_3 & 0 & p_2 & p_1 \cr
	}.\]
\end{ex}

\begin{ex}Il y a $6$ façons de disposer une collection contenant seulement deux paires d'objets identiques. On peut les représenter par des mots sur un alphabet binaire : chaque mot contient deux fois chacune des lettres. Alors, la matrice du mélange doublement aléatoire est
	\[ \bordermatrix{
		&{\scriptstyle \textsf{aabb}} & {\scriptstyle \textsf{abab}} & {\scriptstyle \textsf{baab}} & {\scriptstyle \textsf{abba}} &
		{\scriptstyle \textsf{baba}} & {\scriptstyle \textsf{bbaa}}\cr
		{\scriptstyle \textsf{aabb}} & 8 & 4 & 2 & 2 & 0 & 0 \cr
		{\scriptstyle \textsf{abab}} & 4 & 4 & 3 & 3 & 2 & 0 \cr
		{\scriptstyle \textsf{baab}} & 2 & 3 & 6 & 0 & 3 & 2 \cr
		{\scriptstyle \textsf{abba}} & 2 & 3 & 0 & 6 & 3 & 2 \cr
		{\scriptstyle \textsf{baba}} & 0 & 2 & 3 & 3 & 4 & 4 \cr
		{\scriptstyle \textsf{bbaa}} & 0 & 0 & 2 & 2 & 4 & 8
	}.\]
\label{ex:mat_R2R22}
\end{ex}

Les \textit{valeurs propres d'un opérateur $\varphi$} sont les valeurs propres de la matrice associée, c'est-à-dire les nombres $c$ pour lesquels il existe au moins un vecteur non-nul $\vec{v}$ satisfaisant $\varphi(\vec{v}) = c \cdot \vec{v}$. Un tel vecteur est appelé un \textit{vecteur propre}.

Il y a plusieurs raisons de s'intéresser aux valeurs propres des opérateurs de mélange. Celle qui est la plus souvent citée est que connaître les valeurs propres permet parfois de calculer le temps de mélange, c'est-à-dire le nombre d'itérations nécessaires pour que la collection soit bien mélangée. On en parle davantage à la \autoref{ssec:consequences}.

Comme pour toute matrice diagonalisable, connaître les valeurs propres et les vecteurs propres permet de diagonaliser la matrice, et donc de calculer très rapidement ses puissances. Ici, la puissance $m$-ième de la matrice donne de l'information sur le nombre de façons de passer d'un état à un autre de notre collection après $m$ itérations du mélange. Une interprétation des puissances de la matrice en termes de probabilités est aussi présentée à la \autoref{ssec:chaines_Markov}. On verra au même endroit une interprétation pour au moins un des vecteurs propres. 

\subsection{Opérateurs de Bidigare--Hanlon--Rockmore ($\bhr$)}\label{ssec:bhr}
Une façon de généraliser la bibliothèque de Tsetlin est de prendre un sous-ensemble de $k$ objets dans la suite $w_1\ldots w_n$ puis de les placer sur le dessus de la collection, tout en préservant leur ordre. Dans la littérature, on retrouve ce mélange sous le nom de \textit{$k$ déplacés vers le dessus} (\og$k$-pop shuffle\fg\ ou \og $k$-move-to-front \fg ). Ainsi, si $v_1 \ldots v_k$ est un sous-mot de $w_1 \ldots w_n$ et $u_1 \ldots u_{n-k}$ est son complémentaire, alors $v_1 \ldots v_k u_1 \ldots u_{n-k}$ serait un résultat possible du mélange $k$ déplacés vers l'avant.

\begin{ex}
	Dans le mot \textit{voitures}, \textit{vite} est un sous-mot et \textit{ours} est son complémentaire. Un des résultats possibles du mélange $4$ déplacés vers le dessus est donc \textit{oursvite}. Ici, le dessus est la fin du mot.
	\begin{center}
		\begin{tikzpicture}
		\node (v) at (0,0) {v};
		\node[red] (o) at (0.5,0) {\underline{o}};
		\node (i) at (1,0) {i};
		\node (t) at (1.5,0) {t};
		\node[red] (u) at (2,0) {\underline{u}};
		\node[red] (r) at (2.5,0) {\underline{r}};
		\node (e) at (3,0) {e};
		\node[red] (s) at (3.5,0) {\underline{s}};
		\draw[->, >=latex] (v) to [out=35,in=150] (3.8,0);
		\draw[->, >=latex] (i) to [out=40,in=150] (4,0);
		\draw[->, >=latex] (t) to [out=45,in=150] (4.2,0);
		\draw[->, >=latex] (e) to [out=50,in=150] (4.4,0);
		\end{tikzpicture}
		$\longrightarrow$
		\begin{tikzpicture}
		\node[red] (o) at (0,0) {\underline{o}};
		\node[red] (u) at (0.5,0) {\underline{u}};
		\node[red] (r) at (1,0) {\underline{r}};
		\node[red] (s) at (1.5,0) {\underline{s}};
		\node (v) at (2,0) {v};
		\node (i) at (2.5,0) {i};
		\node (t) at (3,0) {t};
		\node (e) at (3.5,0) {e};
		\end{tikzpicture}
	\end{center}
\end{ex}

De façon encore plus générale, on pourrait partitionner notre suite $w_1\ldots w_n$ en plusieurs sous-mots (disjoints), puis concaténer chacun de ces sous-mots, après avoir déterminé un ordre pour ceux-ci. Une façon de décrire l'ordre dans lequel on concatène les sous-mots est d'utiliser les compositions ensemblistes (expliquées plus loin dans cette section).

C'est d'ailleurs la notation qu'ont utilisé en 1999 Pat Bidigare, Phil Hanlon et Dan Rockmore \cite{BHR}. Ils ont étudié les mélanges \textit{par éclatement} (pop shuffles), dans lesquels les mots sont partitionnés en un certain nombre de sous-mots de longueur fixe, puis concaténés. 

\paragraph{Compositions ensemblistes}
\begin{defn}
	Une \textit{composition ensembliste} d'un ensemble $S$ est une liste (ordonnée) de sous-ensembles non-vides et disjoints dont l'union est $S$.
	
	Autrement dit, une composition de $S$ en est une partition ordonnée.
\end{defn}

\begin{ex}
	$[\{1,3,5\},\{2,4\}]$\index{$[\{\cdot\},\{\cdot\}]$} est une composition de l'ensemble $[5]$ et $[\{2,4\},\{1,3,5\}]$ en est une autre.
\end{ex}

\begin{ex}\label{ex:comp_en_bij_avec_perm}
	$[\{5\},\{3\},\{4\},\{1\},\{2\}]$ est une composition de l'ensemble $[5]$ dans laquelle tous les ensembles sont des singletons. De telles compositions sont en bijection avec les permutations. Ainsi, la composition $[\{5\},\{3\},\{4\},\{1\},\{2\}]$ peut être associée à la permutation $\mathsf{53412}$.
\end{ex}

\begin{ex}
	$[\{1,2\},\{2,3\}]$ n'est pas une composition ensembliste, parce que les sous-ensembles ne sont pas disjoints.
\end{ex}

La \textit{forme d'une composition ensembliste de $S$} est la composition de l'entier $|S|$ dans laquelle on ne retient que la taille des blocs. Autrement dit, la forme de $c$ est la liste (ordonnée) qui contient la taille du premier bloc de $c$, celle du deuxième bloc, et ainsi de suite.

\begin{ex}
	$[\{1,3,5\},\{2,4\}]$ est une composition de $[5]$ de forme $(3,2)$, alors que $[\{2,4\},\{1,3,5\}]$ est de forme $(2,3)$.
\end{ex}


\paragraph{Produit de compositions ensemblistes}
Le produit de deux compositions ensemblistes $b = [b_1, \ldots, b_j]$ et $c = [c_1, \ldots, c_l]$ est défini en raffinant les compositions. Par exemple, si on fait le produit $b\ast c$\index{$\ast$}, $b$ raffine les blocs de $c$ et on obtient l'expression
\[[b_1, \ldots, b_j] \ast [c_1, \ldots, c_l]  = [c_1 \cap b_1, c_1 \cap b_2, \ldots, c_1 \cap b_j,\ldots, c_l\cap b_1, \ldots, c_l \cap b_j], \]
dans laquelle on retire les ensembles vides de l'expression à droite.

\begin{ex}\label{ex:prod_comp}
	Avec des compositions de $[5]$,
	\[[\{1,3\},\{2,4,5\}] \ast [\{2,3,4\},\{5\},\{1\}] = [\{3\},\{2,4\},\{5\},\{1\}],\]
	\[ [\{2,3,4\},\{5\},\{1\}]\ast [\{1,3\},\{2,4,5\}] = [\{3\},\{1\},\{2,4\},\{5\}]\]
	et 
	\[[\{1,3\},\{2,4,5\}] \ast [\{12345\}] = [\{1,3\},\{2,4,5\}]. \]
\end{ex}

Il est facile de voir que ce produit n'est pas commutatif, puisque les blocs ne se retrouveraient pas forcément dans le même ordre dans la composition résultante. Un contre-exemple de la commutativité est donné à l'\autoref{ex:prod_comp}. En revanche, ce produit est associatif.

\begin{rmq}
	Dans plusieurs articles, notamment \cite{BHR}, où les opérateurs de mélange sont définis, le produit de composition s'effectue dans l'autre sens, c'est-à-dire que, si $b\ast c$, alors $c$ raffine les blocs de $b$. Or, comme nous travaillons avec des vecteurs ligne multipliés à droite par des matrices (tel qu'expliqué à la \autoref{rmq:mult_mat_a_droite}), nous avons choisi de modifier la définition pour qu'elle soit cohérente avec le but que nous poursuivons.
\end{rmq}
\paragraph{Monoïde de compositions ensemblistes}
On peut munir l'ensemble des compositions ensemblistes de $[n]$ d'une loi de composition pour le doter d'une structure algébrique. Cette loi est le produit de compositions ensemblistes, qui est associatif. Le \textit{monoïde de compositions} est obtenu en notant que la composition formée d'un seul sous-ensemble agit comme l'élément neutre, comme à l'\autoref{ex:prod_comp}.


%

\subsubsection{Mélanges et compositions}
Bidigare, Hanlon et Rockmore, dans leur étude des mélanges par éclatement, ont étudié le produit de composition sur les permutations \cite{BHR}. En fait, les permutations de $[n]$ sont en bijection avec les compositions ensemblistes de $[n]$ composées uniquement de singletons~: pour passer de la composition à la permutation, il suffit de noter dans l'ordre les éléments dans chacun des singletons, comme à l'\autoref{ex:comp_en_bij_avec_perm}.

De plus, si on multiplie une composition ensembliste composée uniquement de singletons avec une autre composition ensembliste, on obtient de nouveau une liste de singletons. En utilisant la bijection, on peut donc multiplier une permutation avec une composition ensembliste pour obtenir de nouveau une permutation.

Les mélanges par éclatement sont définis sur les permutations comme la multiplication (à droite) par une combinaison linéaire de compositions ensemblistes. Chaque mélange est défini par la distribution de probabilités sur les combinaisons linéaires.

\begin{ex} On peut multiplier \textsf{13254} par une combinaison linéaire de compositions de $[5]$~:
	\[\mathsf{13254} \ast \left( [\{2,3\}, \{1\}, \{4,5\} ] + 2 \cdot  [\{5\},\{1,2,3\}, \{4\} ]\right) = \mathsf{32154} + 2 \cdot \mathsf{51324}.\]
\end{ex}

Pour obtenir les mélanges par éclatement, on choisit avec probabilité $p_c$ la composition $c$, puis on agit avec $c$ sur les permutations. La probabilité de transition de $\sigma$ à $\tau$ est $\sum_{c \mid \sigma \ast c = \tau} p_c$.
\begin{ex}[La bibliothèque de Tsetlin]\label{ex:tsetlin}
	On a défini la bibliothèque de Tsetlin comme le mélange qui déplace le $i$-ième objet vers le dessus de notre collection. On fixe ici que le dessus soit à la droite des mots. Si nous admettons que les objets aient des poids différents, nous définissions $p_i = p_{[[n]\backslash\{i\}, \{i\}]}$ comme la probabilité de déplacer le $i$-ième livre.
	\begin{center}
		\begin{tikzpicture}
		\node[inner sep=0pt] (b1) at (0,0)
		{\includegraphics[height=2cm]{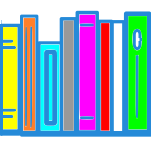}};
		\node[inner sep=0pt] (b2) at (6,0)
		{\includegraphics[height=2cm]{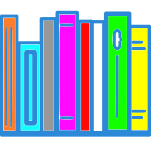}};
		\draw[->,thick] (2,0) -- (4,0) node[above, midway] {$p_1 = p_{[\{2,\ldots, 8\}, \{1\}]}$};
		\end{tikzpicture}
	\end{center}	
\end{ex}
La bibliothèque de Tsetlin correspond ainsi à multiplier à droite par $\sum_{i=1}^n\ p_i [[n]\backslash\{i\}, \{i\}]$.


\subsubsection{Valeurs propres et vecteurs propres des opérateurs de $\bhr$}
Un des grands résultats obtenus par Bidigare, Hanlon et Rockmore est une façon combinatoire de calculer les valeurs propres de tous les mélanges par éclatement \cite[théorème 2.1]{BHR}. Celles-ci sont indexées par les partitions ensemblistes (les listes non-ordonnées de sous-ensembles non-vides et disjoints dont l'union est l'ensemble de départ) et peuvent être calculées à l'aide d'outils combinatoires simples.

\begin{rmq}
	Les opérateurs de \bhr\ sont définis en fonction de la probabilité de choisir une certaine composition ensembliste, et cette valeur est inscrite comme le coefficient devant la composition, dans la combinaison linéaire. Dans le cadre de cette thèse, la probabilité de choisir une certaine combinaison linéaire a toujours été fixée soit à $0$, soit elle correspond à $1$ divisé par le nombre de combinaisons linéaires qu'on pourrait choisir. Pour la suite de la thèse, on utilise donc $0$ ou $1$ comme coefficients pour les combinaisons linéaires.
\end{rmq}

\subsection{Opérateurs de mélange symétrisés}\label{ssec:op_sym}
Les opérateurs de mélange symétrisés sont construits à partir des opérateurs de \bhr. Pour un opérateur de \bhr\ $\pi$, sa \textit{version symétrisée} s'écrit $\pi^\top\circ \pi$. Nous nous intéressons ici à deux familles d'opérateurs symétrisés~:
la première, $\allnuk$, est présentée au \autoref{chap:1efamille}, et est telle que $\nu_k$ est un multiple de $\pi_k^\top\circ \pi_k$ avec $\pi_k$ la somme des compositions de la forme $(n-k, 1, \ldots, 1)$.
Quant à la deuxième famille, $\allgammak$, elle est telle que $\gamma_k$ est un multiple de ${\pi'}_k^\top\circ {\pi'}_k $ lorsque ${\pi'}_k$ est la multiplication par la somme des compositions formées de $k$ blocs de taille $2$ suivis de singletons. Cette famille est étudiée au \autoref{chap:2efamille}.

\begin{ex}
	Le mélange doublement aléatoire, présenté à la \autopageref{par:r2r_def}, est la version symétrisée de la bibliothèque de Tsetlin. En effet, il consiste à appliquer successivement les mélanges de l'aléatoire au dessus et du dessus à l'aléatoire.
\end{ex}

Bien qu'on connaisse les valeurs propres des opérateurs de \bhr, cela est peu utile pour trouver les valeurs propres des opérateurs symétrisés. Malheureusement, très peu de propriétés d'une certaine matrice sont forcément vraies pour sa version symétrisée.

\subsubsection{Suites croissantes}\label{sssec:non-inversions}
Un outil essentiel pour notre étude des opérateurs de mélange symétrisés est la notion de suites croissantes d'une permutation. En effet, on verra qu'une façon de calculer les entrées des matrices est à l'aide des suites croissantes de permutations (pour plus d'informations sur la construction des matrices, voir la \autoref{ssec:def_matrices_nuk} et la \autoref{ssec:gammak_matrices}).

\begin{defn}
	Dans une permutation $\sigma$, une \textit{suite croissante de longueur $k$} est un ensemble $\{i_1, \ldots, i_k \}$, où à la fois les suites $(i_1, i_2, \ldots, i_k)$ et $(\sigma^{-1}(i_1), \ldots,\sigma^{-1}(i_k))$ sont placées en ordre croissant.\\
\end{defn}
\begin{ex}
	Considérons la permutation $\sigma = \mathsf{5623147}$. Sa plus longue suite croissante est de longueur $4$ et est unique~: il s'agit de $\mathsf{2347}$. En plus des sous-suites de cette dernière, $\sigma$ contient comme sous-suites croissantes de longueur $3$ les suites $\mathsf{147}$ et $\mathsf{567}$.
\end{ex}

\subsection{Chaînes de Markov}\label{ssec:chaines_Markov}
\raggedbottom
Un autre point de vue sur les opérateurs de mélange, très répandu, est celui des chaînes de Markov. Une \textit{chaîne de Markov} est un ensemble d'états (ici, les dispositions de notre collection d'objets) muni de transitions entre les états, qui s'effectuent avec une probabilité donnée. Ces probabilités de transition d'un état à un autre restent les mêmes tout au long du processus et ne dépendent pas de ce qui s'est passé auparavant.

Les chaînes de Markov sont souvent représentées par leur matrice de transition. Comme la probabilité d'effectuer une transition reste constante tout au long de l'exécution du processus, la seule information dont on ait besoin pour décrire complètement la chaîne de Markov est la matrice qui donne les probabilités de transition. Cette matrice est indexée par les états, et l'entrée $(i,j)$ de la matrice représente la probabilité de passer de l'état $i$ à l'état $j$.
\vspace{1em}

\begin{minipage}{0.82\linewidth}
	\begin{ex}
		Un premier exemple de chaîne de Markov est souvent celui-ci. Quelqu'un lance successivement une pièce de monnaie biaisée qui a une probabilité $p$ de tomber sur le côté pile et une probabilité $f = 1-p$ de tomber sur le côté face. Il y a deux états~: l'un correspond à avoir, jusqu'à un certain moment, obtenu un nombre pair de fois \textit{pile}, c'est l'état $\mathbf{0} $. L'autre état, $\mathbf{1}$, est atteint lorsqu'on a obtenu un nombre impair de fois \textit{pile}.
	\end{ex}
\end{minipage}
\begin{minipage}{0.14\linewidth}
	\centering
	\includegraphics[height=3.7cm]{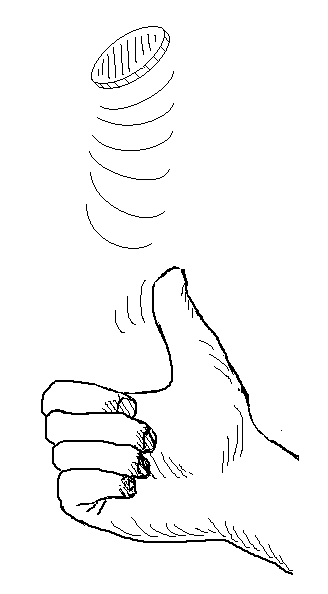}
\end{minipage}

On peut expliciter la chaîne de Markov par la matrice de transition qui suit. Rappelons que, dans cette situation, on change d'état exactement lorsque le joueur obtient le résultat pile~: on change d'état avec probabilité $p$, et on reste dans le même avec probabilité $f$.
	\[ \bordermatrix{
	&{\scriptstyle \mathbf{0}} & {\scriptstyle \mathbf{1}}\cr
	{\scriptstyle \mathbf{0}} & f & p \cr
	{\scriptstyle \mathbf{1}} & p & f  \\
}.\]

	Les deux états pourraient correspondre à deux emplacements, et le lanceur de pièce se déplace d'un emplacement à l'autre seulement s'il obtient pile. Ce joueur est souvent personnalisé comme une grenouille qui se promène d'un nénuphar à l'autre; ici, les nénuphars représentent les états.

À partir de la définition, on peut déduire certaines propriétés des matrices de transition~:
\begin{prop}
	Les matrices de transition des chaînes de Markov n'ont que des entrées positives ou nulles et la somme des entrées d'une même ligne vaut $1$~: c'est la somme des probabilités de transition d'un état donné vers tous les états, et ça correspond donc à tous les événements possibles. De telles matrices sont appelées des \textit{matrices stochastiques}.
\end{prop}

Les opérateurs de mélange ont tous des entrées positives. Cependant, les matrices des opérateurs de mélange ne contiennent que des entiers, et la somme des entrées d'une même ligne est en général beaucoup plus grande que $1$.

En revanche, la somme des entrées d'une même ligne correspond au nombre de suites de mouvements qu'on peut faire. Ceci est indépendant de la configuration de notre collection d'objets, et donc la somme d'une ligne est toujours la même. Pour obtenir une chaîne de Markov, on peut diviser la matrice par la somme des entrées d'une ligne.

\begin{ex}
	À l'\autoref{ex:mat_R2R22}, on a vu la matrice de l'opérateur de mélange doublement aléatoire sur les collections d'objets qui contiennent deux paires d'objets identiques. La somme de chaque ligne est $16$. Ainsi, la matrice de transition associée à ce mélange est
	\[ \frac{1}{16} \times \bordermatrix{
		&{\scriptstyle \textsf{aabb}} & {\scriptstyle \textsf{abab}} & {\scriptstyle \textsf{baab}} & {\scriptstyle \textsf{abba}} &
		{\scriptstyle \textsf{baba}} & {\scriptstyle \textsf{bbaa}}\cr
		{\scriptstyle \textsf{aabb}} & 8 & 4 & 2 & 2 & 0 & 0 \cr
		{\scriptstyle \textsf{abab}} & 4 & 4 & 3 & 3 & 2 & 0 \cr
		{\scriptstyle \textsf{baab}} & 2 & 3 & 6 & 0 & 3 & 2 \cr
		{\scriptstyle \textsf{abba}} & 2 & 3 & 0 & 6 & 3 & 2 \cr
		{\scriptstyle \textsf{baba}} & 0 & 2 & 3 & 3 & 4 & 4 \cr
		{\scriptstyle \textsf{bbaa}} & 0 & 0 & 2 & 2 & 4 & 8\\
	}.\]

\end{ex}

Une chaîne de Markov est un processus qui est en général répété plusieurs fois. On souhaite donc connaître les puissances de la matrice. Les entrées de la puissance $m$-ième d'une matrice de transition correspondent en effet aux probabilités de transition en $m$ étapes.
\begin{rmq}\label{rmq:mult_mat_a_droite}
	Une convention largement répandue dans la littérature sur les chaînes de Markov correspond à ce qu'on fait intuitivement~: on part d'une certaine distribution, qui représente les probabilités de se trouver dans un état donné au temps $t$, puis on multiplie à droite par la matrice de transition. On obtient ainsi une nouvelle distribution, qui symbolise les probabilités de se retrouver dans chacun des états au temps $t+1$. Cette convention, qui s'oppose à celle souvent retenue pour les opérateurs (multiplier une matrice par un vecteur colonne pour obtenir un vecteur colonne), est toutefois celle retenue pour cette thèse.
\end{rmq}

\begin{figure}
	\[ \underbrace{\begin{pmatrix} q_{1,t} & \ldots & q_{n,t}	\end{pmatrix}}_{\text{Distribution initiale}} 
	\cdot \underbrace{\begin{pmatrix} p_{1,1} & \ldots& p_{1,n}\\
		p_{2,1} & \ldots & p_{2,n}\\
		\vdots & \ddots & \vdots\\
		p_{n,1} & \ldots &p_{n,n}
	\end{pmatrix} }_{\text{Matrice de transition}}  =
	\underbrace{\begin{pmatrix} q_{1,t+1} & \ldots &q_{n,t+1}	\end{pmatrix}}_{\text{Distribution finale}}. \]
	\caption{Multiplication à droite par les matrices.}
\end{figure}

\paragraph{Valeurs propres des chaînes de Markov}
Plusieurs choses sont connues sur les valeurs propres des chaînes de Markov. C'est donc de l'information que nous avons \textit{a priori} sur les valeurs propres. Parmi celles-ci, le théorème de Perron-Frobenius est particulièrement utile.
\begin{prop}[Théorème de Perron-Frobenius, reformulé]\label{prop:perron-frobenius}
	Les matrices de transition des chaînes de Markov ont comme valeurs propres des nombres complexes dont le module est inférieur ou égal à $1$, et $1$ est une valeur propre de toutes les matrices de transition. Si de plus la chaîne est irréductible, c'est-à-dire qu'à partir de n'importe quel état on peut rejoindre tous les autres états en un nombre fini d'étapes, et apériodique (une condition généralement satisfaite par les opérateurs de mélanges), alors,
	\begin{itemize}
		\item $1$ est la seule valeur propre dont le module est $1$
		\item le vecteur propre (à gauche) associé à la valeur propre $1$ est l'unique distribution stationnaire, soit la distribution vers laquelle la chaîne de Markov converge~: $\vec{v} \cdot T = \vec{v}$.\footnote{Pour les mélanges d'objets de même poids, la distribution stationnaire est généralement la distribution uniforme, mais ce n'est pas le cas des mélanges où les objets ont des poids distincts.}
	\end{itemize}
\end{prop}

Les relations entre propriétés algébriques des matrices de transition et des questions sur les probabilités sont exhibées au \autoref{tab:alg_vs_proba}.
\begin{table}
	\caption{Correspondances entre propriétés algébriques et questions de probabilités.}\label{tab:alg_vs_proba}
	\centering
	\begin{tabular}{|ccc|}
		\hline
		Questions sur les processus aléatoires & & Propriétés algébriques\\
		\hline
		les probabilités après $m$ étapes? & $\longleftrightarrow$ & les entrées de $T^m$ \\
		\hdashline
		\begin{tabular}{c}le comportement à long terme? \\(la distribution stationnaire) \end{tabular}& $\longleftrightarrow$& \begin{tabular}{c}les vecteurs propres $\vec v$\\
		tels que $\vec v T = \vec v$\end{tabular}\\
		\hdashline
		\begin{tabular}{c} le taux de convergence à\\ la distribution stationnaire? \end{tabular} & $\longleftrightarrow$ & \begin{tabular}{c} contrôlé par les valeurs\\		 propres de $T$\end{tabular}\\
		\hline
	\end{tabular}
\end{table}

\chapter{Théorie de la représentation du groupe symétrique}\label{chap:th_rep}

Un outil fondamental à la recherche des valeurs propres d'opérateurs est ce qu'on appelle la théorie de la représentation. Cette théorie permet en effet de réduire la tâche titanesque de calculer les valeurs propres sur d'immenses espaces vectoriels à celle, plus réaliste, de réaliser des calculs sur des sous-espaces de taille (plus) raisonnable. Nous nous servons aussi de la théorie de la représentation pour procéder par induction des mélanges de $n$ objets à ceux de $n+1$ objets. Le point de vue retenu pour présenter la théorie de la représentation est celui des modules. La démarche adoptée est la suivante~:
\begin{itemize}
	\item On décompose l'espace vectoriel sur lequel on agit (l'algèbre du groupe symétrique) en sous-espaces stables pour une action du groupe symétrique. Ceux-ci sont appelés des \textit{modules} (\autoref{sec:modules}).
	\item On pousse la décomposition jusqu'à ce qu'on obtienne des \textit{modules simples}, c'est-à-dire des modules qui ne contiennent pas de sous-modules (propres). On verra que cela est possible, comme l'indique le \textit{théorème de Maschke} (\autoref{ssec:decomposition_modules_simples}).
	\item Heureusement, la décomposition de l'algèbre du groupe symétrique en modules simples est bien connue. On présente la \textit{règle de Young}, qui nous donne le détail de cette décomposition. Ce théorème donne une correspondance entre les modules simples et des objets combinatoires, qu'on appelle des tableaux standards (\autoref{sec:tableaux}).
	\item On utilise le \textit{lemme de Schur}, qui nous permet d'associer les valeurs propres des opérateurs de mélange aux modules simples (\autoref{thm:lemme_Schur}).
	\item Pour procéder par induction entre les espaces correspondant à des collections de $n$ et $n+1$ objets, on a besoin de la \textit{règle de branchement} pour passer d'un module simple à un autre (\autoref{ssec:regle_branchement}).
	\item Enfin, on présente les \textit{caractères du groupe symétrique}. Les caractères d'un groupe sont un outil particulièrement efficace pour une diversité de calculs. Bien qu'ils ne soient que sporadiquement utilisés dans cette thèse, certaines propriétés des caractères essentielles au projet sont présentées à la \autoref{sec:caracteres}.
\end{itemize}

\section{Actions du groupe symétrique}
\begin{defn}
	Une \textit{action du groupe symétrique} $S_n$ sur un ensemble $E$ de taille $n$ est une fonction
	\[ \ldotp : S_n  \times E \to E \]
	qui permute les éléments de $E$ et qui respecte les contraintes suivantes, pour toutes les permutations $\sigma$ et $\tau$ et pour tous les éléments $x$ et $y$ de l'ensemble $E$~:
	\begin{itemize}
		\item $\Id \ldotp x = x $,
		\item $\sigma \ldotp (\tau \ldotp x) = (\sigma \circ \tau) \ldotp x$,
		\item $\sigma \ldotp (kx+ly) = k (\sigma \ldotp x) + l (\sigma \ldotp y)$, pour tous scalaires $k$ et $l$.
	\end{itemize}
\end{defn}

\section{Modules}\label{sec:modules}
Un  espace vectoriel $V$ est un \textit{$S_n$-module} s'il existe une action de $S_n$ qui préserve $V$, c'est-à-dire telle que $\sigma \ldotp \vec{v} \in V$ pour toute permutation $\sigma$ et pour tout vecteur $\vec{v}$. Pour un groupe $G$ en général, on peut définir un $G$-module de la même façon.
\begin{ex}\label{ex:modules_Sn}
	Soit $V$ un espace vectoriel de dimension $n$, et considérons l'action de $S_n$ par permutation des coordonnées. Sous une telle action, le vecteur \mbox{$\vec{v} = (v_1, \ldots, v_n)$} est envoyé sur $\sigma \ldotp \vec{v} = (v_{\sigma(1)}, \ldots, v_{\sigma(n)})$.
	
	Les deux sous-espaces vectoriels suivants sont des $S_n$-modules~:
	\begin{itemize}
		\item le module trivial, qui est de dimension $1$,
		\[\{\vec{v} \mid v_1 = v_2 = \ldots = v_n\}.\]
		\item le module standard, de dimension $n-1$,
		\[\{\vec{v} \mid v_1 + v_2 + \ldots + v_n = 0\}.\]
	\end{itemize}
\end{ex}

\subsubsection{Morphismes de modules}
\begin{defn}
	Une application linéaire qui préserve l'action d'un groupe est appelée un \textit{homomorphisme de modules} (ou tout simplement \textit{un morphisme de modules}). Une telle application $\theta$ satisfait $\theta(g\ldotp v) = g \ldotp \theta(v)$ pour tout élément $g$ du groupe et pour tout vecteur $v$ de l'espace sur lequel $g$ agit.
\end{defn}
\begin{ex}
	Supposons que le groupe est $S_3$ et que l'espace sur lequel $S_3$ agit est $\mathbb{R}_3$, par permutation des coordonnées. Alors, une homothétie de rapport $k$ est un homomorphisme, car \mbox{$k\cdot (v_{\sigma(1)},v_{\sigma(2)},v_{\sigma(3)}) = ((kv)_{\sigma(1)},(kv)_{\sigma(2)},(kv)_{\sigma(3)})$}.
	En revanche, l'application qui multiplie uniquement la première coordonnée par le scalaire $k$ n'est pas un homomorphisme, car, pour la permutation $(12)$, \mbox{$(k\cdot v_{2},v_{1},v_{3}) \neq (v_{2},k\cdot v_1,v_3)$}.
\end{ex}

\subsection{Décomposition en modules simples}\label{ssec:decomposition_modules_simples}
Comme mentionné dans le préambule du chapitre, nous nous intéressons aux modules parce qu'ils nous permettent de limiter les calculs de valeurs propres à de beaucoup plus petits espaces vectoriels. Deux résultats sont fondamentaux dans l'étude de la décomposition des modules~: le premier, dû à Heinrich Maschke, est originalement paru en 1899 \cite{Maschke} \footnote{La théorie de la représentation était naissante à l'époque de Maschke, et l'article original de Maschke ne parle ni de représentations, ni de modules, mais plutôt d'actions linéaires de groupes.}, alors que le second est dû à Issaï Schur et est paru en 1907 \cite{Schur}.  Les théorèmes cités ici sont une adaptation dans le langage moderne de la théorie de la représentation, et peuvent être retrouvés dans divers ouvrages : \cite[théorème 15.1 et lemme 15.11]{DF}, \cite[théorèmes 1.5 et 1.7]{FH}, \cite[théorèmes 1.5.3 et 1.6.5]{sagan}, entre autres. La version donnée ici s'applique aux algèbres de groupe dont le corps de base est algébriquement clos et de caractéristique $0$, les nombres complexes par exemple.

\begin{thm}[Théorème de Maschke]
	Soit $G$ un groupe fini et soit $V$ un $G$-module non-nul. Alors, $V$ se décompose en une somme directe de sous-modules simples~:
	\[V = U_1 \oplus U_2 \oplus \ldots \oplus U_k. \]
\end{thm}

\begin{thm}[Lemme de Schur]\label{thm:lemme_Schur}
	Soit $G$ un groupe et soit $U$ et $V$ deux $G$-modules simples. Un homomorphisme $\phi : U \to V$ est soit nul, soit un isomorphisme. De plus, si $U = V$, $\phi$ est une homothétie, c'est-à-dire qu'il existe un nombre $c \in \C$ tel que $\phi = c  \cdot \Id$.
\end{thm}

La conséquence du lemme de Schur qui nous intéresse le plus est que, pour chaque homomorphisme d'un module simple vers lui-même, il existe une unique valeur propre. La raison pour laquelle ceci est particulièrement utile est que les opérateurs de mélange sont des actions de groupe, et donc des homomorphismes de modules. Pour plus de précisions, voir la \autoref{ssec_contenant_par:melanges_alg_du_groupe}. 
Le théorème de Maschke nous indique d'ailleurs qu'il est possible de décomposer pleinement les modules sur lesquels nous agissons, jusqu'à l'obtention d'une somme directe de modules simples. Le détail de la décomposition des modules sur l'algèbre du groupe symétrique (notre principal sujet d'intérêt) est donné à la \autoref{sec:CSn}.

\section{Algèbre du groupe symétrique}\label{sec:CSn}
Soit $G$ un groupe. L'algèbre du groupe (sur le corps des nombres complexes), notée $\C G$,\index{Modules!$\mathbb{C}G$} est formée des combinaisons linéaires formelles d'éléments de $G$ à coefficients parmi les nombres complexes.

\begin{ex}
	La combinaison linéaire $4 \cdot \textsf{312} + \sqrt{2} \cdot \textsf{213}$ appartient à l'algèbre du groupe symétrique $S_3$ sur le corps $\C$.
\end{ex}

Lorsque nous étudions les mélanges d'une collection d'objets, nous agissons par une combinaison linéaire de permutations. Les opérateurs de mélange sont donc des éléments de l'algèbre du groupe symétrique, $\CSn$\index{Modules!$\CSn$}, (la preuve est donnée à la \autoref{ssec_contenant_par:melanges_alg_du_groupe}) et nous nous penchons ici sur la décomposition de cette algèbre en modules simples. Cette décomposition est bien connue depuis la première moitié du \mbox{XX$^e$ siècle}, et a d'abord été énoncée par Ferdinand Georg Frobenius et Issaï Schur, avant d'être reprise par Alfred Young, qui a exprimé la décomposition dans les termes que nous reprenons \cite[avant-propos par Paul Moritz Cohn]{JK}. L'outil combinatoire utilisé pour la comprendre est ce qu'on appelle les tableaux. Nous introduirons brièvement les tableaux, puis donnons la décomposition en modules simples de l'algèbre du groupe symétrique.

\subsection{Tableaux}\label{sec:tableaux}

\paragraph{Partages et diagrammes}
Pour un entier positif $n$, un \textit{partage $\lambda$ de $n$}\index{$\vdash$} est une suite décroissante d'entiers strictement positifs $(\lambda_1, \ldots, \lambda_l)$ dont la somme est $n$. À chaque partage $\lambda = (\lambda_1, \ldots, \lambda_l) \vdash n$, on associe un \textit{diagramme} formé de $\lambda_1$ boîtes dans la première rangée, de $\lambda_2$ boîtes placées dans la seconde (juste au-dessus de la première), et ainsi de suite. On place toutes les boîtes comme dans un escalier et on les aligne à gauche (comme dans l'\autoref{ex:diagram}); cet escalier compte donc $n$ boîtes au total et $l$ rangées. Des exemples de diagrammes et de partages sont donnés dans l'\autoref{ex:diagram}.

\begin{rmq}
	Au \autoref{chap:chap0}, nous nous sommes penchés notamment sur les compositions d'un entier $n$, soit des listes ordonnées de nombres naturels dont la somme est $n$. Les partages diffèrent des compositions puisqu'ils sont des listes placées en ordre décroissant. On peut associer à une composition un partage en triant la composition en ordre décroissant. Par exemple, la composition $(1,3,3,1) \vDash 8$ \index{$\vDash$} est envoyée sur le partage $(3,3,1,1) \vdash 8$ lorsqu'on la trie.
	Évidemment, ceci n'est pas une injection. Toujours en prenant l'exemple de $(1,3,3,1)$, six compositions sont associées au partage $(3,3,1,1)$~:
	\[ (1,1,3,3), \ (1,3,1,3),\ (1,3,3,1),\ (3,1,1,3),\ (3,1,3,1),\ (3,3,1,1).\]
	En général, à un partage $\lambda = (\lambda_1, \ldots, \lambda_l)$, on peut associer $\frac{l!}{\prod_{j\in [n]} \#\{i \mid \lambda_i = j\}!}$ compositions d'entiers.
\end{rmq}

\paragraph{Tableau}
Pour un diagramme donné, on peut remplir ses boîtes par des entiers positifs; le résultat est appelé un \textit{tableau} ou un \textit{tableau de Young}. Un tableau est un dit \text{standard} si les nombres inscrits dans les $n$ boîtes sont les entiers $1$ à $n$ et si, de plus, ils sont placés de la manière suivante~: lorsqu'on lit une ligne de gauche à droite ou de bas en haut, les nombres sont placés en ordre (strictement) croissant. Tous les tableaux standards de taille $2$ et $3$ sont listés à l'\autoref{ex:SYT}. 

\begin{ex}\label{ex:diagram}
	Le tuple $\lambda = (4,3,2,1)$ est un partage de $10$ et peut être représenté comme le premier diagramme ci-dessous. La deuxième illustration montre un tableau (pas standard) de forme $(3,1,1)$.
	\[ \YFrench \Yboxdim{16pt} \yng(4,3,2,1)\ , \qquad \young(933,1,4)\ .\]
\end{ex}

\begin{ex}\label{ex:SYT}
	Les tableaux standards de taille $2$ et $3$ sont 
	\[ \YFrench\Yboxdim{16pt}\young(12)\ ,\quad \young(1,2)\ , \quad \young(123)\ , \quad \young(12,3)\ , \quad \young(13,2)\ , \quad  \young(1,2,3)\ .\]
\end{ex}


\paragraph{Tableaux standards et suites croissantes d'une permutation}
Comme mentionné à la \autoref{ssec:op_sym}, les opérateurs de mélange symétrisés peuvent être décrits en termes de suites croissantes d'une permutation. En plus d'être utiles pour décomposer l'algèbre du groupe symétrique, les tableaux standards sont pratiques pour étudier les suites croissantes d'une permutation, grâce à une correspondance entre les permutations et les paires de tableaux standards d'une même forme, explicitée par Craige Schensted et Gilbert de Beauregard Robinson \cite{schensted}. 
Toutefois, c'est principalement pour leur lien avec la théorie de la représentation du groupe symétrique qu'ils sont présentés ici.

\subsection{Modules de Specht}\label{ssec:defn_Specht}
Les modules simples du groupe symétrique ont été définis par Wilhelm Specht \cite{specht}, et portent son nom. Chaque module de Specht est associé à un partage, et celui associé à $\lambda \vdash n$ est noté $S^\lambda$\index{Modules!$S^\lambda$}, par exemple. Les $\{S^\lambda\}_{\lambda \vdash n}$ forment la collection complète de tous les modules simples du groupe symétrique $S_n$, à isomorphisme près.

Quant à la construction de ces modules, elle est omise ici, mais se trouve dans divers ouvrages de référence, notamment au chapitre 2 de \cite{sagan}. Cependant, nous aurons besoin de connaître la dimension de chacun des modules~:
\begin{thm}[Conséquence du théorème 2.5.2 de \cite{sagan}]\label{thm:base_specht}
	Il existe une base du module $S^\lambda$ qui soit indexée par les tableaux standards de forme $\lambda$. 
\end{thm}
Ce dernier résultat nous donne en particulier la dimension du module de Specht $S^\lambda$~: c'est $f^\lambda$\index{$f^\lambda$}, le nombre de tableaux standards de forme $\lambda$.


\begin{ex}
	Le module trivial (présenté à l'\autoref{ex:modules_Sn}) est représenté par le diagramme formé d'une seule part (de taille $n$). Sa dimension est $1$, ce qui correspond au nombre de tableaux standards de cette forme; l'unique tableau est $\tab(123{\scriptstyle \ldots}n)$. Quant au module standard (aussi à l'\autoref{ex:modules_Sn}), il correspond au diagramme $(n-1,1)$ et les $n-1$ tableaux standards de cette forme ressemblent à 
	\[\Yboxdim{16pt}\tab(1{\scriptstyle\ldots}<{\scriptstyle i-1}><{\scriptstyle i+1}>{\scriptstyle \ldots}n,i)\ ,\]
	avec $i \in \{2, 3, \ldots, n\}$.
\end{ex}

Si nous avons présenté les modules de Specht, c'est surtout pour pouvoir décomposer l'algèbre du groupe symétrique en modules simples. Nous sommes maintenant en mesure de le faire.
\begin{thm}[Règle de Young - Décomposition de $\CSn$]\label{thm:regle_de_Young}
	L'algèbre du groupe symétrique sur $n$ éléments peut se décomposer en une somme directe de modules de Specht. Plus précisément,
	\[ \CSn \cong \bigoplus_{\lambda \vdash n} f^\lambda S^\lambda  \cong \bigoplus_{\substack{t \text{ tableau standard}\\\text{de taille }n}} S^{\mathrm{ forme}(t)}, \]
	où $f^\lambda$ est le nombre de tableaux standards de forme $\lambda$ et $S^\lambda$ est un module de Specht.
\end{thm}

\subsubsection{Exemple~: La décomposition de l'algèbre du groupe symétrique $S_3$}
L'algèbre du groupe symétrique $S_3$, de dimension $6$, peut être décomposée avec les tableaux standards. Il y a quatre tableaux standards contenant trois boîtes, dont $2$ de forme $(2,1)$~:
\[\tab(123)\ , \ \tab(12,3)\ , \ \tab(13,2)\ \text{ et }\tab(1,2,3)\ . \]
Ainsi, 
\[ \C S_3 \cong S^{\smalldia(3)} \oplus 2\ S^{\smalldia(2,1)} \oplus S^{\smalldia(1,1,1)}. \]

L'algèbre $\C S_3$ contient donc trois modules simples~:
\begin{itemize}
	\item un premier, $S^{\smalldia(3)}$, de dimension $1$, car c'est le nombre de tableaux standards de forme $(3)$. Il correspond au module trivial.
	\item un deuxième, $S^{\smalldia(2,1)}$, de dimension $2$, car il existe deux tableaux standards de cette forme. Il correspond au module standard. Sa multiplicité dans la décomposition est $2$, aussi parce qu'il y a deux tableaux standards de cette forme.
	\item un troisième, $S^{\smalldia(1,1,1)}$, aussi de dimension $1$. Il correspond au module qu'on appelle le module signe.
\end{itemize}

\paragraph{Valeurs propres et modules de Specht}
Une conséquence de la règle de Young est qu'il y autant de copies de modules simples dans $\CSn$ qu'il existe de tableaux standards de taille $n$. Dans les prochains chapitres, on se sert donc de cette bijection, notamment en conjonction avec un corollaire du lemme de Schur.

\begin{figure}
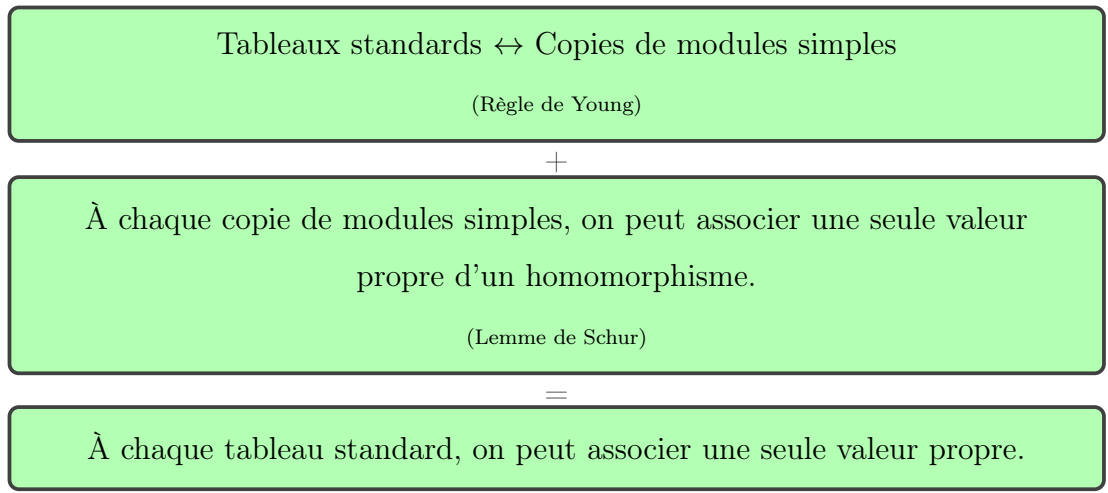

	\begin{center}
		\begin{tcolorbox}[colback=green!30!white]
			\begin{center}
				Tableaux standards $\leftrightarrow$ Copies de modules simples\\
				{\scriptsize (Règle de Young)}
			\end{center}
		\end{tcolorbox}
		\vspace{-4mm}
		+
		\begin{tcolorbox}[colback=green!30!white]
			\begin{center}
				À chaque copie de modules simples, on peut associer une seule valeur propre d'un homomorphisme.\\
				{\scriptsize (Lemme de Schur)}
			\end{center}
		\end{tcolorbox}
		\vspace{-4mm}
		=
		\begin{tcolorbox}[colback=green!30!white]
			\begin{center}
				À chaque tableau standard, on peut associer une seule valeur propre.
			\end{center}
		\end{tcolorbox}
	\end{center}
	\caption{Les valeurs propres de tous les homomorphismes du groupe symétrique sont indexées par les tableaux standards.}\label{fig:recap_vp_tableaux}
\end{figure}

Le lemme de Schur dit que les morphismes d'un module simple vers lui-même sont des homothéties, et ces morphismes n'ont donc qu'une seule valeur propre. Or, dans le cas du groupe symétrique, chaque copie d'un module simple est associée (de façon bijective) à un tableau standard. On peut donc associer à chaque tableau standard une valeur propre. La multiplicité d'une valeur propre associée à un tableau standard est alors la dimension du module de Specht correspondant. Celle-ci est donnée par le \autoref{thm:base_specht}, et est égale au nombre de tableaux standards ayant pour forme le partage auquel le module est associé. Le tout est schématisé à la \autoref{fig:recap_vp_tableaux}.

\subsection{Règle de branchement}\label{ssec:regle_branchement}
Rappelons que la technique que nous souhaitons utiliser pour calculer les valeurs propres utilise l'induction~: nous souhaitons en effet trouver une manière de calculer les valeurs propres pour les mélanges d'une collection de $n+1$ objets à partir des valeurs propres pour les mélanges d'une collection de $n$ objets. Ainsi, puisque nous savons ce qui se passe pour un nombre (très) petit d'objets, nous pouvons utiliser cette information pour comprendre les mélanges d'une plus grande collection. Pour faire l'induction, nous utilisons le fait qu'on peut projeter $S_n$ dans $S_{n+1}$ pour en faire un sous-groupe. Ce faisant, on  utilise l'induction pour passer de $S_n$ à $S_{n+1}$, et les modules de $\CSn$ \og deviennent \fg, en quelque sorte, des modules de $\C S_{n+1}$.

Le processus inverse, la \textit{restriction}, est aussi utilisé.

La \textit{règle de branchement} est le théorème qui nous permet de comprendre ces transformations entre $\CSn$ et $\C S_{n+1}$~:
\begin{thm}[Règle de branchement, \cite{James}]
	Soit $\lambda$ un partage de $n$. Le module induit de $S^\lambda$ dans $\C S_{n+1}$ se décompose comme la somme directe des modules $S^\mu$, où $\mu$ est un partage obtenu de $\lambda$ en ajoutant une boîte~:
	\[ S^\lambda\uparrow^{S_{n+1}}\index{Modules!$S^\lambda\uparrow^{S_{n+1}}$} \cong \bigoplus_{\text{Partages }\mu = \lambda+\Box}\index{$\lambda+\Box$} S^\mu. \]
	Le module restreint de $S_\lambda$ dans $\C S_{n-1}$ se décompose comme la somme directe des $S^\mu$, où $\mu$ est un partage de $\lambda$ obtenu en retirant une boîte~:
	\[ S^\lambda\downarrow_{S_{n-1}}\index{Modules!$S^\lambda\downarrow_{S_{n-1}}$} \cong \bigoplus_{\text{Partages }\mu = \lambda-\Box\index{$\lambda-\Box$}} S^\mu. \]
\end{thm}

\begin{ex}
	Considérons l'induction de $\C S_4 $ vers $\C S_5$ du module $S^{\smalldia(2,1,1)}$. À partir du diagramme $\smalldia(2,1,1)$\ , on peut ajouter une boîte dans la première, la deuxième ou quatrième ligne. On ne peut pas en ajouter une dans la troisième ligne, car $(2,1,2)$ n'est pas un partage. Ainsi,
	\[S^{\smalldia(2,1,1)}\uparrow^{S_5} \cong
	S^{\smalldia(3,1,1)} \oplus S^{\smalldia(2,2,1)} \oplus 
	S^{\smalldia(2,1,1,1)}.\]
	De la même façon, la restriction de $S^{\smalldia(2,1,1)}$ à $\C S_3$ se fait en enlevant une boîte dans la première ou la troisième ligne; enlever la boîte dans la deuxième ligne ne donne pas un diagramme~:
	\[ S^{\smalldia(2,1,1)}\downarrow_{S_3} \cong S^{\smalldia(2,1)} \oplus S^{\smalldia(1,1,1)}. \]
\end{ex}

\begin{figure}
	\begin{center}
		\begin{tikzpicture}
		\node (Csn1) [text width=5cm] at (-3,4) {Valeurs propres lorsqu'on agit sur $\C S_{n+1}$};
		\node (CSn) [text width=5cm] at (-3,0) {Valeurs propres lorsqu'on agit sur	$\C S_{n}$};
		\node (Slamb) [text width=5cm] at (5,0) {Valeurs propres lorsqu'on restreint à $S^\lambda, \ \lambda \vdash n$};
		\node (Slamb1) [text width=5cm] at (5,4) {Valeurs propres lorsqu'on restreint à $S^{\lambda+e_j}$};
		\draw[<->] [thick] (CSn) -- node[right]{\footnotesize Objectif général} (Csn1);
		\draw[<->] [thick, dotted] (Slamb) --node[right]{\footnotesize $\stackunder{\text{Objectif atteint }}{\stackunder{\text{(satisfait l'objectif}}{\text{général)}}}$} (Slamb1);
		\draw[->] (CSn) -- (1.5, -0.5) node{};
		\draw[->] (CSn) -- (Slamb) node{};
		\draw[->] (CSn) -- (1.5, 0.5) node{};
		\draw[->] (Csn1) -- (1.5, 3.5) node{};
		\draw[->] (Csn1) -- (Slamb1) node{};
		\draw[->] (Csn1) -- (1.5, 4.5) node{};
		\end{tikzpicture}
	\end{center}
	\caption{Schéma de la recherche des valeurs propres par induction et par restriction.}\label{fig:vp_induction}
\end{figure}
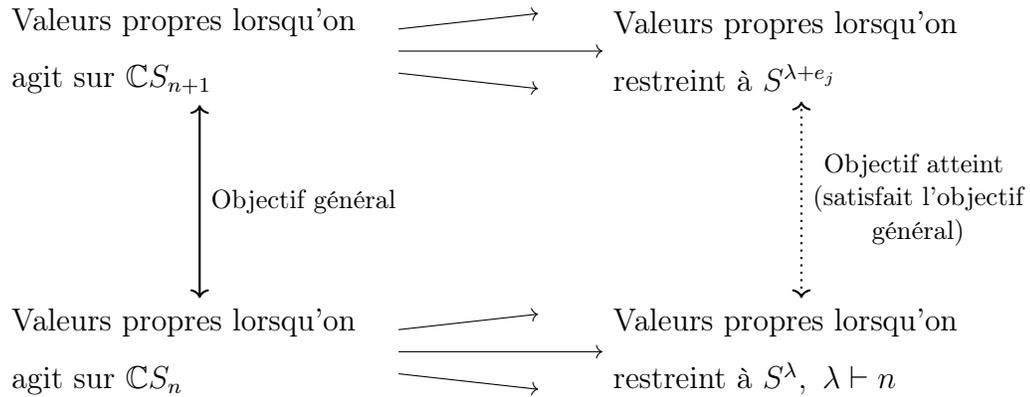

La règle de branchement nous dit donc que, pour l'induction, on pourra se contenter de trouver une procédure d'induction de $S^\lambda$ à $S^{\lambda+\Box}$, plutôt que sur toute l'algèbre. Le schéma de la quête de valeurs propres par l'induction ressemble plutôt à ce qui est à la \autoref{fig:vp_induction}.

\section{Caractères du groupe symétrique}\label{sec:caracteres}
Les caractères d'un groupe sont un outil très puissant pour comprendre la structure du groupe; c'est d'ailleurs la raison pour laquelle ils ont été inventés. En plus de cette application, ils peuvent décrire les valeurs propres de certains opérateurs de mélange, sous certaines conditions. La principale condition est que la probabilité d'agir par un élément du groupe soit constante au sein même d'une classe de conjugaison; voir, par exemple, \cite{DiSh} pour un exemple d'un tel opérateur de mélange, ou \cite{SC_review} pour plus de détails sur les techniques. Malheureusement, les opérateurs de mélange symétrisés ne satisfont pas cette propriété, mais les caractères seront tout de même utiles.

\begin{defn}
	Soit $\varphi$ un morphisme de groupes de $G$ vers $GL(\C)$,\index{Modules!$GL(\C)$} le groupe des matrices inversibles à coefficients complexes.
	Le \textit{caractère} $\chi^\varphi$\index{$\chi^\varphi$} de ce morphisme est la fonction qui associe à chaque élément $g$ du groupe la trace de la matrice $\varphi(g)$.
\end{defn}

\begin{prop}
	Le caractère d'un groupe est constant sur une classe de conjugaison~: c'est ce qu'on appelle une fonction centrale.
\end{prop}
\begin{proof}
	En effet, $\varphi$ est un morphisme de groupes, et \mbox{$\varphi(hgh^{-1}) = \varphi(h)\varphi(g)(\varphi(h))^{-1}$}. De plus, la trace étant invariante par similitude, \mbox{$\chi^\varphi(hgh^{-1})=\tr(\varphi(h)\varphi(g)(\varphi(h))^{-1}) = \tr(\varphi(g)) = \chi^ \varphi(g)$}.\index{$\tr(\cdot)$}
\end{proof}

\paragraph{Que des nombres entiers~: la règle de Murnaghan-Nakayama}

Les caractères du groupe symétrique ont ceci de particulier qu'ils sont tous des nombres réels et entiers et qu'ils peuvent être décrits de façon combinatoire, bien qu'ils ne soient pas toujours positifs. Ils sont en effet décrits par le retrait d'un ensemble de boîtes dans un diagramme. La procédure qui décrit le calcul des caractères du groupe symétrique est la règle de Murnaghan-Nakayama et peut facilement être trouvée dans les ouvrages de référence (par exemple, \cite[\S 4.10]{sagan} et \cite[\S 7.17]{EC2})

\begin{corl}[de la règle de Murnaghan-Nakayama]\label{corl:caracteres_Sn_entiers}
	Les caractères du groupe symétrique sont des entiers.
\end{corl}

À titre d'exemples, tous les caractères irréductibles des groupes $S_3$ et $S_4$ sont représentés aux tableaux \ref{tab:_table_caracteres_S_3} et \ref{tab:_table_caracteres_S_4}. Comme les caractères sont des fonctions centrales, une seule permutation par classe de conjugaison est représentée dans les tables.
Les caractères des groupes  $\mathbb{Z}/ 4\mathbb{Z}$ et  $\mathbb{Z}/ 5\mathbb{Z}$ sont représentés aux tableaux \ref{tab:table_caracteres_Z4} et \ref{tab:table_caracteres_Z5}, où l'on peut observer qu'ils ne sont pas tous des entiers.

Le produit de caractères et les relations d'orthogonalité entre les caractères sont présentés plus loin, à la \autoref{ssec:proj_iso}.

\begin{table}[h]
	\centering
	\begin{minipage}{0.42\textwidth}
		\caption{Table de caractères de $S_3$.}\label{tab:_table_caracteres_S_3}
		\begin{equation*}
		\begin{tabular}{l|rrrrr}
		& $\Id$ & $(1\,2)$ & $(1\,2\,3)$ \\
		\hline
		$\chi^{\smalldia(3)}$ & $1$ & $1$ & $1$ \\
		$\chi^{\smalldia(1,1,1)}$ & $1$ & $-1$ & $1$ \\
		$\chi^{\smalldia(2,1)}$ & $2$ & $0$   & $-1$  \\
		\end{tabular}
		\end{equation*}
	\end{minipage}
	\quad
	\begin{minipage}{0.52\textwidth}
		\caption{Table de caractères de $S_4$.}\label{tab:_table_caracteres_S_4}
		\begin{equation*}
		\begin{tabular}{l|rrrrr}
		& ${\scriptstyle \Id }$ & {\scriptsize $(1\,2)$} & {\scriptsize $(1\,2)(3\,4)$} & {\scriptsize $(1\,2\,3)$} & {\scriptsize $(1\,2\,3\,4)$} \\
		\hline
		$\chi^{\smalldia(4)}$ & $1$ & $1$ & $1$ & $1$ & $1$ \\
		$\chi^{\smalldia(1,1,1,1)}$ & $1$ & $-1$ & $1$ & $1$ & $-1$ \\
		$\chi^{\smalldia(2,2)}$ & $2$ & $0$ & $2$ & $-1$ & $0$ \\
		$\chi^{\smalldia(2,1,1)}$ & $3$ & $-1$ & $-1$ & $0$ & $1$ \\
		$\chi^{\smalldia(3,1)}$ & $3$ & $1$ & $-1$ & $0$ & $-1$ \\
		\end{tabular}
		\end{equation*}
	\end{minipage}
\end{table}

\begin{table}[h]
	\begin{minipage}{0.45\textwidth}
		\centering
		\caption{Table de caractères de $\mathbb{Z}/ 4\mathbb{Z}$.}\label{tab:table_caracteres_Z4}
		\begin{equation*}
		\begin{array}{c|rrrr}
		& 0 & 1  & 2  & 3  \\ \hline
		\chi^1 & 1 & 1  & 1  & 1  \\
		\chi^2 & 1 & i  & -1 & -i \\
		\chi^3 & 1 & -1 & 1  & -1 \\
		\chi^4 & 1 & -i & -1 & i  \\
		\end{array}
		\end{equation*}
	\end{minipage}
	\hfill
	\begin{minipage}{0.45\textwidth}
		\caption[Table de caractères de $\mathbb{Z}/5\mathbb{Z}$.]{Table de caractères de $\mathbb{Z}/5\mathbb{Z}$, où $\zeta = e^{2 \pi i / 5}$.}\label{tab:table_caracteres_Z5}
		\begin{equation*}
		\begin{tabular}{l|rrrrr}
		& $0$ & $1$ & $2$ & $3$ & $4$ \\
		\hline
		$\chi^1$ & $1$   & $1$       & $1$       & $1$       & $1$       \\[0.5ex]
		$\chi^2$ & $1$   & $\zeta^1$ & $\zeta^2$ & $\zeta^3$ & $\zeta^4$ \\[0.5ex]
		$\chi^3$ & $1$   & $\zeta^2$ & $\zeta^4$ & $\zeta^1$ & $\zeta^3$ \\[0.5ex]
		$\chi^4$ & $1$   & $\zeta^3$ & $\zeta^1$ & $\zeta^4$ & $\zeta^2$ \\[0.5ex]
		$\chi^5$ & $1$   & $\zeta^4$ & $\zeta^3$ & $\zeta^2$ & $\zeta^1$ \\
		\end{tabular}
		\end{equation*}
	\end{minipage}
\end{table}

\chapter{Une première famille d'opérateurs}\label{chap:1efamille}

On entre maintenant dans le vif du sujet. Des deux familles d'opérateurs sur lesquels porte principalement cette thèse, la première famille est certainement celle sur laquelle nous savons le mieux expliquer les phénomènes. Elle était déjà bien étudiée dans le chapitre VI de \cite{RSW}, mais certaines conjectures étaient laissées en plan. Ce chapitre donne la réponse à celles-ci. La plus importante est que les valeurs propres de ces opérateurs sont toutes entières, et c'est le contenu du \autoref{corl:vp_nuk_entiers}. Le résultat principal, qui décrit comment calculer explicitement les valeurs propres des opérateurs $\allnuk$, est décrit ici (\autoref{thm:main}), mais sa preuve est présentée au \autoref{chap:preuves}.
\section{Définition}
On peut définir les opérateurs de la famille $\allnuk$ de plusieurs façons. Trois sont décrites dans cette section, et l'équivalence des définitions est démontrée. La forme matricielle de ces opérateurs linéaires est aussi donnée. Une quatrième définition, plus technique, est présentée au \autoref{chap:preuves}. 

\subsection{Comme un humain qui mélange...}
La première définition est la plus utile, à la fois pour raconter les opérateurs de mélange symétrisés, mais aussi pour comprendre ce qui se passe. Il s'agit en fait de sa définition en tant qu'opérateur de mélange. Si on voulait mélanger une collection finie d'objets avec un opérateur $\nu_k$, c'est ce qu'on ferait, concrètement.

\begin{defn} \label{defn:nu_k_cartes}\index{Opérateurs!$\nu_k$}
Soit $w = (w_1,\ldots, w_n)$ une suite ordonnée d'objets, des lettres dans un mot, par exemple. L'opérateur $\Nu{k}$ retire $k$ des $n$ objets de la suite $w$, puis les réinsère un après l'autre de telle façon que l'ordre de ces $k$ objets n'est pas nécessairement préservé. On encode le résultat de cette opération dans une combinaison linéaire formelle des permutations de  $w_1, \ldots, w_n$ que l'on peut obtenir en réalisant ce mélange. Dans cette combinaison linéaire, le coefficient d'une suite donnée est le nombre de façons dont on peut obtenir cette suite en une itération du mélange.
\end{defn}

\begin{figure}
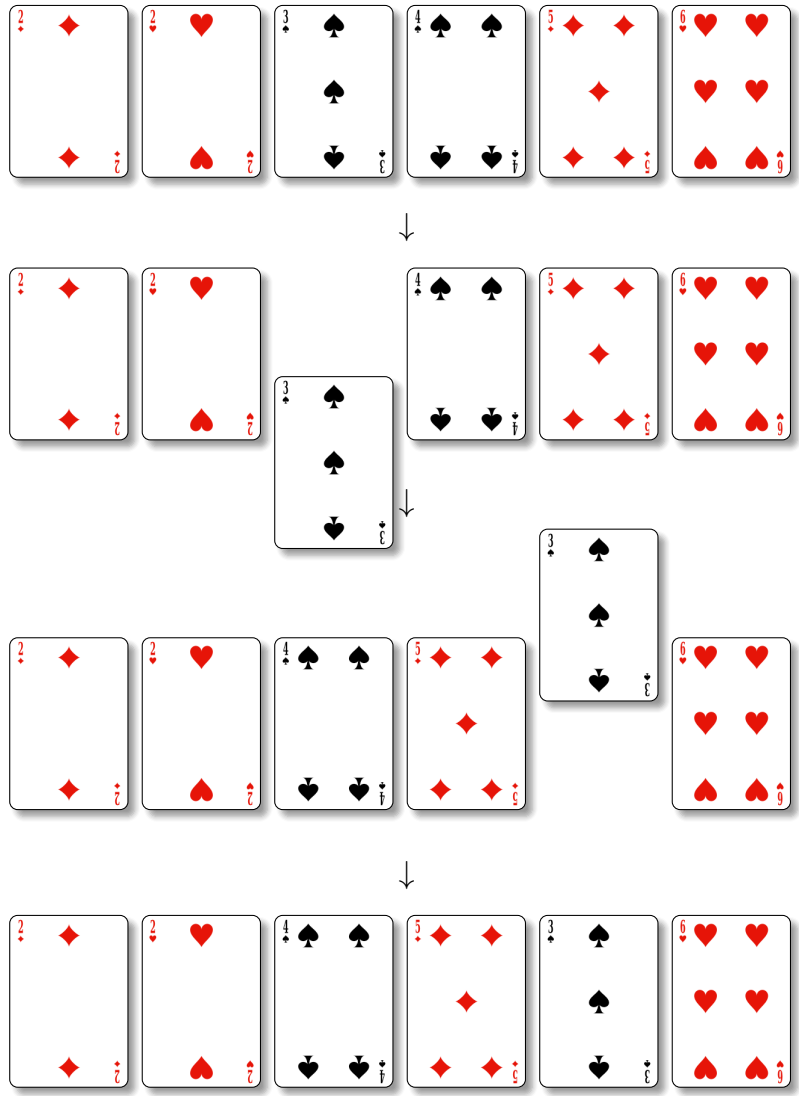

	\[\begin{array}{c}
		\begin{array}{c}\carda\cardb\cardc\cardd\carde\cardf\end{array}\\
		\vspace{2mm}
		\downarrow\\
		\vspace{-20mm}
		\begin{array}{c}\carda\cardb{\raisebox{-3.5em}{\cardc}}\cardd\carde\cardf\end{array}\\
		\vspace{2mm}\\
		\downarrow\\
		\vspace{2mm}
		\begin{array}{c}\carda\cardb\cardd\carde{\raisebox{+3.5em}{\cardc}}\cardf\end{array}\\
		\vspace{2mm}
		\downarrow\\
		\vspace{2mm}
		\begin{array}{c}\carda\cardb\cardd\carde{\raisebox{-0.0em}{\cardc}}\cardf
		\end{array}
	\end{array}\]
	\caption{Le mélange doublement aléatoire déplace une seule carte à la fois.}
\end{figure}

\begin{ex}
	Appliquer $\Nu{1}$ à la suite \textsf{abc} revient à retirer une des trois lettres et à la réinsérer dans un des trois espaces restants. Par exemple, en retirant la lettre \textsf{a}, les suites \textsf{abc}, \textsf{bac} et \textsf{bca} peuvent chacune être obtenue d'une façon. En retirant plutôt \textsf{b} ou \textsf{c}, on obtient
	\[\Nu{1}(\textsf{abc}) = 3 \cdot \textsf{abc} + 2 \cdot \textsf{bac} + 2 \cdot \textsf{acb} + \textsf{cab} + \textsf{bca}.\]
\end{ex}

\begin{ex}
	L'opérateur $\Nu{2}$ agir sur le mot \textsf{abcd} en déplaçant deux lettres, les lettres pouvant être \og déplacées \fg \ à l'endroit où elles étaient déjà. En supposant que les lettres \textsf{a}, \textsf{b}, \textsf{c} et \textsf{d} soient toutes distinctes, le coefficient de \textsf{cabd} dans $\Nu{2}(\textsf{abcd})$ est $4$. Les paires de lettres à déplacer pour obtenir \textsf{cabd} sont soit $\{\textsf{a,b}\}$, $\{\textsf{a,c}\}$, $\{\textsf{b,c}\}$ ou $\{\textsf{c,d}\}$.

	\begin{center}
		\begin{tikzpicture}[xscale=0.4]
		\node (ZERO) at (0.5,1) {};
		\node (A) at (1,1) {\textsf{a}};
		\node (B) at (2,1) {\textsf{b}};
		\node (C) at (3,1) {\textsf{c}};
		\node (D) at (4,1) {\textsf{d}};
		\node (C2) at (3.8,1) {};
		
		\draw (0.5,0.8) rectangle (2.5,1.2);
		\draw[->, >=latex] (1.5,1.2) to [out=30,in=150] (C2);
		\end{tikzpicture}
		\quad
		\begin{tikzpicture}[xscale=0.4]
		\node (ZERO) at (-0.3,1) {};
		\node (A) at (1,1) {\textsf{a}};
		\node (B) at (2,1) {\textsf{b}};
		\node (C) at (3,1) {\textsf{c}};
		\node (C1) at (3.5,1) {};
		\node (D) at (4,1) {\textsf{d}};
		
		\draw[->, >=latex] (C1) to [out=150,in=30] (ZERO);
		\draw[->, >=latex] (A) to [loop above, looseness=2] 
		node{} (A);
		\draw (0.5,0.8) rectangle (1.5,1.2);
		\draw (2.5,0.8) rectangle (3.5,1.2);
		\end{tikzpicture}
		\quad
		\begin{tikzpicture}[xscale=0.4]
		\node (ZERO) at (-0.3,1) {};
		\node (A) at (1,1) {\textsf{a}};
		\node (B) at (2,1) {\textsf{b}};
		\node (C) at (3,1) {\textsf{c}};
		\node (C1) at (3.5,1) {};
		\node (D) at (4,1) {\textsf{d}};
		
		\draw[->, >=latex] (C1) to [out=150,in=30] (ZERO);
		\draw[->, >=latex] (B) to [loop above, looseness=2] 
		node{} (B);
		\draw (1.5,0.8) rectangle (2.5,1.2);
		\draw (2.5,0.8) rectangle (3.5,1.2);
		\end{tikzpicture}
		\quad
		\begin{tikzpicture}[xscale=0.4]
		\node (ZERO) at (-0.3,1) {};
		\node (A) at (1,1) {\textsf{a}};
		\node (B) at (2,1) {\textsf{b}};
		\node (C) at (3,1) {\textsf{c}};
		\node (C1) at (3.5,1) {};
		\node (D) at (4,1) {\textsf{d}};
		
		\draw[->, >=latex] (C1) to [out=150,in=30] (ZERO);
		\draw[->, >=latex] (D) to [loop above, looseness=3] 
		node{} (D);
		\draw (3.5,0.8) rectangle (4.5,1.2);
		\draw (2.5,0.8) rectangle (3.5,1.2);
		\end{tikzpicture}
	\end{center}
	Ainsi,
	\[\Nu{2}(\textsf{abcd}) = \ldots + 4 \cdot \textsf{cabd} + \ldots \]
\end{ex}

L'opérateur $\Nu{0}$ est l'identité; il y a en effet une seule façon de ne retirer aucun objet dans la suite. Quant à $\nu_1$, c'est le mélange doublement aléatoire, présenté à la \autopageref{par:r2r_def}.

\begin{rmq}
	 Ces opérateurs ont aussi été notés sous la forme $\nu_{(n-k,1^k)}$ \index{Opérateurs!$\nu_{(n-k,1^k)}$}
	 \cite{RSW}. Dans le contexte de leur article, il est naturel que les auteurs utilisent cette notation puisqu'ils couvrent le cas des opérateurs $\{\nu_{\lambda}\}_{\lambda \vdash n}$ indexés par les partages d'entiers. 

     Dans ce chapitre et le suivant, nous ne nous intéresserons qu'aux partages de la forme $(n-k, 1^k)$. Au besoin, nous préciserons la valeur de $n$. Au \autoref{chap:2efamille}, nous étudions les partages d'une autre forme, mais les opérateurs portent un autre nom.
\end{rmq}

\subsection{Suites croissantes d'une permutation}
La définition de ces opérateurs en termes de mélanges est pratique pour développer une certaine intuition, mais son utilité décroît rapidement 
lorsqu'on souhaite écrire une définition formelle des matrices de ces 
opérateurs linéaires. Pour pallier cette lacune, on peut utiliser une 
définition alternative des opérateurs $\allnuk$ utilisant les suites 
croissantes d'une permutation, définies à la \autopageref{sssec:non-inversions}. On note le nombre de suites croissantes de longueur $i$ dans une permutation $\sigma$ par $\noninv_i(\sigma)$\index{$\noninv_i$}. Notons que ces suites ne doivent pas nécessairement être disjointes.

\begin{rmq}
	La notation $\noninv_i$ provient de \cite{RSW} et est nommée ainsi en référence aux $i$-non-inversions, comme elles sont appelées dans ce texte. Une suite décroissante de longueur $2$ dans une permutation étant une inversion, une non-inversion serait une suite croissante. C'est pourquoi une suite croissante de longueur $i$ est appelée une $i$-non-inversion.
\end{rmq}

Rappelons que, dans une permutation $\sigma$, une \textit{suite croissante de longueur $k$} est un ensemble $\{i_1, \ldots, i_k \}$, où à la fois les suites $(i_1, i_2, \ldots, i_k)$ et $(\sigma^{-1}(i_1), \ldots,\sigma^{-1}(i_k))$ sont placées en ordre croissant. En d'autres mots, il s'agit d'un sous-mot formé de lettres placées en ordre croissant dans le mot qu'est la permutation.
\begin{prop}\label{prop:noninv_sym}
	Soit $\sigma$ une permutation. Le nombre de suites croissantes de longueur $k$ est le même dans $\sigma$ et dans son inverse, $\sigma^{-1}$.
\end{prop}

\begin{proof}
	Pour le démontrer, il suffit de construire une suite croissante de longueur $k$ dans $\sigma^{-1}$ pour chacune de celles dans $\sigma$.
    
    Soit $i_1 < i_2 < \ldots < i_k$ une suite croissante de $\sigma$. Alors $\sigma^{-1}(i_1) < \ldots < \sigma^{-1}(i_k)$.
    
    Posons $j_1 = \sigma^{-1}(i_1), \ldots, j_k = \sigma^{-1}(i_k)$.
    Alors, $(\sigma(j_1), \ldots, \sigma(j_k)) = (i_1, \ldots, i_k)$ et $(j_1, \ldots, j_k)$ est une suite croissante de longueur $k$ dans $\sigma^{-1}$.
\end{proof}

Ce dernier résultat nous permet, plus loin, de démontrer que les opérateurs $\allnuk$ sont symétriques. Pour cela, voici une autre définition de ces opérateurs, donnée par \cite{RSW}.

\begin{defn}\label{defn:nu_k_noninv_perm}\index{Opérateurs!$\nu_k$}
	 Soit $\sigma$ une permutation de $n$ éléments et notons $\noninv_i(\sigma)$ le nombre de séquences croissantes de longueur $i$ dans $\sigma$.
	 Alors, $\Nu{k}$ agit sur $\sigma$ de la façon suivante~:
	\[ \Nu{k}(\sigma) = \sum_{\tau \in S_n} \noninv_{n-k}(\sigma^{-1}\tau) \tau. \]
\end{defn}

\begin{rmq}
	Cette définition fonctionne bien pour $\nu_k(\sigma)$, mais, contrairement à la \autoref{defn:nu_k_cartes}, elle ne permet pas de décrire $\nu_k(w)$ sur un mot $w$ qui admettrait des répétitions de lettres. Cette lacune sera comblée par la \autoref{defn:nu_k_noninv_mots}.
\end{rmq}

\begin{proof}[Preuve de l'équivalence des définitions \ref{defn:nu_k_cartes} et \ref{defn:nu_k_noninv_perm}.]
	Pour simplifier la preuve, on suppose que notre mot de départ est la permutation $\sigma = \Id$. Notre preuve pourra alors être étendue à toutes les permutations étant donné que la \autoref{defn:nu_k_noninv_perm} tient compte de la permutation de départ. Pour les mots qui ne sont pas des permutations, la définition est donnée plus loin (\autoref{defn:nu_k_noninv_mots}) et l'équivalence est alors démontrée.
	
	On peut alors écrire formellement la \autoref{defn:nu_k_cartes}. Pour ce faire, on rappelle une notation. Le \textit{battage}, noté $u \shuffle v$, est la combinaison linéaire formelle de tous les mots qui contiennent exactement $u$ et $v$ comme sous-mots disjoints. 
	Alors,
	\begin{align*}
	\Nu{k}(\mathsf{12\ldots n}) = \sum_{1\leq i_1 < \ldots < i_k \leq n}  \sum_{\eta \in S_k} \mathsf{12\ldots (i_1-1)(i_1+1)}&\mathsf{\ldots (i_k-1) (i_k+1) \ldots n} \\& \shuffle\ \mathsf{i_{\eta(1)}\cdots i_{\eta(k)}}.
	\end{align*}
	
	Pour montrer l'équivalence avec la \autoref{defn:nu_k_noninv_perm}, il faut que le coefficient de la permutation $\tau$ dans cette équation soit $\noninv_{n-k}(\tau)$. Or, le nombre de suites croissantes de longueur $(n-k)$ dans $\tau$ est calculé ainsi~: c'est le nombre de sous-mots $\mathsf{12\ldots (i_1-1)(i_1+1)\ldots (i_k-1)(i_k+1)\ldots n}$ tels que \mbox{$\tau(1)\tau(2)\ldots \tau(i_1-1)\tau(i_1+1)\ldots \tau(i_k-1)\tau(i_k+1)\ldots \tau(n)$} est également un sous-mot. Ensuite, pour une telle sous-suite croissante, il n'y a qu'un choix pour $\eta$ et les positions auxquelles insérer $\mathsf{i_{\eta(1)}, \ldots, i_{\eta(k)}}$ pour obtenir $\tau$. Ainsi, on peut réécrire
	\[ \Nu{k}(\Id) =  \sum_{\tau} \noninv_{n-k}(\tau) \tau, \]
	ce qui est la \autoref{defn:nu_k_noninv_perm}.
\end{proof}

\subsubsection{Une définition analogue sur les mots}

Comme noté précédemment, la \autoref{defn:nu_k_noninv_perm} ne fonctionne que sur les permutations, n'étant pas bien définie sur les mots qui auraient des répétitions puisque ces mots ne peuvent pas être inversés. Nous présentons une définition alternative.
\begin{defn}[Alternative à \autoref{defn:nu_k_noninv_perm}]\label{defn:nu_k_noninv_mots}\index{Opérateurs!$\nu_k$}
	Soit $w$ un mot de taille $n$. Alors,
	\[\Nu{k}(w) = \sum_{\tau \in S_n} \noninv_{n-k}(\tau)\ w \cdot \tau, \]
	où $w \cdot \tau = w_{\tau(1)}\ldots w_{\tau(n)}$.
\end{defn}

\begin{proof}[Cohérence des définitions \ref{defn:nu_k_noninv_perm} et \ref{defn:nu_k_noninv_mots}.]
	Cette définition découle de la \autoref{defn:nu_k_noninv_perm}. En effet, on peut réécrire, pour $\sigma \in S_n$,
	\begin{align*}
	\nu_k(\sigma) &= \sum_{\tau \in S_n} \noninv_{n-k}(\sigma^{-1}\tau)\ \tau\\
	& = \sum_{\sigma\tau \in S_n} \noninv_{n-k}(\sigma^{-1}\sigma\tau)\ \sigma\tau\\
	& = \sum_{\sigma\tau \in S_n} \noninv_{n-k}(\tau)\ \sigma\tau\\
	& = \sum_{\tau \in S_n} \noninv_{n-k}(\tau)\ \sigma\tau.\\
	\end{align*}
	La dernière égalité découle du fait qu'on somme sur toutes les permutations, alors que la deuxième provient d'un changement de variable.
	
	Cette définition ne demandant pas qu'on \og inverse\fg\ un mot, elle peut être directement étendue des permutations aux mots~:
	\[ \Nu{k}(w) =  \sum_{\tau \in S_n} \noninv_{n-k}(\tau)\ w \cdot \tau. \qedhere\]
\end{proof}

\subsection{Matrices}\label{ssec:def_matrices_nuk}
Ces dernières définitions de $\nu_k$ nous permettent de définir facilement les matrices représentant les opérateurs linéaires $\nu_k$. Rappelons que l'espace vectoriel dans lequel on se place est l'algèbre du groupe symétrique, $\CSn$, et que les éléments de cet espace sont tous des combinaisons linéaires de permutations des nombres de $1$ à $n$.
Nous définissons d'abord les matrices pour les permutations, qui forment une base de $\CSn$, puis nous généralisons aux mots avec répétitions possibles (sur un autre espace vectoriel). La matrice $M_k$ est telle que $(\vec{v} \cdot M_k)_\sigma = \nu_k(\sigma)$ lorsque $\vec{v}$ est le vecteur ligne de toutes les permutations de $S_n$, dans le même ordre que les lignes de la matrice.

L'entrée $(\sigma, \tau)$ de la matrice $M_k$\index{$M_k$} associée à $\nu_k$ vaut $\noninv_{n-k}(\sigma^{-1}\tau)$. On a alors bien que le vecteur de toutes les permutations  multiplié par la colonne $\sigma$ vaut $\sum_{\tau \in S_n} \noninv_{n-k}(\sigma^{-1}\tau) \tau = \nu_k(\sigma)$.\\

\begin{ex}\label{ex:matrice_R2R}
	La matrice de l'opérateur linéaire $\Nu{1}$ sur les permutations de longueur $3$ est
	\[\bordermatrix{
		&{\scriptstyle \textsf{abc}} & {\scriptstyle \textsf{acb}} & 
		{\scriptstyle \textsf{bac}} & {\scriptstyle \textsf{bca}} & 
		{\scriptstyle \textsf{cab}} & {\scriptstyle \textsf{cba}} \cr
		{\scriptstyle \textsf{abc}} & 3 & 2 & 2 & 1 & 1 & 0 \cr
		{\scriptstyle \textsf{acb}} & 2 & 3 & 1 & 0 & 2 & 1 \cr
		{\scriptstyle \textsf{bac}} & 2 & 1 & 3 & 2 & 0 & 1 \cr
		{\scriptstyle \textsf{bca}} & 1 & 0 & 2 & 3 & 1 & 2 \cr
		{\scriptstyle \textsf{cab}} & 1 & 2 & 0 & 1 & 3 & 2 \cr
		{\scriptstyle \textsf{cba}} & 0 & 1 & 1 & 2 & 2 & 3 \cr
	}.\]
\end{ex}

Pour les mots avec répétitions, on indexe les lignes (et les colonnes) par tous les mots qu'on peut obtenir en permutant les positions d'un mot. Ces mots sont tous les mots qui comportent le même nombre d'occurrences de chacune des lettres. Encore une fois, si $\vec{v}$ représente le vecteur (ligne) de tous les états (les mots) possibles et si $M_k$ est la matrice de $\Nu{k}$, on souhaite que l'unique élément de la colonne $w$ de $\vec{v} \cdot M_k$ soit l'élément $\Nu{k}(w)$; un exemple est donné avec $\nu_2$ et le mot $\mathsf{abab}$ (\autoref{ex:matrice_mots}). Étant donné la définition de $\Nu{k}(w)$, l'entrée $(w,u)$ de $M_k$ vaut $ \sum_{\substack{\tau \in S_n,\\ w \cdot \tau = u}} \noninv_{n-k}(\tau)$.

\begin{ex}\label{ex:matrice_mots}
	Cette matrice décrit l'opérateur $\Nu{2}$ sur les mots contenant deux occurrences de $\mathsf{a}$ et deux occurrences de $\mathsf{b}$~:
	\[ \Nu{2} = \bordermatrix{
	&{\scriptstyle \textsf{aabb}} & {\scriptstyle \textsf{abab}} & {\scriptstyle \textsf{baab}} & {\scriptstyle \textsf{abba}} &
	{\scriptstyle \textsf{baba}} & {\scriptstyle \textsf{bbaa}}\cr
	{\scriptstyle \textsf{aabb}} & 20 & 16 & 12 & 12 & 8 & 4 \cr
	{\scriptstyle \textsf{abab}} &16 & 14 & 12 & 12 & 10 & 8 \cr
	{\scriptstyle \textsf{baab}} & 12 & 12 & 12 & 12 & 12 & 12 \cr
	{\scriptstyle \textsf{abba}} & 12 & 12 & 12 & 12 & 12 & 12 \cr
	{\scriptstyle \textsf{baba}} & 8 & 10 & 12 & 12 & 14 & 16 \cr
	{\scriptstyle \textsf{bbaa}} & 4 & 8 & 12 & 12 & 16 & 20\\}.\]

	De plus, soit le vecteur $\vec{v} = (\textsf{aabb}, \textsf{abab}, \textsf{baab}, \textsf{abba}, \textsf{baba}, \textsf{bbaa})$.\\
	Ici, l'unique élément de la deuxième colonne du vecteur ligne $\vec{v} \cdot \Nu{2}$ vaut
	\[16 \cdot \textsf{aabb} + 14 \cdot \textsf{abab} + 12 \cdot  \textsf{baab} + 12 \cdot \textsf{abba} + 10 \cdot \textsf{baba} + 8 \cdot \textsf{bbaa},\]
	soit $\nu_2(\textsf{abab})$.
\end{ex}

Ces matrices étant définies, on peut analyser leurs propriétés. Quelques-unes d'entre elles se retrouvent ci-dessous.
\paragraph{Matrices entières}

Les matrices ainsi formées contiennent uniquement des entiers positifs, étant donné que leurs entrées comptent le nombre de suites croissantes d'une permutation.

\paragraph{Matrices symétriques}
 On a remarqué, à la \autoref{prop:noninv_sym} que $\noninv_k(\sigma^{-1}\tau) = \noninv_k(\tau^{-1}\sigma)$. Or, $\noninv_k(\sigma^{-1}\tau)$ est l'entrée $(\sigma, \tau)$ de la matrice de $\Nu{k}$. Ainsi, les entrées $(\sigma,\tau)$ et $(\tau,\sigma)$ de la matrice sont égales et la matrice est symétrique. Comme elle est aussi réelle, elle est hermitienne.

\paragraph{Matrices stochastiques... ou presque!}
Celles et ceux à qui les chaînes de Markov sont familières remarqueront peut-être que la matrice ci-dessus n'est pas une matrice de transition. Dans une matrice de transition, les entrées sont des probabilités de transition, et prennent donc des valeurs entre $0$ et $1$. De plus, les entrées d'une même ligne d'une matrice de transition ont toujours pour somme $1$, puisqu'il s'agit de la somme des probabilités pour tous les événements qui peuvent se produire lorsque l'on se trouve dans un certain état.

Cela dit, dans l'\autoref{ex:matrice_R2R}, la somme de chaque ligne est toujours $9$ et toutes les entrées sont comprises entre $0$ et $9$. Pour retrouver une matrice de transition, on pourrait donc diviser chacune des entrées par $9$ (elles ne seraient alors plus des entiers). 

Chaque ligne de la matrice $\Nu{k}$ représente ce qui peut se produire à partir d'un état donné. Chaque état correspond à une suite $w = (w_1,\ldots, w_n)$, et les événements obtenus après une transition le sont en~:
\begin{itemize}
	\item sélectionnant $k$ lettres à retirer parmi $w_1,\ldots, w_n$ ($\binom{n}{k}$ possibilités);
	\item permutant ces $k$ lettres ($k!$ possibilités);
	\item réinsérant les $k$ lettres pour former une nouvelle suite $w' = (w_1',\ldots, w_n')$ ($\binom{n}{k}$ possibilités).
\end{itemize}  
Le nombre d'événements possibles, en comptant les répétitions, est ainsi $\binom{n}{k}^2 k!$, ce qui correspond au produit du nombre de façons d'accomplir chacune des trois étapes ci-haut. Ainsi, $\frac{\Nu{k}}{\binom{n}{k}^2 k!}$ est la matrice de transition de la chaîne de Markov associée à $\Nu{k}$.

Il est connu que les matrices des chaînes de Markov ont des valeurs propres complexes dont le module est inférieur ou égal à 1. Il est de plus connu que la multiplicité de la valeur propre $1$ est toujours $1$ pour les chaînes de Markov qui satisfont deux propriétés (irréductibilité et apériodicité), satisfaites par les opérateurs de mélange symétrisés, de même que par beaucoup d'opérateurs de mélange. 
Normaliser les entrées de la matrice nous permet d'avoir des valeurs propres entières, dont la plus grande est $\binom{n}{k}^2 k!$.

\paragraph{Matrices semi-définies positives}\label{par:nuk_semi_def_pos}
Une matrice $M$ est \textbf{semi-définie positive} si, pour tout vecteur colonne à coefficients complexes $z$, le nombre $z^* Mz$ est réel et positif, où $z^*$ désigne le vecteur ligne dans lequel $z^*_i$ est le conjugué complexe de $z_i$.\\

Les opérateurs de mélange symétrisés sont des matrices semi-définies positives. En effet, pour chaque opérateur de mélange symétrisé $M$, il existe une matrice réelle $\pi$ telle que $\pi\pi^\top = M$~: à la \autoref{ssec:op_sym}, on disait justement que les opérateurs de mélange symétrisés s'écrivaient toujours comme $\pi^\top\circ\pi$ pour un certain opérateur (réel) de \bhr\ $\pi$. Comme on multiplie les matrices à droite (tel qu'expliqué dans la \autoref{rmq:mult_mat_a_droite}), ceci revient à dire que $M=\pi\pi^\top$. L'opérateur $\pi$ est décrit dans la prochaine sous-section. Ainsi, pour tout vecteur colonne $z$ à coefficients complexes,
\[ z^* Mz = z^*\pi\pi^\top z = (\pi^\top\overline{z})^\top(\pi^\top z) = \sum_{i=1}^{\dim(M)} (\overline{\pi^\top z})_i(\pi^\top z)_i  = \sum_{i=1}^{\dim(M)} \parallel(\pi^\top z)_i \parallel^2	\ \in \mathbb{R}_+.\]

\subsection{Opérateurs de $\bhr$ symétrisés}\label{ssec:1efamille_bhr}
On a énoncé plus tôt (\autoref{ssec:op_sym} et \autoref{ssec:def_matrices_nuk}) que les matrices des opérateurs $\allnuk$ pouvaient s'écrire comme la composition d'une matrice $\pi_k$ et de sa transposée ($\pi_k^\top$). En effet, chacun des $\nu_k$ consiste en deux grandes étapes~: retirer $k$ objets et les mélanger de toutes les façons possibles, puis les placer sur le dessus est la première grande étape; la seconde consiste à retirer les $k$ objets du dessus de notre collection, à les mélanger, puis à les battre avec le reste de la collection. 

Ainsi, la première grande étape est une généralisation de l'opérateur de l'aléatoire vers le dessus dans laquelle on retire $k$ objets plutôt qu'un seul; on note cet opérateur $\pi_k$\index{Opérateurs!$\pi_k$}. Pour la notation ici, on suppose que le dessus est à la droite des mots. Pour une permutation $\sigma = \sigma_1 \ldots \sigma_n$, on obtient alors 
\begin{equation}
	\pi_k(\sigma) = \sum_{1\leq i_1< \ldots i_k\leq n}\sum_{\tau \in S_k}\sigma_1 \ldots \widehat{\sigma_{i_1}} \ldots \widehat{\sigma_{i_k}} \ldots \sigma_{n} \sigma_{i_{\tau_1}}\ldots \sigma_{i_{\tau_k}}. \label{eq:pi_k_sigma}
\end{equation}
On peut voir sans trop de difficulté que $\pi_k$ est un opérateur de \bhr, puisqu'on l'obtient en multipliant $\sigma$ à droite par une certaine somme de composition. En effet, le produit de $\sigma$ par la somme de toutes les compositions de $k+1$ blocs, le premier étant de taille $n-k$ et suivi de $k$ singletons, donne exactement le résultat de l'équation \eqref{eq:pi_k_sigma}.

\begin{ex}
	L'opérateur $\pi_2$ appliqué à \textsf{1234} donne le résultat \begin{align*} 
		\pi_2(\mathsf{1234}) = \mathsf{1234}& + \mathsf{1243}  + \mathsf{1324} + \mathsf{1342} + \mathsf{1423} + \mathsf{1432}\\ & + \mathsf{2314} + \mathsf{2341} + \mathsf{2413} + \mathsf{2431} + \mathsf{3412} + \mathsf{3421}.
	\end{align*}
	C'est aussi le résultat de $\sigma \cdot \left(\sum_{c_1, c_2} [\{1,2,3,4\}\backslash\{c_1, c_2\},\{c_1\}, \{c_2\}]\right)$.
\end{ex}

La deuxième grande étape des opérateurs de mélange $\nu_k$ est le battage des $k$ objets du dessus avec les autres objets. Ceci n'est pas un opérateur de \bhr, puisque l'action effectuée ne dépend pas des objets à déplacer, mais plutôt de leur position. On fait toutefois le mouvement inverse de celui à la première étape~: plutôt que de déplacer $k$ objets de n'importe quelle position vers le dessus, on déplace les $k$ objets du dessus vers n'importe quelle position. C'est d'ailleurs une généralisation d'un mélange bien connu, celui du dessus vers l'aléatoire. Ainsi, le nombre de façons d'obtenir $\sigma$ à partir de $\tau$ à la première étape correspond au nombre de façons d'obtenir $\tau$ à partir de $\sigma$ à la seconde étape.
La matrice de la seconde étape est donnée par $\pi_k^\top$.

Enfin, il est utile de remarquer que, dans la définition de chacune des deux étapes, nous avons dit que nous mélangions les $k$ objets retirés. Dans les faits, nous n'avons besoin de le faire qu'une fois, plutôt que deux (à la fin de la première étape et au début de la deuxième). En mélangeant deux fois, nous n'avons pas changé les résultats possibles. Nous n'avons qu'artificiellement introduit $k!$ fois plus de résultats possibles, et nous devons donc diviser par ce nombre pour obtenir une définition adéquate de $\nu_k$.

On trouve donc la définition suivante pour l'opérateur $\nu_k$~:

\begin{defn}\label{def:nu_k_bhr}
	Soit $\pi_k$ la multiplication à droite par la somme des compositions de la forme $(n-k,1, \ldots, 1)$. Alors, $\nu_k$ est la multiplication à droite par $\frac{\pi_k\pi_k^\top}{k!}$.
\end{defn}
Notons que la dernière définition ne s'applique qu'aux permutations, et non aux mots.

\section{Commutativité}
\begin{thm}[Théorème I.1.1 de \cite{RSW}]\label{thm:commutativite_nu}
	Les opérateurs de la famille $\allnuk$ commutent entre eux.
\end{thm}

Victor Reiner, Franco Saliola et Volkmar Welker remarquaient, dans leur article, que cette propriété était très particulière à la famille $\allnuk$. En effet, parmi les opérateurs de mélange symétrisés, seules les familles $\allnuk$ et $\allgammak$ (présentée au \autoref{chap:2efamille}) voient leurs opérateurs commuter entre eux (à l'exception de cas triviaux). Dans leur étude de la famille $\allnuk$, ils notaient toutefois un défi~: trouver une preuve plus éclairante que la leur du \autoref{thm:commutativite_nu}. Une telle preuve se trouve au \autoref{chap:preuves}.

\section{Valeurs propres}
\label{sec:eval}

Dans leur article sur les opérateurs de mélange symétrisés, Reiner, Saliola et Welker ont largement discuté des valeurs propres des opérateurs de mélange symétrisés, ce qui a permis d'en dégager quelques propriétés.

Ils ont d'abord démontré que les valeurs propres de tous les opérateurs de mélange symétrisés étaient réelles et positives. Dans la \autoref{ssec:positivite}, nous en détaillons les raisons.

Ils ont aussi remarqué que les valeurs propres semblaient être entières pour les opérateurs $\allnuk$ et $\allgammak$. Dans le second cas, une formule utilisant les caractères du groupe symétrique existe pour calculer les valeurs propres; ceci est d'ailleurs expliqué à la \autoref{ssec:gammak_caracteres}. Ils ont de plus conjecturé qu'il s'agit des deux seules familles d'opérateurs de mélange symétrisés dont les valeurs propres sont entières.

Dans cette section, on démontre que les valeurs propres de $\Nu{k}$ sont bel et bien entières. Pour les valeurs propres de $\gammak{k}$, elles sont détaillées au \autoref{chap:2efamille}. Quant au fait que les valeurs propres des autres opérateurs ne sont pas entières, ce sujet n'est pas abordé ici.

En plus de démontrer que les valeurs propres de $\Nu{k}$ sont entières, on donne une formule récursive permettant de les calculer explicitement, c'est-à-dire qu'on peut calculer les valeurs propres de $\Nu{k}$ à partir de celles de $\Nu{k-1}$. Notons que les valeurs propres pour le mélange doublement aléatoire, $\Nu{1}$, ont déjà été calculées par Anton Dieker et Franco Saliola \cite{DS}; leur résultat est présenté à la \autoref{ssec:r2r}. La façon de faire combinatoire présentée ici permet notamment de calculer les valeurs propres quand les suites de lettres sont longues~: la complexité du calcul est en effet donnée plus loin (\autoref{ssec:consequences}).

\subsection{Positivité} \label{ssec:positivite}
Dans leur étude des opérateurs $\allnuk$, Reiner, Saliola et Welker ont remarqué que leurs valeurs propres étaient toujours réelles et positives. Le coeur de leur argument est que chaque matrice $\Nu{k}$ est réelle, symétrique et semi-définie positive. 

En effet, le théorème spectral pour les matrices indique que toute matrice réelle et symétrique $M$ est diagonalisable par une matrice unitaire $P$ (c'est-à-dire que \mbox{$P^\top = P^{-1}$}; voir, par exemple, \cite[théorème 2.5.6]{HJ}). Les valeurs propres de $M$ sont les éléments sur la diagonale de $P^\top MP$ et, comme $M$ est semi-définie positive, elles sont des nombres réels et positifs.

\subsection[Valeurs propres et tableaux]{Valeurs propres et tableaux~: statistiques importantes}

Les notions suivantes sont essentielles pour exprimer les valeurs propres.
\paragraph{Diagramme gauche}
Considérons deux diagrammes $\lambda$ et $\mu$ de sorte que toutes les boîtes de $\mu$ soient aussi dans $\lambda$. On définit alors le \textit{diagramme gauche} $\lambda/\mu$\index{$\lambda/\mu$} comme l'ensemble des boîtes qui sont dans $\lambda$ sans être dans $\mu$. Comme dessiné dans l'\autoref{ex:diagramme_gauche}, ces boîtes ne sont pas nécessairement justifiées à gauche.

\begin{ex}\label{ex:diagramme_gauche}
	Le tuple $\lambda = (4,3,2,1)$ est un partage de $10$, alors que $\mu = (2,1,1)$ est un partage de $4$. La troisième illustration est le diagramme gauche $\lambda/\mu$, et il contient $10-4 = 6$ boîtes. 
	\[ \YFrench \Yboxdim{16pt} \yng(4,3,2,1)\ , \qquad \yng(2,1,1)\ ,\qquad 
	\gyoung(::;;,:;;,:;,;) \ .\]
\end{ex}

\paragraph{Indice diagonal}
Une statistique sur les diagrammes (incluant les diagrammes gauches) qui sera utile pour nos calculs est son \textit{indice diagonal}. À chaque boîte dans le diagramme, on associe un nombre en prenant ses coordonnées $(i,j)$, comme dans un plan cartésien, et on calcule la différence $j-i$ (voir l'\autoref{ex:diagonal_index}). Pour un diagramme (gauche) $\lambda$, son indice diagonal est la somme des indices diagonaux de ses boîtes et est noté $\diag(\lambda)$\index{$\diag(\lambda)$}. 
\begin{ex} \label{ex:diagonal_index}
	L'indice diagonal des boîtes dans les diagrammes ci-dessous est inscrit dans chacune des boîtes. Le diagramme de gauche a pour indice diagonal $7$, alors que celui de droite a un indice diagonal de $17$.
	\[\YFrench \Yboxdim{16pt} \yngres(0,5,4,4,1)\ , \qquad  \yngres(0,6,4)\ .\]
\end{ex}
\begin{rmq}
	Certains textes font référence à l'indice diagonal en termes de \textit{contenu} d'une boîte. Toutefois, le contenu d'un tableau peut aussi qualifier les nombres que l'on trouve dans les boîtes. Pour éviter toute ambiguïté, on n'utilise le mot contenu que dans le dernier cas, et cela est détaillé à la \autoref{ssec:vp_mots}.
\end{rmq}

\subsubsection{Opérateur $\Delta$ de Schützenberger}\label{sssec:op_delta_schutz}
La stratégie présentée pour calculer les valeurs propres est une induction à la fois sur l'indice de l'opérateur, $k$, et sur le nombre d'objets dans la collection, $n$. L'idée de base est qu'il n'est pas difficile de calculer les valeurs propres des matrices pour de très petites collections. Par contre, les dimensions augmentent de façon factorielle, et calculer les valeurs propres de façon naïve paraît rapidement complètement irréaliste.

Pour calculer une valeur propre pour un opérateur $\Nu{k}$ et une collection de taille $n$ donnés, on procède ainsi~: on attribue à chaque tableau de taille $n$ un autre tableau, cette fois de taille $n-1$. On fait donc le pari que la valeur propre associée au tableau de taille $n-1$  sera plus facile à calculer. En faisant ce processus, on fait donc la restriction de $\CSn$ à $\C S_{n-1}$. Avec la règle de branchement, présentée à la \autoref{ssec:regle_branchement}, on sait que le plus petit tableau a la même forme que le plus grand, auquel on aurait enlevé une boîte.

Pour trouver une façon d'associer à chaque tableau de taille $n$ le bon tableau de taille $n-1$, c'est-à-dire celui qui nous permet de calculer les valeurs propres, on utilise l'\textit{opérateur $\Delta$ de Schützenberger} \cite{schutzenberger}\index{Opérateurs!$\Delta$}. Pour l'exécuter sur un tableau $t$ de taille $n$~:
\begin{enumerate}
	\item On retire l'élément $1$ et on le remplace par une boîte vide.
	\item Tant que la boîte vide n'est pas dans un coin extérieur du tableau, on la déplace vers l'extérieur en faisant des mouvements de jeu-de-taquin~: ceci veut dire que, si $(i,j)$ est la position de la boîte vide, on l'échange avec la boîte qui contient le minimum des boîtes aux positions $(i+1,j)$ et $(i,j+1)$.
	\item Une fois que la boîte vide est dans un coin extérieur, on obtient un nouveau tableau, sans boîte vide, qui contient les nombres $2, \ldots, n$. Pour retrouver un tableau standard, on remplace ces nombres par $1, \ldots, n-1$, en préservant l'ordre, et on enlève la boîte vide.
\end{enumerate}
Un exemple de ce processus apparaît à l'\autoref{ex:jdt}.

\begin{ex}\label{ex:jdt}
	Appliquer l'opérateur $\Delta$ de Schützenberger au tableau standard $t$ de taille $7$ retourne un tableau standard de taille $6$.
	\[\YFrench t = \Yboxdim{14pt}\young(125,347,6) \to 
	\young(\bullet25,347,6) \to \young(2\bullet5,347,6) \to 
	\young(245,3\bullet7,6) \to \young(245,37\bullet,6) \to 
	\young(134,26,5) = \Delta(t). \]
\end{ex}

Ce processus nous permet de déduire un unique tableau standard de taille $n-1$ de chaque tableau standard de taille $n$. Ce faisant, on fait la restriction de $\CSn$ à $\C S_{n-1}$.

\begin{rmq}\label{rmq:inv_schutz_delta}
	Il existe aussi une procédure inverse à l'opérateur $\Delta$. Pour passer d'un tableau de taille $n$ à un tableau de taille $n+1$, on doit cependant spécifier dans quelle rangée on souhaite ajouter une boîte. On déplace ensuite la boîte vide vers la première ligne et la première colonne en faisant des mouvements de jeu-de-taquin à l'envers.
\end{rmq}

\subsubsection{Bandes horizontales et tableaux de désarrangement}

Dans un tableau standard de taille $n$, l'entrée $i$ est une \textit{montée} si $i=n$ ou si la boîte contenant $i+1$ est située au sud, à l'est ou au sud-est de celle contenant $i$.
\begin{ex}
	Dans les tableaux ci-dessous, $3$ est une montée dans les tableaux de gauche et de droite. Au centre, par contre, $3$ est une descente, puisque que $4$ se trouve au nord. Dans tous les cas, $4$ est une montée puisque c'est la plus grande entrée dans ces tableaux.
	\[\YFrench \gyoung(;;,34)\ ,\quad\gyoung(;;,3,4)\ ,\quad 
	\gyoung(;;;4,3)\ .\]
\end{ex}

\begin{ex}\label{ex:montees}
	Dans ce tableau de forme $(4,2,2,1)$, les montées sont identifiées par un arrière-plan de couleur.
	\newcommand{\ylw}{\Yfillcolour{yellow}}
	\newcommand{\wh}{\Yfillcolour{white}}
	\[
	\tab(!\ylw12!\wh3!\ylw9,!\wh46,!\ylw58,7)\]
\end{ex}

\begin{rmq}
	Les notions de montées et de descentes sur un tableau peuvent être un peu difficile à visualiser avec la notation française des tableaux. En anglais, la première ligne d'un tableau est située tout en haut, et les lignes suivantes sont placées plus bas, comme dans une matrice. Pour un tableau en anglais, une montée apparaît lorsque l'entrée suivante est située plus haut, ce qui correspond à notre intuition. Toutefois, les descentes d'un tableau standard sont liées aux descentes d'une permutation~: le tableau $Q(\sigma)$ obtenu de $\sigma$ par la correspondance de Robinson-Schensted a les mêmes descentes que $\sigma$ (pour des informations sur cet algorithme, voir, par exemple, \cite[section 3.1]{sagan} ou \cite[section 7.11]{EC2}). C'est pour cette raison que les noms descentes et montées d'un tableau sont préservées en français, malgré la confusion visuelle.
\end{rmq}

Un tableau standard est un \textit{tableau de désarrangement} s'il est vide ou si sa plus petite descente est à une position paire.\footnote{Le mot \textit{désarrangement} n'apparaît pas dans le dictionnaire. C'est plutôt une contraction du nom Désarménien et des mots dérangement et descentes. Michelle Wachs et Jacques Désarménien ont démontré que les permutations dont la plus petite montée est paire sont en bijection avec les permutations sans point fixe (les \textit{dérangements}) \cite{DW}. Comme les descentes d'une permutation sont les mêmes que les descentes du tableau $Q$ qui lui est associé par la correspondance de Robinson-Schensted, ils ont démontré un certain lien entre les tableaux dont la plus petite montée est paire et les dérangements.} 

\begin{ex}\label{ex:desar_tableaux}
	Parmi les tableaux de taille $2$ et $3$, il n'y a que deux tableaux de désarrangement~:
	\[\YFrench \begin{tikzpicture}
	\node (img) {\young(12)};
	\draw[red, line width=.2mm] 
	(img.south west) -- (img.north east)
	(img.south east) -- (img.north west);
	\end{tikzpicture}\ , \quad \young(1,2)\ , \quad 
	\begin{tikzpicture}
	\node (img) {\young(123)};
	\draw[red, line width=.2mm] 
	(img.south west) -- (img.north east)
	(img.south east) -- (img.north west);
	\end{tikzpicture}\ , \quad \begin{tikzpicture}
	\node (img) {\young(12,3)};
	\draw[red, line width=.2mm] 
	(img.south west) -- (img.north east)
	(img.south east) -- (img.north west);
	\end{tikzpicture} \ ,\quad \young(13,2)\ , \quad 
	\begin{tikzpicture}
	\node (img) {\young(1,2,3)};
	\draw[red, line width=.2mm] 
	(img.south west) -- (img.north east)
	(img.south east) -- (img.north west);
	\end{tikzpicture}.\]
\end{ex}

Un diagramme (ou un diagramme gauche) est une \textit{bande horizontale} s'il n'y a pas deux boîtes dans la même colonne (voir l'\autoref{ex:horiz_strip}).

\begin{ex}\label{ex:horiz_strip}
	Le diagramme gauche $(4,2,2,1)/(3,2,1)$ est une bande horizontale, comme le montre le diagramme à la gauche, dans lequel les boîtes qui forment la bande sont identifiées par une croix.
	\[\YFrench \gyoung(;;;;\times,;;,;;\times,\times)\]
\end{ex}

\begin{prop}[Proposition 6.25 de \cite{RSW}]\label{prop:desa_vs_bandes}
	À chaque tableau standard $t$, on peut associer un unique tableau de désarrangement en appliquant l'opérateur $\Delta$ jusqu'à ce que le résultat soit un tableau de désarrangement. Lorsqu'on doit appliquer  $j$ fois l'opérateur $\Delta$ pour obtenir un tableau de désarrangement, alors $t/\Delta^j(t)$ est une bande horizontale de taille $j$.
\end{prop}

Grâce à la \autoref{prop:desa_vs_bandes}, on peut définir le \textit{type} d'un tableau standard $t$ comme le plus petit nombre $j$ tel que $\Delta^j(t)$ est un tableau de désarrangement.\index{$\type(t)$}

\begin{ex}
	Le tableau de l'\autoref{ex:montees} a un type égal à $3$~:
	\Yvcentermath1
	\[\Yboxdim{14pt}\Delta^3\left(\tab(1239,46,58,7)\right) = \Delta^2\left(\tab(128,35,47,6)\right) = \Delta\left(\tab(147,26,3,5)\right) =\tab(136,25,4)\ .\]
	Le diagramme gauche de la forme $t/\Delta^3(t)$ est la bande horizontale de l'\autoref{ex:horiz_strip}~:
	\[\YFrench \gyoung(;;;;\times,;;,;;\times,\times)\ .\]
\end{ex}

\subsection{Le mélange doublement aléatoire}\label{ssec:r2r}
Anton Dieker et Franco Saliola ont récemment décrit les valeurs propres de l'opérateur de mélange doublement aléatoire \cite{DS}. Cet opérateur, tel que décrit à la \autoref{par:r2r_def}, est non seulement un mélange plutôt naturel, c'est aussi l'opérateur $\Nu{1}$, ce qui en fait le cas de base de notre induction pour les valeurs propres de $\allnuk$. Il est utile de rappeler leur description des valeurs propres~:
\begin{thm}[Théorème 5 de \cite{DS}]\label{thm:vp_r2r_ds}
	Chaque valeur propre de l'opérateur de mélange doublement aléatoire (en tant que chaîne de Markov) est de la forme $\frac{1}{n^2}\eig(\lambda/\mu)$, où $\lambda$ est un partage de  $n$ et $\lambda/\mu$ est une bande horizontale. De plus, $\eig$ est la statistique combinatoire définie sur les partages gauches comme
	\[\eig(\lambda/\mu) = \binom{|\lambda| + 1}{2} - \binom{|\mu| + 1}{2} + \diag(\lambda/\mu). \]
	
	La multiplicité de la valeur propre $\frac{1}{n^2}\varepsilon$ est
	\[ \sum_{\substack{ \lambda/\mu\text{ est une bande horixontale,}\\
			\eig(\lambda/\mu)=\varepsilon}} f^\lambda d^\mu,\]
	où $f^\lambda$ est le nombre de tableaux standards de forme $\lambda$, et $d^\mu$ est le nombre de tableaux de désarrangement de forme $\mu$.
\end{thm}

\begin{rmq}
	Dans la forme dans laquelle le théorème est formulée, les valeurs propres sont toutes réelles et se situent entre $0$ et $1$. Dans ce contexte, $\frac{1}{n^2}\eig(\lambda/\mu)$ n'est pas un entier, mais plutôt la valeur propre de la chaîne de Markov du mélange doublement aléatoire. La valeur propre entière de $\Nu{1}$ est $\eig(\lambda/\mu)$.
\end{rmq}

\begin{ex}
	Considérons les valeurs propres du mélange doublement aléatoire sur $4$ éléments. La matrice de cet opérateur est carrée et de taille $24$; on peut donc calculer en temps raisonnable ses valeurs propres avec un logiciel comme SageMath \cite{sagemath}.

	On obtient les valeurs propres affichées au \autoref{tab:vp_R2R4_sage}.
	\begin{table}
		\caption{Valeurs propres du mélange doublement aléatoire avec quatre objets distincts.}\label{tab:vp_R2R4_sage}
		\centering
		\begin{tabular}{|l|c|c|c|c|c|c|}
			\hline
			Valeurs propres & 16 & 10 & 6 & 4 & 2 & 0\\
			\hline
			Multiplicité & 1 & 3 & 6 & 2 & 3 & 9\\
			\hline
		\end{tabular}
	\end{table}
	 On peut aussi les calculer avec le \autoref{thm:vp_r2r_ds}, comme au \autoref{tab:vp_R2R4_thm}. Naturellement, on obtient les mêmes valeurs propres qu'au \autoref{tab:vp_R2R4_sage}.
	\begin{table}\caption[Calcul des valeurs propres pour les collections de quatre objets distincts.]{Calcul des valeurs propres pour les collections de quatre objets distincts. Les dernières colonnes n'ont pas été remplies s'il n'y a pas de valeur propre de toute façon (si $d^\mu = 0$ ou si $\lambda/\mu$ n'est pas un tableau de désarrangement). Ici, $\binom{|\lambda|+1}{2} = 10$.}\label{tab:vp_R2R4_thm}
		\centering
		\begin{math}
		\Yvcentermath1
		\arraycolsep=6pt\def\arraystretch{1.1}
		\begin{array}{ccccccc}
		\lambda & \mu \text{ si } d^\mu\neq 0 & \binom{|\mu|+1}{2} & d^\mu 
		& \lambda/\mu \text{ est b.h.} & f^\lambda & 
		\eig(\lambda/\mu)\\
		\hline          
		\smalldia(4) & \emptyset & 0 & 1 & \checkmark & 
		1 & 16 \\
		\hline
		\smalldia(3,1) & \emptyset & 0 & 1 & \times\\
		& \smalldia(1,1) & 3 & 1  & \checkmark & 3 & 10\\
		& \smalldia(2,1) & 6 & 1 & \checkmark & 3 & 6\\
		& \smalldia(3,1) & 10 & 1 & \checkmark & 3 & 0\\
		\hline
		\smalldia(2,2) & \emptyset & 0 & 1 & \times \\
		& \smalldia(1,1) & 3 & 1  & \times\\
		& \smalldia(2,1) & 6 & 1 & \checkmark & 2 & 4\\
		& \smalldia(2,2) & 10 & 1 & \checkmark & 2 & 0\\
		\hline
		\smalldia(2,1,1) & \emptyset & 0 & 1 & \times \\
		& \smalldia(1,1) & 3 & 1  & \checkmark & 3 & 6\\
		& \smalldia(2,1) & 6 & 1 & \checkmark & 3 & 2\\
		& \smalldia(2,1,1) & 10 & 1 & \checkmark & 3 & 0\\
		\hline
		\smalldia(1,1,1,1) & \emptyset & 0 & 1 & \times \\
		& \smalldia(1,1) & 3 & 1  & \times\\
		& \smalldia(1,1,1,1) & 10 & 1 & \checkmark & 1 & 0
		\end{array}
		\end{math}
	\end{table}
\end{ex}

À la \autoref{ssec:defn_Specht}, on a mentionné qu'on pouvait associer à chaque tableau standard une valeur propre. Dans le cas du mélange doublement aléatoire, on peut associer une paire de partages dont la différence est une bande horizontale avec une valeur propre. On aimerait donc trouver une façon d'associer cette paire de partages avec un tableau standard. Or, la \autoref{prop:desa_vs_bandes} nous donne un indice sur comment faire ceci~: à partir d'un tableau standard $t$, on applique $\type(t)$ fois l'opérateur $\Delta$, soit jusqu'à ce qu'on obtienne un tableau de désarrangement. Alors, $t/\Delta^{\type(t)}(t)$ est une bande horizontale, et on peut calculer la valeur propre $\eig(t/\Delta^{\type(t)}t)$. On peut donc associer à un tableau standard une unique paire de diagrammes dont la différence est une bande horizontale, et donc une valeur propre du mélange doublement aléatoire. Cette reformulation sera utilisée dans la \autoref{ssec:thm_principal} pour calculer les valeurs propres des autres opérateurs. Pour un tableau standard $t$~:
\begin{itemize}
	\item Si $t$ est un tableau de désarrangement, la valeur propre associée à $t$ est $0$.
	\item Sinon, la valeur propre est $\varepsilon + n + \diag(t/\Delta(t))$, où $\varepsilon$ est la valeur propre associée à $\Delta(t)$.
\end{itemize}
Puisque tous les tableaux standards mènent à un tableau de désarrangement en itérant successivement $\Delta$ (\autoref{prop:desa_vs_bandes}), ce processus termine toujours.

\section{Une formule récursive pour les valeurs propres}
On est maintenant en mesure de présenter le théorème principal concernant les opérateurs $\allnuk$, soit une description de toutes leurs valeurs propres. Pour le mélange d'une collection d'objets tous distincts, on présente le résultat à la \autoref{ssec:thm_principal}; le cas des collections avec objets répétés est présenté à la \autoref{ssec:vp_mots}. Quant aux preuves des principaux théorèmes présentés ici (\autoref{thm:main} et \autoref{thm:val_pr_mots}), elles se trouvent au \autoref{chap:preuves}.
\subsection{Le cas particulier des permutations}\label{ssec:thm_principal}
Pour les collections d'objets tous distincts, on procède en deux temps~: on rappelle d'abord un résultat de Reiner, Saliola et Welker concernant le noyau des opérateurs $\allnuk$, puis on décrit les autres valeurs propres.

\begin{thm}[Théorème VI.9.5 de \cite{RSW}]\label{thm:ker_nuk_rsw}
	Le noyau de $\nu_k$ est isomorphe à une somme directe $\bigoplus S^{\mathrm{forme}(t)}$, dans laquelle $S^\lambda$ représente le module de Specht associé aux partages de forme $\lambda$, et la somme directe se fait sur tous les tableaux standards de type strictement inférieur à $k$.
\end{thm} 
\begin{ex}\label{ex:syt4_noyaux}
	\begin{table}
		\caption{Les tableaux standards de taille $4$ et leur type.}\label{tab:syt(4)_type}
		\[\Yvcentermath1 \YFrench \Yboxdim{11pt}
		{\setlength\arraycolsep{3.5pt}
			\begin{array}{|r|c|c|c|c|c|c|c|c|c|c|} 
			\hline
			&\young(1234) & \young(123,4) & \young(124,3) & \young(134,2) & \young(12,34) & \young(13,24) & \young(12,3,4) &\young(13,2,4) & \young(14,2,3) & \young(1,2,3,4)\\
			\hline
			\textrm{{\small Type}} & 
			4 & 2 & 1 & 0 & 1 & 0 & 2 & 0 & 1 & 0\\
			\hline
			\end{array}}\]
	\end{table}
	Les tableaux standards de taille $4$ sont dans le \autoref{tab:syt(4)_type}, avec leur type. Les noyaux des opérateurs $\allnuk$ pour une collection de quatre objets distincts sont les suivants~:
	\[\begin{array}{|c|c|}
		 \hline
		\Nu{1} & S^{\smalldia(3,1)} \oplus S^{\smalldia(2,2)} \oplus S^{\smalldia(2,1,1)} \oplus S^{\smalldia(1,1,1,1)}\\
		\hline
		\Nu{2} & 2 S^{\smalldia(3,1)} \oplus 2 S^{\smalldia(2,2)} \oplus 2 S^{\smalldia(2,1,1)} \oplus S^{\smalldia(1,1,1,1)}\\
		\hline
		\Nu{3} & \C S_4 \backslash S^{\smalldia(4)}\\
		\hline
		\Nu{4} & \C S_4 \backslash S^{\smalldia(4)}\\
		\hline
		\Nu{\geq 5} & \C S_4\\
		\hline
	\end{array}	\]
\end{ex}

Ce qu'on doit retenir du théorème \autoref{thm:ker_nuk_rsw}, c'est qu'à chaque tableau de type au moins $k$, on doit associer une valeur propre non-nulle pour l'opérateur $\Nu{k}$. La façon dont nous calculons cette valeur nous permet de montrer que c'est un entier positif (\autoref{corl:vp_nuk_entiers}).

\begin{thm}\label{thm:main}
	Les valeurs propres de l'opérateur de mélange $\Nu{k}$ sur des collections de $n$ objets tous distincts sont indexées par les tableaux standards de taille $n$.
	Pour un tel tableau $t$, la valeur propre $v_k(t)$ est donnée par
	\[v_k(t) = \begin{cases}
	v_k(\Delta(t)) + \left(n+1-k+\diag(t/\Delta(t))\right)\ 
	v_{k-1}(\Delta(t)) & \text{ si } \type(t) \geq k,\\
	0 & \text{sinon}.
	\end{cases}\]
	Le tableau $t$ contribue pour $f^{\mathrm{forme}(t)}$ à la multiplicité de la valeur propre $v_k(t)$, où $f^{\mathrm{forme}(t)}$ est le nombre de tableaux standards de la même forme que $t$.
\end{thm}

Remarquons que si on pose $k=1$, on retrouve les valeurs propres calculées pour le mélange doublement aléatoire ($\Nu{1}$) à la  \autoref{ssec:r2r}.

\begin{ex}\label{ex:main}
	Avec le \autoref{thm:main}, on peut par exemple calculer les valeurs propres pour $\Nu{2}$ sur une collection de quatre éléments distincts. Il y a alors $4! = 24$ états possible (les permutations des quatre éléments); on peut donc encore calculer les valeurs propres de la matrice avec les techniques habituelles. C'est ce qu'on a fait avec l'aide de SageMath \cite{sagemath}, et on sait donc que les valeurs propres sont $0$, $4$, $20$ et $72$, de multiplicités $17$, $3$, $3$ et $1$, respectivement. Ces valeurs sont confirmées par le \autoref{thm:main}, comme nous le démontrons ici.\\
	Il y a dix tableaux standards de taille $4$; ils sont présentés au \autoref{tab:syt(4)_type} avec leur type. 
	
	\Yvcentermath1  
	On peut déjà calculer la dimension du noyau, et donc la multiplicité de la valeur propre $0$~: à l'\autoref{ex:syt4_noyaux}, nous avons calculé que le noyau de $\Nu{2}$ est \mbox{$2 S^{\smalldia(3,1)} \oplus 2 S^{\smalldia(2,2)} \oplus 2 S^{\smalldia(2,1,1)} \oplus S^{\smalldia(1,1,1,1)}$.} Nous savons aussi que la dimension du module $S^\lambda$ est égale au nombre de tableaux standards de forme $\lambda$. Ainsi, la dimension du noyau est 
	\[2 f^{\smalldia(3,1)} + 2 f^{\smalldia(2,2)} + 2 f^{\smalldia(2,1,1)} + f^{\smalldia(1,1,1,1)} = 2 \cdot 3 + 2 \cdot 2 + 2\cdot 3 + 1 =17. \]

	Pour trouver les valeurs propres non-nulles de $\Nu{2}$, on doit calculer $v_k(t)$ pour les tableaux standards de type au moins $2$. Il y a trois tels tableaux~:
	\begin{align*}
	\Yboxdim{10pt}
	v_2\left(\tab(1234)\right)
	&= \Yboxdim{10pt}v_2\left(\tab(123)\right) + (4+1-2+3)\cdot 
	v_1\left(\tab(123)\right)\\
	&= \Yboxdim{10pt} \left(  v_2\left(\tab(12)\right) + (3+1-2+2)\cdot 
	v_1\left(\tab(12)\right) \right) + 6 \cdot 
	v_1\left(\tab(123)\right)\\
	&= \Yboxdim{10pt} \left(\left(v_2\left(\tab(1)\right) + (2+1-2+1) 
	\cdot v_1\left(\tab(1)\right) \right) + 4\cdot 
	v_1\left(\tab(12)\right) 
	\right) + 6 \cdot 9\\
	& = \Yboxdim{10pt} \left(\left( 0 + 2 \cdot v_1\left(\tab(1)\right)\right)+ 4 \cdot 4 \right) + 54\\
	& = 2 \cdot 1 + 16+ 54 = 72.\\
	\Yboxdim{10pt} v_2\left(\tab(123,4)\right)
	&= \Yboxdim{10pt} v_2\left(\tab(12,3)\right) + (4+1-2+2)\cdot 	v_1\left(\tab(12,3)\right)\\
	&= \Yboxdim{10pt} 0 + 5 \cdot v_1\left(\tab(12,3)\right)\\
	&= 0 + 5 \cdot 4 = 20.\\
	\Yboxdim{10pt} v_2\left(\tab(12,3,4)\right) &= \Yboxdim{10pt} v_2\left(\tab(1,2,3)\right) +(4+1-2+1) \cdot v_1\left(\tab(1,2,3)\right)\\
	& \Yboxdim{10pt} = 0 + 4 \cdot v_1\left(\tab(1,2,3)\right)\\
	& = 4 \cdot 1 = 4.
	\end{align*}
	
	Dans le dernier calcul, la deuxième égalité vient du fait que  $\tab(1,2,3)$ est de type strictement inférieur à 2, tandis que la troisième a été calculée en utilisant le résultat de la \autoref{ssec:r2r}, concernant les valeurs propres du mélange doublement aléatoire. On aurait aussi pu faire ce dernier calcul avec le \autoref{thm:main}, en utilisant le fait que $v_0(t) = 1$ pour tout tableau $t$, étant donnée que $v_0$ est l'identité.
\end{ex}

\subsection{Sur les mots en général}\label{ssec:vp_mots}
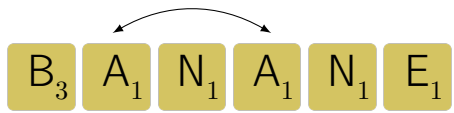
\begin{figure}
	\centering
	\begin{tikzpicture}
		\node (B) at (0,0) {\Scrabble{B}{3}};
		\node (A1) at (1,0) {\Scrabble{A}{1}};
		\node (gauche) at (0.8,0.5) {};
		\node (N1) at (2,0) {\Scrabble{N}{1}};
		\node (A2) at (3,0) {\Scrabble{A}{1}};
		\node (droite) at (3.2,0.5) {};
		\node (N2) at (4,0) {\Scrabble{N}{1}};
		\node (E) at (5,0) {\Scrabble{E}{1}};
		\draw[<->, >=latex] (droite) to [out=150, in=30] (gauche) {};
	\end{tikzpicture}
	\caption[Échanger les deux occurrences de $a$ dans le mot \textit{banane} redonne le mot initial.]{Échanger les deux occurrences de $a$ dans le mot \textit{banane} redonne le mot initial. C'est une opération qu'on effectue, par exemple, au cours d'une partie de Scrabble.}\label{fig:scrabble}
\end{figure}

Jusqu'ici, nous nous avons présenté les résultats sur les valeurs propres pour des suites d'objets que nous avons supposés tous distincts. Par exemple, un paquet de cartes à jouer, s'il n'est pas truqué, est un ensemble d'objets où toutes les cartes sont différentes. Par contre, les techniques énoncées plus haut s'appliquent aussi aux collections qui admettent des répétitions d'objets. Ainsi, on pourrait mélanger les lettres du mot \textit{banane} (dans un jeu de Scrabble, par exemple), et on souhaiterait que le mot initial et celui obtenu en échangeant les deux \textit{a} soient considérés comme le même mot. Ceci est illustré à la \autoref{fig:scrabble}.

Tel qu'expliqué à la \autoref{ssec:def_matrices_nuk}, et illustré à l'\autoref{ex:matrice_mots}, les états de ces chaînes de Markov sont moins nombreux, et correspondent à tous les mots qu'il est possible de faire tout en gardant le nombre d'occurrences de chaque lettre constant.

On appelle le \textit{contenu} d'une suite d'objets le tuple formé du nombre de copies de chacun des objets dans la suite. Ces nombres sont placés en ordre décroissant, ce qui fait qu'on obtient un partage. Ainsi, nous pouvons le comparer aux autres partages.

\begin{ex}
	Le contenu du mot de langue française \textit{ananas} est $(3,2,1)$, alors que celui du mot \textit{banane} est $(2,2,1,1)$.
\end{ex}
\begin{ex}
	Les permutations sont des mots de contenu $(1,1,\ldots,1)$.
\end{ex}

\paragraph{Ordre de dominance}
Si $\lambda = (\lambda_1, \ldots, \lambda_l)$ et $\mu = (\mu_1, \ldots,  \mu_m)$ sont deux partages du même nombre, on dit que $\lambda$ \textit{domine} $\mu$, noté $\lambda \unrhd \mu$\index{$\unrhd$}, si, pour tout $i \in \{1, 2, \ldots, \min(l,m)\}$,
\[ \lambda_1 + \ldots + \lambda_i \geq \mu_1 + \ldots + \mu_i. \]

\begin{ex}
	L'ordre de dominance des partages de $6$ s'illustre par le diagramme de Hasse à la \autoref{fig:diag_hasse_dominance_6}.
	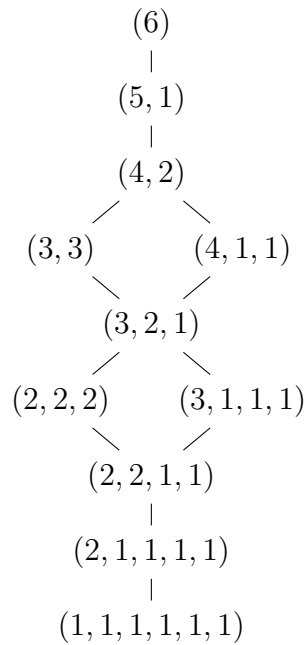
\begin{figure}
		\centering
		\begin{tikzpicture}
			\node (6) at (0,0) {$(6)$};
			\node[below of=6] (51) {$(5,1)$};
			\node[below of=51] (42) {$(4,2)$};
			\node[below left of=42,xshift=-5mm,yshift=-3mm] (33) {$(3,3)$};
			\node[below right of=42,xshift=5mm,yshift=-3mm] (411) {$(4,1,1)$};
			\node[below left of=411,xshift=-5mm,yshift=-3mm] (321) {$(3,2,1)$};
			\node[below left of=321,xshift=-5mm,yshift=-3mm] (222) {$(2,2,2)$};
			\node[below right of=321,xshift=5mm,yshift=-3mm] (3111) {$(3,1,1,1)$};
			\node[below left of=3111,xshift=-5mm,yshift=-3mm] (2211) {$(2,2,1,1)$};
			\node[below of=2211] (21111) {$(2,1,1,1,1)$};
			\node[below of=21111] (16) {$(1,1,1,1,1,1)$};
			\draw (6) -- (51);
			\draw (42) -- (51);
			\draw (33) -- (42);
			\draw (411) -- (42);
			\draw (321) -- (411);
			\draw (321) -- (33);
			\draw (321) -- (222);
			\draw (321) -- (3111);
			\draw (2211) -- (3111);
			\draw (2211) -- (222);
			\draw (2211) -- (21111);
			\draw (21111) -- (16);
		\end{tikzpicture}
		\caption{Ordre de dominance sur les partages de 6.}\label{fig:diag_hasse_dominance_6}
	\end{figure}
\end{ex}

\paragraph{Tableaux semi-standards}\label{par:SSYT}
Un tableau semi-standard est d'abord un tableau et c'est donc un remplissage d'un diagramme par des nombres entiers positifs. La notion de tableau semi-standard est plus flexible que celle de tableau standard, puisqu'elle autorise les répétitions dans les entrées du tableau.
\begin{defn}
	Un tableau est dit \textit{semi-standard} si les boîtes sur une même ligne contiennent des nombres placés en ordre faiblement croissant et si les boîtes d'une même colonne contiennent des nombres tous distincts placés en ordre croissant.
	
	Le \textit{contenu} d'un tableau standard est un partage qui compte le nombre d'occurrences de chacune des valeurs des entrées du tableau.
\end{defn}

\begin{ex}
	Le tableau suivant, de forme $(3,2)$, est semi-standard. Son contenu est $(2,2,1)$.
	\[ \tab(112,23) \]
\end{ex}
\begin{ex}
	Les tableaux standards sont des tableaux semi-standards de contenu $(1,1,1,\ldots, 1)$.
\end{ex}

\begin{table}
	\caption{Règle pour les tableaux standards et semi-standards.}
	\centering
	\begin{tabular}{cc}
		Tableau standard & Tableau semi-standard\\	
		$\Yvcentermath1 \begin{array}{rc}
		 & \stackrel{<}{\longrightarrow}\\
		 <\bigg\uparrow& \dia(3,2,1,1)
		\end{array}$ &
		$\Yvcentermath1 \begin{array}{rc}
		& \stackrel{\leq}{\longrightarrow}\\
		<\bigg\uparrow& \dia(3,2,1,1)
		\end{array}$
	\end{tabular}
\end{table}
Nous connaissons maintenant toutes les notions pour énoncer l'analogue du \autoref{thm:main}, pour les mots. Comme les mots sont une généralisation des permutations, le \autoref{thm:main} est un cas particulier de ce qui suit.
\begin{thm}\label{thm:val_pr_mots}
	Les valeurs propres non-nulles de $\Nu{k}$ sur les mots de contenu $\mu \vdash n$ sont indexées par les tableaux standards de taille $n$, de type au moins $k$ et de forme $\lambda \unrhd \mu$. Pour un tel tableau $t$,
	\[v_k(t) = v_k(\Delta(t)) + (n+1-k+\diag(t/\Delta(t)))\ 
	v_{k-1}(\Delta(t)).\]
	De plus, le tableau $t$ contribue $m_{\lambda\mu}$\index{$m_{\lambda\mu}$}, le nombre de tableaux semi-standards de forme $\lambda$ et remplis avec le contenu $\mu$,   à la multiplicité de la valeur propre $v_k(t)$. Ces nombres sont aussi appelés nombres de Kostka.
\end{thm}
Enfin, cela est cohérent avec le \autoref{thm:main}, puisque les permutations ont pour contenu $(1,1,\ldots,1)$, ce qui est minimal pour l'ordre de dominance.

\begin{ex}
	La matrice pour $\Nu{2}$ sur les mots de contenu $(2,2)$, comme $\textsf{aabb}$, est celle exhibée à l'\autoref{ex:matrice_mots}. Elle a pour valeurs propres $72$, $20$ et $0$ avec multiplicité $1$, $1$ et $4$, respectivement, ce qui est cohérent avec le \autoref{thm:val_pr_mots}, comme on peut le voir au \autoref{tab:vp_nu2_ssyt}. Ici, le contenu des mots étant $(2,2)$, les partages qui le domine sont $(4)$, $(3,1)$ et $(2,2)$. Les valeurs propres sont celles qui ont été calculées à l'\autoref{ex:main}, mais leur multiplicité doit être de nouveau calculées. C'est ce à quoi sert le \autoref{tab:vp_nu2_ssyt}.
	\begin{table}
		\caption{Calcul des valeurs propres pour l'opérateur $\Nu{2}$ sur une collection de taille $4$ contenant deux paires d'éléments identiques.}\label{tab:vp_nu2_ssyt}
		\[\Yvcentermath1 \YFrench
	{\setlength\arraycolsep{5pt}
		\renewcommand{\arraystretch}{1.9}
		\begin{array}{cccc}
		t \text{ \small de forme $\lambda$ }&\stackunder[3mm]{\text{Tableaux semi-standards}}{\text{de forme $\lambda$ et contenu $(2,2)$}} & v_2(t) & m_{\lambda,(2,2)}\\
		\young(1234) & \tab(1122) & 72 & 1\\
		\hdashline
		\young(123,4) & \tab(112,2) & 20 & 1\\
		\young(124,3) & \tab(112,2) & 0 & 1\\
		\young(134,2) & \tab(112,2) & 0 & 1\\
		\hdashline
		\young(12,34) & \tab(11,22) & 0 & 1\\
		\young(13,24) & \tab(11,22) & 0 & 1\\
		\end{array}}\]
	\end{table}
\end{ex}

Le \autoref{thm:val_pr_mots} découle du \autoref{thm:main}, puisque les valeurs propres sont les mêmes, mais les multiplicités des valeurs propres diffèrent. Celles-ci peuvent cependant être retrouvées avec la règle de Young pour les modules de permutation (\autoref{thm:regle_de_young_mod_perm}).

\subsection{Valeurs propres pour la composante isotypique \mbox{$(n-1,1)$}}\label{ssec:vp_Sn-1,1}
Aussi intéressant que soit le résultat de cette section, il n'est pas le plus pratique. Pourquoi? Parce que la formule récursive fait qu'il est difficile d'obtenir des bornes inférieures ou supérieures sur les valeurs propres, ou encore de faire la somme de toutes les valeurs propres. On ne connaît toutefois pas de formule close pour les valeurs propres de $\allnuk$, sauf pour le mélange doublement aléatoire, où elles ont été énoncées par \cite{DS}. C'est d'ailleurs un des principaux problèmes ouverts pour la famille $\allnuk$ (\autoref{pb:formule_close}).

Cependant, on est capable de calculer les valeurs propres de la restriction des opérateurs $\Nu{k}$ à des copies du module $S^{(n-1,1)}$. On appelle la \textit{composante isotypique} $\lambda$ de l'algèbre $\CSn$ la somme directe de toutes les copies de $S^\lambda$ dans la décomposition de $\CSn$. Pour plus de détails sur cette décomposition, voir la \autoref{ssec:defn_Specht}.
\begin{thm}[confirme la conjecture VI.12.2 de \cite{RSW}]\label{thm:vp_S(n-1,1)}
	Les valeurs propres de l'opérateur $\Nu{k}$ restreint à la composante isotypique $(n-1,1)$ sont données par la formule
	\begin{equation}
	\Yvcentermath1
	v_k\left( \Yboxdim{16pt}\tab(1{\ldots}<{\scriptstyle i-1}><{\scriptstyle i+1}>\ldots n,i) \right) = k!\binom{i-2}{k}\binom{2n-i+1}{k}\label{eq:vp_(n-1,1)}
	\end{equation}
	pour $i$ un entier entre $2$ et $n$.
\end{thm}

Pour le démontrer, on procède par induction sur $n$ et on utilise le \autoref{thm:main}. 
On constate aussi que $i-2$ est le type du tableau
\[\Yboxdim{16pt}\tab(1{\ldots}<{\scriptstyle i-1}><{\scriptstyle i+1}>\ldots n,i)\ .\]
Si $i=2$, c'est en particulier un tableau de désarrangement, et la valeur propre $\Nu{k}$ est $0$ pour toute valeur de $k$. Sinon,
\Yvcentermath1
\[\Delta\left(\Yboxdim{16pt}\tab(1{\ldots}<{\scriptstyle 	i-1}><{\scriptstyle i+1}>\ldots n,i) \right) =
\Yboxdim{16pt}\tab(1{\ldots}<{\scriptstyle i-2}><{\scriptstyle i}>\ldots <{\scriptstyle n-1}>,<{\scriptstyle i-1}>). \]

Avec la notation
$\tabw{i,n}  := \Yboxdim{16pt}\tab(1{\ldots}<{\scriptstyle i-1}><{\scriptstyle i+1}>\ldots n,i)$ \index{$\tabw{i,n}$},
on déduit du \autoref{thm:main} ce qui suit.

\begin{lem} \label{lem:rec_eval}
	Les valeurs propres de $\nu_k$ sur la composante isotypique $(n-1,1)$ sont indexées par les tableaux standards de forme $(n-1,1)$ et sont\index{Opérateurs!$v_k$}
	\[v_k(\tabw{i,n}) = \begin{cases}
	0  & \text{ si } i-2 < k,\\
	v_k ( \tabw{i-1, n-1}) + (2n-k-1) \cdot v_{k-1}(\tabw{i-1,n-1}) & 	\text{sinon}.
	\end{cases}\]    
\end{lem}
Ce lemme est la clef pour l'induction dans la preuve du  \autoref{thm:vp_S(n-1,1)}.

\begin{proof}[Preuve du \autoref{thm:vp_S(n-1,1)}]
	On peut vérifier que ça fonctionne pour l'unique tableau standard de forme $(1,1)$, puisque c'est un tableau de désarrangement.
	
	Supposons maintenant que tous les tableaux de taille $n-1$ (et de forme $(n-2,1)$) satisfont l'équation \eqref{eq:vp_(n-1,1)}. On calcule les valeurs propres associées aux tableaux standards de taille $n$ et de type au moins $k$~:
	\allowdisplaybreaks
	\begin{align*}
	v_k(\tabw{i,n}) &=\quad  v_k ( \tabw{i-1, n-1}) + (2n-k-1) \cdot 
	v_{k-1}(\tabw{i-1,n-1})\\
	& \eqexp{\text{Hyp. ind.}}{}k!\binom{i-3}{k}\binom{2n-i}{k} + (2n-k-1) \cdot 
	(k-1)!\binom{i-3}{k-1}\binom{2n-i}{k-1}\\
	&= \quad  k!\binom{i-3}{k}\binom{2n-i}{k} + (2n-k-1) \\ & \qquad \cdot (k-1)!\left(\frac{k}{i-2-k}\binom{i-3}{k}\right) 
	\left(\frac{k}{2n-i-k+1}\binom{2n-i}{k}\right)\\
	& = \quad  k!\binom{i-3}{k}\binom{2n-i}{k} \left( 1 + 
	\frac{k(2n-k-1)}{(2n-i-k+1)(i-2-k)}\right)\\
	& = \quad k!\left(\frac{i-2-k}{i-2}\binom{i-2}{k}\right)
	\left(\frac{2n-i-k+1}{2n-i+1}\binom{2n-i+1}{k}\right) \\
	&\quad \cdot  \left( 1 + \frac{k(2n-k-1)}{(2n-i-k+1)(i-2-k)}\right)\\
	& = \quad k!\binom{i-2}{k} \binom{2n-i+1}{k}
	\frac{(i-2-k)(2n-i-k+1) + k(2n-k-1)}{(2n-i+1)(i-2)}\\
	& = \quad k!\binom{i-2}{k} \binom{2n-i+1}{k},
	\end{align*}
	où la dernière égalité est obtenu en vérifiant que le numérateur de la fraction est égal à $(2n-i+1)(i-2)$.
\end{proof}
La dernière preuve n'est pas des plus élégantes, c'est pourquoi il serait intéressant de trouver une preuve plus combinatoire du \autoref{thm:vp_S(n-1,1)}.

De plus, cette nouvelle expression pour les valeurs propres associées à certains tableaux nous permet de déduire la caractéristique suivante~:
\begin{prop}
	Les valeurs propres associées aux tableaux $t_{i,n}$ sont monotones croissantes en $i$.
\end{prop}
\begin{proof}
	Rappelons d'abord que, lorsque $i < k+2$, la valeur propre de $\Nu{k}$ associée à $t_{i,n}$ est nulle, car le type (ici, $i-2$) est strictement inférieur à $k$.
	Comme les valeurs propres sont toutes positives, $v_k(t_{i,n}) \leq v_k(t_{i+1,n})$ si $i < k+2$.
	
	Sinon, on vérifie que  $v_k(t_{i,n}) < v_k(t_{i+1,n})$. Ceci est facile une fois qu'on a démontré le \autoref{thm:vp_S(n-1,1)}, et ce qu'on souhaite vérifier peut se traduire ainsi~:
    {\allowdisplaybreaks
	\begin{align*}
		k!\binom{i-2}{k}\binom{2n-i+1}{k} &< k!\binom{i-1}{k}\binom{2n-i}{k}\\
		\iff\quad \frac{2n-i+1}{2n-i+1-k} &< \frac{i-1}{i-1-k}\\
		\iff\quad (2n-i+1)(i-1-k) &< (2n-i+1-k)(i-1)\\
		\iff\quad (2n-i+1)(i-1)-k (2n-i+1) &< (2n-i+1)(i-1)-k(i-1)\\
		\iff\quad 2n-i+1 & > i-1\\
		\iff\quad n &> i-1,
	\end{align*}}
	et la dernière inégalité est toujours vraie, confirmant que  $v_k(t_{i,n}) \leq v_k(t_{i+1,n})$.
\end{proof}
\begin{corl}
   La valeur propre $v_k(t_{n,n})$ est maximale parmi celles indexées par les tableaux de forme $(n-1,1)$. 	
\end{corl}
\paragraph{Deuxième valeur propre}
Comme expliqué à la \autopageref{par:temps_de_melange_Sn-1,1}, il est parfois utile de connaître la deuxième plus grande valeur propre d'une chaîne de Markov (la plus grande est toujours $1$, par le théorème de Perron-Frobenius, \autoref{prop:perron-frobenius}). Pour les mélanges d'une collection d'objets tous distincts, la valeur propre associée au tableau $(n-1,1)$ capture le mouvement d'une seule carte \cite{Dia_rep_th}, et permet ainsi de calculer une valeur minimale pour le temps de mélange. Ce n'est pas étranger au fait que, pour plusieurs opérateurs de mélange, la deuxième valeur propre est associée à un tableau de forme $(n-1,1)$. 

Un candidat naturel pour la deuxième plus grande valeur propre de $\allnuk$ semble donc $v_k(t_{n,n})$. Avec SageMath et le \autoref{thm:main}, on a pu calculer toutes les valeurs propres pour tous les opérateurs $\allnuk$ pour des collections de 16 objets et moins. Dans chacun des cas, la deuxième plus grande valeur propre est celle associée au tableau $t_{n,n}$.
\begin{cj}
	La plus grande valeur propre, outre $k!\binom{n}{k}^2$, du mélange $\Nu{k}$ est 
	$v_k(t_{n,n}) = k!\binom{n-2}{k}\binom{n+1}{k}$.
\end{cj}


\subsection{Conséquences du \autoref{thm:main}}\label{ssec:consequences}
Si l'accent a été autant mis sur la recherche des valeurs propres, c'est notamment parce que ce résultat peut être utile et qu'il résout quelques conjectures. La formule décrivant toutes les valeurs propres (\autoref{thm:main}) a notamment les conséquences suivantes~: elle permet de conclure que toutes les valeurs propres sont des entiers positifs, elle permet de les calculer rapidement et elle pourrait nous donner des estimations sur les valeurs propres en général, ce qui pourrait même mener à un calcul du temps de mélange.

\subsubsection{Valeurs propres entières et positives}
On savait déjà que les valeurs propres des opérateurs $\allnuk$ étaient positives (\autoref{ssec:positivite}). De plus, Reiner, Saliola et Welker ont conjecturé que toutes les valeurs propres sont des entiers, ce que permettent de démontrer les résultats de cette section. Cette conjecture était d'ailleurs une des principales motivations pour étudier les mélanges $\allnuk$ (par rapport aux autres opérateurs de mélange s'exprimant comme la version symétrisée d'un opérateur de BHR, voir \autoref{ssec:gammak_bhr}).

Ici, nous donnons une nouvelle preuve de la positivité des valeurs propres, en montrant que le coefficient dans la formule du \autoref{thm:main} est toujours un entier positif.

\begin{prop}
	Pour tout tableau $t$ de taille $n$, de forme $\lambda$ et de type supérieur ou égal à $k$, $(n+2-k+\diag(t/\Delta(t)))$ est un entier positif.
\end{prop}

Pour démontrer cette proposition, on utilise toutefois une autre définition du type d'un tableau, qui est apparue d'abord dans un document non publié de Victor Reiner et Michelle Wachs. Les deux définitions apparaissent dans \cite{RSW} et sont équivalentes.
\begin{defn}[Définition alternative du \textit{type} d'un tableau]
	Soit $t$ un tableau standard de taille $n$. Le \textit{type} de $t$ est l'unique valeur $j \in \{0, 1, \ldots, n-2, n\}$ telle que
	\begin{itemize}
		\item $1,2,\ldots, j-1$ sont des montées de $t$, et
		\item si $t$ a au moins une descente, la première montée parmi $j+1, \ldots, n$ est à une position $j+l$, avec $l$ pair.
	\end{itemize}
\end{defn}

On peut maintenant donner la preuve de la proposition.
\begin{proof}
	On se base sur la relation entre le type d'un tableau et sa forme. À cet effet, notons l'observation suivante~: le type d'un tableau ne peut être plus grand que la longueur de la première ligne.
	Pour s'en convaincre, il suffit de choisir une valeur $i$ entre $1$ et $\lambda_1$ et un tableau standard de forme $\lambda$. Rappelons qu'une descente survient lorsque la case contenant $i+1$ est placée dans la $j$-ième ligne, avec $i>j$. Comme il n'y a que $\lambda_1$ boîtes dans la première ligne, il doit exister une descente parmi $\{1,2, \ldots, \lambda_1\}$. Sinon, le tableau n'a qu'une ligne et son type est alors $\lambda_1 = n$.
	
	On a donc montré que $\lambda_1 \geq \type(t)$ et, par hypothèse, $\type(t) \geq k$. 
	De plus,
	pour n'importe quel nombre $a$ tel que $\lambda_a > 0$, le diagramme de forme $\lambda$ a au moins $\lambda_1 + a - 1$ boîtes; celles-ci sont les $\lambda_1$ boîtes de la première ligne et les $a$ boîtes de la première colonne jusqu'à la ligne $a$. 
	\[\YFrench
	\newcommand{\ver}{\Yfillcolour{green}}
	\newcommand{\ylw}{\Yfillcolour{yellow}}
	\newcommand{\blu}{\Yfillcolour{cyan}}
	\gyoung(!\ver;!\ylw;;;;;;;!,!\blu;,;,;,;,;,!)\]
	Ainsi, si $a$ est la ligne dans laquelle se trouve l'unique boîte de $t/\Delta(t)$,
	\begin{align*}
	n+2-k+\diag(t/\Delta(t)) & = n+2-k + \lambda_a-a\\
	& > n+2-k-a\\
	& \geq n + 2 -\lambda_1 - a\\
	& \geq n+2- (n+1) > 0. \qedhere
	\end{align*}
\end{proof}

\begin{corl}\label{corl:vp_nuk_entiers}
	Les valeurs propres de tous les opérateurs $\allnuk$ sont des entiers positifs.
\end{corl}
\subsubsection{Un algorithme rapide}
Une autre conséquence du \autoref{thm:main} est que l'on est capable de calculer (plutôt rapidement) les valeurs propres des matrices associées aux $\allnuk$. Jusqu'ici, le mieux que nous savions faire est un calcul naïf utilisant les algorithmes habituels pour trouver les valeurs propres d'une matrice. Or, pour une matrice de taille $n!$, comme celles des opérateurs $\allnuk$, il est impensable d'obtenir un résultat dès que $n$ commence à grandir. On peut désormais calculer des valeurs propres rapidement. On peut même énumérer toutes les valeurs propres, jusqu'à un certain point. Après tout, il y a $n!$ valeurs propres, avec multiplicités.

\subsubsection{Bornes sur les valeurs propres}
Bien qu'on ne soit pas capable de donner une formule close pour les valeurs propres, on aimerait déterminer des bornes pour identifier là où elles se trouvent. De telles bornes seraient notamment utiles pour calculer le temps de mélange. Naturellement, on sait que $0 \leq v_k(t) \leq \binom{n}{k}^2k!$, mais nous aimerions avoir quelque chose de plus circonscrit.

La conjecture suivante est vraie pour tous les tableaux de taille jusqu'à 15; elle a été vérifiée à l'aide du \autoref{thm:main} et du logiciel SageMath \cite{sagemath}, et les valeurs pour les tableaux standards de taille 2 à 4 sont au \autoref{tab:cj_bornes_v_k}.
\begin{cj} \label{cj:bornes_v_k}
	La valeur propre de l'opérateur $\nu_k$ associée au tableau $t$ est bornée~:
	\[ v_k(t) \leq \binom{\type(t)}{k} (n+\lambda_1-\type(t))^k. \]
\end{cj}

\begin{table}
	\caption[Comparaison entre les valeurs propres et la borne donnée par la \autoref{cj:bornes_v_k}.]{Comparaison entre les valeurs propres et la borne donnée par la \autoref{cj:bornes_v_k} pour $n=2,3,4$. Les valeurs propres nulles n'apparaissent pas.}\label{tab:cj_bornes_v_k}
	\[\begin{array}{cccc}\setlength\arraycolsep{15pt}
		t & k & \text{Borne conjecturée} & v_k(t)\\
		\hline 
		\tab(12) & 1 & 4 & 4 \\
		 & 2 & 4 & 2\\
		\hdashline
		\tab(123) & 1 & 9 & 9 \\
		 & 2 & 27 & 18 \\
		 & 3 & 27 & 6 \\
		\hdashline
		  
		\tab(12,3) & 1 & 4 & 4 \\
		\hdashline
		\tab(1,2,3) & 1 & 3 & 1 \\
		\hdashline
		\tab(1234) & 1 & 16 & 16 \\
		 & 2 & 96 & 72 \\
		 & 3 & 256 & 96 \\
		 & 4 & 256 & 24\\
		\hdashline
		\tab(123,4) & 1 & 10 & 10 \\
		& 2 & 25 & 20 \\
		\hdashline
		\tab(12,3,4)& 1 & 8 & 6 \\
		& 2 & 16 & 4 \\
		\hdashline
		\tab(124,3) & 1 & 6 & 6 \\
		\hdashline
		\tab(12,34) & 1 & 5 & 4 \\
		\hdashline
		\tab(14,2,3) & 1 & 5 & 2 \\
		
%
%
%
%
%
	\end{array}\]
\end{table}
Sans pouvoir donner une preuve complète, on peut donner une idée de comment on pourrait prouver une telle affirmation. Remarquons d'abord quelques faits~:
\begin{enumerate}
	\item Le \autoref{thm:main} peut se traduire de la façon suivante~:
	la valeur propre $v_k(t)$ peut se réécrire en fonction de certaines valeurs propres de $\Nu{k-1}$. Précisément~:
	\[ v_k(t) = \sum_{i=1}^{\type(t)-k+1} (n+2-i-k+\diag(\Delta^{i-1}t/\Delta^i t)) \Nu{k-1}(\Delta^i t). \]
	\label{item:thm_main_rec}
	\item La \autoref{prop:desa_vs_bandes}, qui affirme que $t/\Delta^{\type(t)}t$ est une bande horizontale, permet de déduire que l'indice diagonal de la bande horizontale $\diag(t /\Delta^{\type(t)}t)$ est au plus l'indice diagonal des $\type(t)$ cellules en bas à droite; ceci est d'ailleurs représenté graphiquement à la \autoref{fig:diagonales_et_diag}.
\begin{figure}
	\centering
	\begin{minipage}{0.3\linewidth}
		\centering
		\Yvcentermath1
		\[\YFrench \Yboxdim{6mm} \yngres(0,3,2,1,1) \]
	\end{minipage}
	\quad
	\begin{minipage}{0.3\linewidth}
		\centering
		\begin{tikzpicture}
		\node[draw,minimum height=5mm, minimum width=5mm] (11) at(0,0) {};
		\node[draw, minimum height=5mm, minimum width=5mm] (12) at(0.5,0) {};
		\node[draw,minimum height=5mm, minimum width=5mm] (13) at(1,0) {};
		\node[draw,minimum height=5mm, minimum width=5mm] (21) at(0,.50) {};
		\node[draw,minimum height=5mm, minimum width=5mm] (22) at(0.5,0.5) {};
		\node[draw,minimum height=5mm, minimum width=5mm] (31) at(0,1) {};
		\node[draw,minimum height=5mm, minimum width=5mm] (41) at(0,1.5) {};
		\draw[->, color=red, dotted, thick] (1.5,0.5) to (.5,-.5) node[below left] {2};
		\draw[->, color=orange, dotted, thick] (1,0.5) to (0,-.5) node[below left] {1};
		\draw[->, color=brown, dotted, thick] (1,1) to (-.5,-.5) node[below left] {0};
		\draw[->, color=green, dotted, thick] (0.5,1) to (-0.5,0) node[below left] {-1};
		\draw[->, color=blue, dotted, thick] (0.5,1.5) to (-0.5,.5) node[below left] {-2};
		\draw[->, color=purple, dotted, thick] (0.5,2) to (-.5,1) node[below left] {-3};
		\end{tikzpicture}
	\end{minipage}
	\caption[Indice diagonal des cellules du diagramme $(3,2,1,1)$, dessiné avec les diagonales.]{Les cellules ayant le plus gros indice diagonal sont placées en bas à droite. Sur l'image à droite, on voit bien les diagonales.}\label{fig:diagonales_et_diag}
\end{figure}
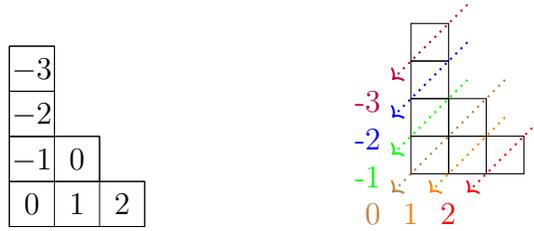

	Ainsi, pour tout $j$ tel que $0< j \leq \type(t)$,
	\[ \diag(t/\Delta^{j}t)\leq \sum_{i=1}^{j} \lambda_1-i = \lambda_1 j - \frac{j(j+1)}{2}. \]\label{item:sum_diag_bande_h}
\end{enumerate}

La première des deux affirmations suggère qu'on puisse procéder par induction sur $k$. Le cas de base, $k=1$, est alors le mélange doublement aléatoire.
Pour ce mélange, la valeur propre \mbox{$v_1(t)\leq \type(t)\left(\lambda_1 - \frac{\type(t)+1}{2}\right)$}, comme le démontre le calcul suivant~:
\begin{align*}
v_1(t) &\eqexp{fait \ref{item:thm_main_rec}}{} \sum_{i=1}^{\type(t)} \left( n+1-i+ \diag(\Delta^{i-1} t/\Delta^{i} t)  \right)\\
&=\qquad (n+1)\type(t) - \frac{\type(t) (\type(t)+1)}{2} + \sum_{i=1}^{\type(t)} \diag(\Delta^{i-1} t/\Delta^{i} t) \\
&\overset{\mathclap{\normalfont\mbox{\tiny \text{fait \ref{item:sum_diag_bande_h}}}}}{\leq\qquad}(n+1) \type(t) - \frac{\type(t) (\type(t)+1)}{2} + \type(t)\left(\lambda_1 - \frac{\type(t)+1}{2}\right)\\
&=\qquad \type(t) \left(n - \type(t) + \lambda_1 \right).
\end{align*}

Si nous savions démontrer la \autoref{cj:bornes_v_k}, nous aurions de meilleures chances de trouver le temps de mélange des opérateurs $\allnuk$, dont nous discutons tout de suite.
	

\subsubsection{Temps de mélange}
Le \textit{temps de mélange} d'une chaîne de Markov peut être compris comme le temps nécessaire pour qu'une distribution de probabilités sur laquelle on exécute une chaîne de Markov soit  assez proche de la distribution stationnaire.
Par exemple, le temps de mélange d'un paquet de cartes à jouer est le temps qu'il faut pour que toutes les permutations aient à peu près la même probabilité de représenter l'ordre des cartes dans le paquet après ce temps.

Formellement, le temps de mélange se définit en fonction de la \textit{distance de variation totale}, une mesure entre deux distributions de probabilités qui permet de dire à quel point elles sont éloignées l'une de l'autre. Cette distance prend des valeurs entre $0$ et $1$~: deux distributions sont très éloignées si leur distance de variation totale est près de $1$, et elles sont près l'une de l'autre lorsque la distance se rapproche de $0$.

\begin{defn}
	Pour deux distributions $p$ et $q$ sur un ensemble fini d'états $X$, la distance de variation totale entre $p$ et $q$ est 
	\[ || p - q ||_{VT} = \frac{1}{2}\sum_{x \in X} |p(x) - q(x)|. \]
\end{defn}

Pour les chaînes de Markov, on souhaite souvent comparer la distribution stationnaire (vers laquelle la chaîne converge) et la distribution après un certain temps, puisque ça nous indique le rythme auquel on s'en rapproche. Le temps de mélange est alors le nombre d'itérations de la chaîne de Markov qui sont nécessaires pour que la distance de variation totale soit inférieure à un seuil fixé (souvent entre $\frac{1}{4}$ et $\frac{1}{2}$).
Lorsque le nombre d'états est très grand, calculer la distance de variation totale avec précision peut être laborieux. Il existe différentes techniques pour passer outre cette difficulté. Une qui est particulièrement intéressante pour nous est celle-ci, d'abord donnée spécifiquement pour un autre mélange dans \cite{DiSh}, puis généralisée. La version ci-dessous s'apparente à celle qu'on retrouve dans \cite{LPW}.

\begin{lem}[Lemme 14 de \cite{DiSh}, lemme 12.18 de \cite{LPW}]
	Pour les marches aléatoires sur un groupe $G$, la distance de variation totale après $s$ itérations du mélange est bornée supérieurement par
	\begin{equation}
		\sqrt{\frac{1}{4}\sum_{j=2}^{|G|} \lambda_j^{2s}},\nonumber
	\end{equation}
	où $1 > \lambda_2 \geq \ldots \geq \lambda_{|G|} 
	\geq -1$ sont les valeurs propres de l'opérateur de mélange (vu comme une chaîne de Markov).
	
	Dans le cas du groupe symétrique, on peut réécrire cette équation avec $\lambda_t$ qui désigne la valeur propre associée au tableau $t$~:
	\begin{equation}
	\sqrt{\frac{1}{4}\sum_{\text{tableaux standards $t$}} f^{\mathrm{forme}(t)} \lambda_t^{2s}},\label{eq:DV_Sn}
	\end{equation}
	où $f^{\mathrm{forme}(t)}$ est le nombre de tableaux standards de la même forme que $t$.
\end{lem}

C'est notamment pour calculer cette somme que nous aimerions avoir une formule close ou des bornes supérieures suffisamment parlantes pour exprimer les valeurs propres; on serait alors peut-être en mesure de donner une borne supérieure sur le temps de mélange.

On peut toutefois avoir une idée de ce à quoi pourrait ressembler cette somme; c'est le contenu du prochain paragraphe.
\paragraph{Distance de variation et valeurs propres des tableaux de forme $(n-1,1)$}\label{par:temps_de_melange_Sn-1,1}
Grâce au \autoref{thm:vp_S(n-1,1)}, on connaît une formule close pour certaines valeurs propres; on a même conjecturé que la deuxième valeur propre pour chacun des opérateurs est parmi celles pour lesquelles on connaît cette expression. On pourrait vouloir calculer la somme des puissances de ces valeurs propres, comme dans l'équation \eqref{eq:DV_Sn}~: cette équation est
\[\sqrt{\frac{1}{4}\sum_{\text{$t$ de forme $(n-1,1)$}} (n-1) \lambda_t^{2s}} = \sqrt{\frac{1}{4}\sum_{i=2}^{n} (n-1) \left( \frac{\binom{i-2}{k}\binom{2n-i+1}{k}}{\binom{n}{k}^2}\right)^{2s}}.\]
Ceci ne nous donne pas tout à fait une borne supérieure satisfaisante, puisque la grande majorité des valeurs propres sont omises. Cependant, la contribution de ces valeurs propres est la plus grande parmi toutes les valeurs propres et le temps de mélange est généralement du même ordre de grandeur que celui donné par les plus grandes valeurs propres (pour plus de détails, voir par exemple la notion de temps de relaxation dans \cite[Section 12.2]{LPW}).

En faisant la somme uniquement sur les valeurs propres associées aux tableaux de forme $(n-1,1)$, on peut arriver à la conjecture suivante~:
\begin{cj}[Megan Bernstein, communication personnelle]
	Lorsque $k$ est petit par rapport à $n$, le temps de mélange des opérateurs $\nu_k$ est, à constante près, au plus $\frac{3}{4k}n \log(n)$; formellement, on dit que le temps de mélange est $\mathcal{O}(\frac{3}{4k}n \log(n))$.
\end{cj}
Cette dernière conjecture est démontrée lorsque $k=1$; c'est le temps de mélange du mélange doublement aléatoire, établi par Megan Bernstein et Evita Nestoridi \cite{BN}.

\section{Vecteurs propres}
Si l'étude des valeurs propres est en soi intéressante, connaître les vecteurs propres permettrait de calculer rapidement les puissances des opérateurs $\allnuk$. Anton Dieker et Franco Saliola ont décrit une procédure pour obtenir une base de tous les espaces propres de l'opérateur de mélange doublement aléatoire, $\Nu{1}$, à partir des vecteurs propres dans le noyau (qui sont inconnus). Il serait intéressant de voir si cette description peut s'étendre aux autres opérateurs de la famille $\allnuk$.

Sans savoir explicitement comment construire les vecteurs propres, on peut démontrer que les opérateurs $\allnuk$ partagent les mêmes vecteurs propres. Plus explicitement~:

\begin{thm}\label{thm:vec_prop_communs_allnuk}
	Il existe une base de $\CSn$ dont les vecteurs sont des vecteurs propres pour chacun des opérateurs $\allnuk$.
\end{thm}

\begin{proof}
	On montrera en fait que n'importe quelle famille de matrices symétriques qui commutent peut être associée à une base de vecteurs propres pour tous ses éléments. Rappelons que les opérateurs $\allnuk$ commutent (\autoref{thm:commutativite_nu}) et sont symétriques, tel qu'expliqué à la \autoref{ssec:def_matrices_nuk}.
	
	Le théorème spectral affirme que toute matrice réelle et symétrique est diagonalisable (voir, par exemple, \cite{HJ}, théorème 2.5.6). C'est donc le cas des opérateurs $\nu_k$.
	
	De plus, des opérateurs diagonalisables qui commutent sont simultanément diagonalisables (voir par exemple le théorème 1.3.19 de \cite{HJ}).
	
	Enfin deux matrices simultanément diagonalisables partagent une base de vecteurs propres~: soit $P$ une matrice telle que $P^{-1}MP$ soit diagonale pour toute matrice $M$ de la famille. Les vecteurs colonnes de $P$ sont en effet les vecteurs propres de chacune des matrices de la famille.
\end{proof}

En plus de présenter un intérêt intrinsèque et de permettre de calculer les puissances des matrices, les vecteurs propres sont utiles notamment pour connaître le temps de mélange. Pour les opérateurs de mélange en général, on peut aussi utiliser les vecteurs propres, par exemple, pour estimer la distribution stationnaire, lorsqu'il ne s'agit pas forcément de la distribution uniforme. Pour connaître davantage d'usages des vecteurs propres, voir \cite[section 2]{DPR}.

\chapter{Preuves des théorèmes principaux}\label{chap:preuves}
\begin{figure}[H]
	\begin{center}
		\begin{tikzcd}[ampersand replacement=\&, column sep=12em, row 
	sep=5.5em]
			M^{\lambda+\Box} \arrow[r, "\textrm{\large Projection}", two 
			heads] \& S^{\lambda+\Box}\\
			M^\lambda \arrow[u, "\textsf{\large(A)}"] \& S^\lambda \arrow 
			[l, "\textrm{\large Inclusion}", hook'] \arrow[ul, 
			"\textsf{\large(B)}"] \arrow[u, "\textsf{\large(C)}"]
		\end{tikzcd}
		\caption[Schéma de la preuve du \autoref{thm:main}.]{Schéma de la preuve du \autoref{thm:main}. Rappelons que nous cherchons à comprendre comment les valeurs propres de $S^\lambda$ et de $S^{\lambda + \Box}$ sont liées.}\label{fig:schema_pr_val_pr}
	\end{center}
\end{figure}
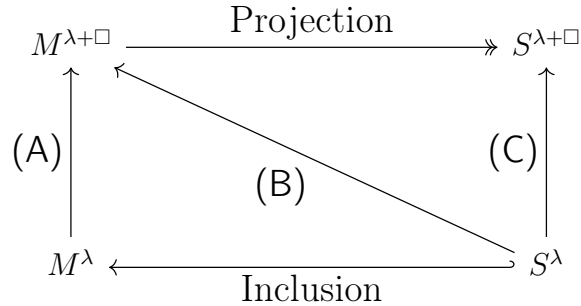

Ce chapitre est dédié à prouver le \autoref{thm:main}, donnant toutes les valeurs propres de $\Nu{k}$ pour une collection de $n$ éléments, et le \autoref{thm:commutativite_nu}, qui explique que les opérateurs $\allnuk$ commutent. Pour parvenir à démontrer le premier, on utilise la stratégie suivante, synthétisée à la \autoref{fig:schema_pr_val_pr}~:
\begin{enumerate}[label=(\textsf{\Alph*})]
	\item On s'attarde à ce qui se passe sur les mots. On se place dans le contexte de l'algèbre des mots, et plus précisément des \textit{modules de permutation} (\autoref{ssec:modules_perm}). Ceux-ci ne sont pas des modules simples, et ils contiennent donc plusieurs copies de modules de Specht, mais ils ont l'avantage d'avoir les mots comme base; comprendre l'action d'un opérateur y est donc beaucoup plus facile.
	
	On veut savoir ce qu'il se passe entre $\CSn$ et $\C S_{n+1}$; en plus de mélanger, on doit donc insérer une lettre dans notre mot. On compare donc les actions d'insérer puis mélanger et celles de mélanger d'abord et d'insérer ensuite.
	
	C'est le travail que l'on retrouve à la \autoref{sec:alg_mots}
	\item Comme les valeurs propres de l'algèbre $\C S_{n+1}$ sont obtenues de celles de $\CSn$ de la même façon que celles des modules simples $S^{\lambda+\Box}$ proviennent de celles de $S^\lambda$, on restreint les équations aux modules de Specht $S^\lambda$. Malheureusement, le résultat de plusieurs équations se situe dans le module de permutations plutôt que dans le module de Specht. On en parle à la \autoref{sec:restr_Specht}.
	\item On résout ce problème en projetant le résultat sur les modules de Specht; on procède en projetant sur une composante isotypique $S^{\lambda}$, à la \autoref{ssec:proj_iso}.
	
	On obtient ainsi toutes les valeurs propres, dont on parle à la \autoref{sec:vp_preuve}.
\end{enumerate}

L'autre théorème important qui est démontré est que les opérateurs $\allnuk$ commutent entre eux. Pour démontrer celui-ci, on réutilise le \autoref{thm:thm1}, qui explique les relations entre l'insertion d'une lettre et les opérateurs $\allnuk$ sur les modules de permutation. Puis, on décrit une relation similaire avec la suppression d'une lettre. Conjointement avec une définition alternative de $\Nu{k}$ (\autoref{thm:lem8}), cela nous permettra de démontrer par induction la commutativité.

Afin de faciliter la lecture de ce chapitre, la liste de tous les théorèmes (incluant les propositions, lemmes et corollaire) qui y sont présentés se trouve à l'\autoref{chap:annexe_thm}.

\section{L'algèbre des mots}\label{sec:alg_mots}
Pour ce chapitre, nous travaillons sur l'algèbre des mots (définie à la \autopageref{par:algebre_mots}). Rappelons qu'un \textit{mot} est une suite de lettres présentées dans un certain ordre. Dans le contexte de cette thèse, tous les mots sont finis et on peut donc définir leur longueur~: il s'agit du nombre de \textit{lettres} (incluant les répétitions) que contient le mot. Les lettres sont tirées d'un ensemble fini totalement ordonné, appelé \textit{alphabet} et noté $A$. Pour un mot de longueur $n$, on utilise l'alphabet des lettres $1$ à $n$ (sans forcément utiliser toutes les lettres). Notons que tous les mots de longueur $n$ peuvent se rapporter à cet alphabet; l'injection explicite sera détaillée sous peu.

Certains mots ont la particularité que toutes leurs lettres sont distinctes. Ces mots sont peu étudiés en combinatoire des mots parce que la recherche de motifs y est difficile, mais ils ont ceci de très intéressants qu'il s'agit des permutations. Appliquer un mélange sur une telle suite est intéressant, puisqu'il s'agit de mélanger des éléments tous distincts. Un paquet de cartes dont on aurait retiré les jokers est un exemple d'un mot dont chaque lettre a une unique occurrence.

\paragraph{Passage à l'alphabet de $1$ à $n$}\label{par:alph_contenu} Pour passer d'un mot quelconque de longueur $n$ à un mot sur l'alphabet des nombres de $1$ à $n$, on restreint l'alphabet aux seules lettres dans le mot. Ensuite, on réordonne les lettres de la plus fréquente à celle qui l'est le moins; en cas d'égalité, on conserve l'ordre alphabétique. Les deux exemples qui suivent détaillent ce processus.

\begin{ex}
	Le nom québécois \textit{efface} est composé de quatre des six premières lettres de l'alphabet~: $a$, $c$, $e$ et $f$, les deux dernières ayant deux occurrences dans le mot. En réordonnant l'alphabet, on trouve l'injection suivante~:
	\begin{center}
		{\setlength\tabcolsep{15pt}
		\begin{tabular}{ccccccc}
			a & b & c & d & e & f & \ldots \\
			$\updownarrow$ & - & $\updownarrow$ & - & $\updownarrow$ & $\updownarrow$ & \ldots \\
			3 & - & 4 & - & 1 & 2 & \ldots
		\end{tabular}}
	\end{center}
	Notons que $e$ et $f$ précèdent $a$ et $c$ étant donné leur nombre d'occurrences supérieur.
	
	\begin{center}
		\begin{tabular}{|l|c|c|c|c|c|c|}
			\hline
			& e&f&f&a &c &e\\
			\hline
			Valeur de la lettre après l'injection & 1 & 2 & 2 & 3 & 4 & 1 \\
			\hline
		\end{tabular}
	\end{center}
	Ainsi, l'alphabet du mot étant $\{a,c,e,f\}$, le mot \textit{efface} est réécrit $122341$.
\end{ex}

\begin{ex}
	Le mot \textit{banane} est composé des première, deuxième, cinquième et quatorzième lettres de l'alphabet, et $a$ et $n$ sont plus fréquentes.
	L'injection est 
	\begin{center}
		{\setlength\tabcolsep{15pt}
		\begin{tabular}{cccc}
			a & b & e & n\\
			$\downarrow$ & $\downarrow$ & $\downarrow$ & $\downarrow$ \\
			1 & 3 & 4 & 2
		\end{tabular}}
	\end{center}
	
	Le mot \textit{banane} devient $312124$ sur l'alphabet composé des nombres de $1$ à $6$.
\end{ex}

\paragraph{Algèbre des mots}\label{par:algebre_mots}
À partir d'un alphabet $A$, les mots de $A^*$ sont formés par concaténation, c'est-à-dire que, pour une certaine suite de lettres $(a_1, a_2,\ldots, a_n)$, on obtient un mot en collant les lettres de la suite tout en préservant leur ordre~: $a_1a_2\ldots a_n$. La concaténation est une opération associative, mais elle n'est pas commutative. Si l'on considère les lettres comme des variables non-commutatives, on peut définir l'\textit{algèbre des mots}, notée $\A$\index{Modules!$\A$}, formée des polynômes dont les variables sont les lettres de $A$~: les monômes sont alors des mots (avec leur coefficient). Le produit non-commutatif de deux monômes est la concaténation et la somme est définie de la façon habituelle sur les espaces vectoriels.

On peut décomposer l'algèbre $\A$ de plusieurs façons. Une d'elle est de décomposer $\A$ en sous-espaces vectoriels engendrés par les mots de même longueur~: on note alors $\C A^n$\index{Modules!$\C A^n$} le sous-espace vectoriel dont les éléments sont les combinaisons linéaires de mots de longueur $n$. Alors, $\A = \bigoplus_{n\in\N}\ \C A^n$.

\subsection{Actions du groupe symétrique sur $\C A^n$}\label{ssec_contenant_par:melanges_alg_du_groupe}
Le sous-espace vectoriel $\C A^n$ est un module pour l'action du groupe symétrique. Il existe en fait deux actions du groupe symétrique sur $\C A^n$~: une action (à gauche) qui permute les lettres de l'alphabet, puis une autre action (à droite) qui permute les positions des lettres. Les mélanges sont un exemple de cette action à droite.

\paragraph{Action à gauche~: les \textit{coccyx} des \textit{nonnes}}
L'action à gauche du groupe symétrique sur les mots permute les lettres de l'alphabet. Par exemple, si on travaille sur l'alphabet latin habituel, la permutation (notée en cycles) \mbox{$\sigma = \{(cn),(ey),(sx)\}$} échange, dans le mot \textit{coccyx}, les $c$ et les $n$, les $e$ et les $y$, et enfin les $s$ et les $x$, pour chaque occurrence de ces lettres. Elle fixe toutes les autres lettres de l'alphabet (le $o$, par exemple). Ainsi, $\sigma \cdot\textit{coccyx} = \textit{nonnes}$.

Plus formellement, l'action sur les lettres de l'alphabet est définie comme ceci, pour un mot $w = a_1\ldots a_n$ et une permutation $\sigma \in S_n$~:
\[\sigma \cdot w  = \sigma(a_1)\ldots \sigma(a_n).  \]
Évidemment, cette action ne change pas la longueur des mots et préserve donc les modules $\{\C A^n\}_{n \in \N}$.

\paragraph{Action à droite~: \textit{imaginer} une \textit{migraine}}
Une autre action par permutation peut être définie sur $\C A^n$. Chaque mot garde le même ensemble de lettres qui le constituent, mais l'ordre de ces lettres est changé. Cette action à droite est définie comme suit pour un mot $w = a_1\ldots a_n$ et une permutation $\sigma$~:
\[ w \cdot \sigma = a_{\sigma(1)} \ldots a_{\sigma(n)}.\]
Comme cette action ne change pas la longueur des mots, elle préserve les modules $\{\C A^n\}_{n\in \N}$.

\begin{ex}
	Le mot \textit{imaginer} est une anagramme du mot \textit{migraine}, c'est-à-dire qu'ils ont le même nombre d'occurrences de chaque lettre.
	On peut obtenir \textit{imaginer} à partir de \textit{migraine} par la permutation $21536784$~:
	\[\textit{migraine}\cdot 21536784 = \textit{imaginer}. \]
	
\end{ex}

\begin{ex}
	Soit $w = \underline{13154}$, un mot sur l'alphabet $\{1,2,3,4,5\}$, et soit la permutation $\sigma = 43251$. Notons que, pour cet exemple, les mots ont été soulignés pour les distinguer de la permutation qui agit sur le mot.\\
	L'action à gauche de $\sigma$ sur (l'alphabet de) $w$ donne
	\[ 43251\cdot \underline{13154} = \underline{42415} .\]
	En revanche, l'action à droite de $\sigma$ sur (les positions de) $w$ donne
	\[ \underline{13154} \cdot 43251  = \underline{51341} .\]
\end{ex}

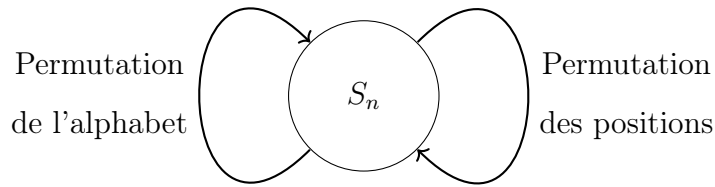
\begin{figure}[H]
	\vspace{-1cm}
	\centering
	\begin{tikzpicture}
		\node[draw, circle, minimum size=2cm] (A) at (0,5) {$S_n$};
		\path[->, thick] (A) edge [loop left, looseness=5, in=135, out=225] node {$\stackunder[5mm]{\text{{\normalsize Permutation}}}{\text{{\normalsize de l'alphabet}}}$}  (A);
		\path[->,thick] (A) edge [loop right, looseness=5, in=-45, out=45] node {$\stackunder[5mm]{\text{{\normalsize Permutation}}}{\text{{\normalsize des positions}}}$}  (A);
	\end{tikzpicture}
	\vspace{-1cm}
	\caption{Action par permutation, à gauche et à droite.}
\end{figure}

\paragraph{Les mélanges sont des éléments de l'algèbre du groupe symétrique qui agissent sur $\C A^n$. } On peut interpréter les opérateurs de mélange comme des éléments de $\C S_n$ qui agissent sur $\C A^n$~: en effet, chaque mot de $A^n$ est envoyé, par un opérateur de mélange, sur une combinaison linéaire d'anagrammes de ce mot. 

La notation en termes de sommes de permutations des opérateurs $\allnuk$ n'est pas particulièrement élégante. Toutefois, on peut exprimer les éléments de l'algèbre du groupe $\CSn$ correspondant aux opérateurs $\{\pi_k\}_{k\in\N}$, dont les $\allnuk$ sont des multiples des symétrisés~: 
\[\pi_k = \sum_{1 \leq i_1 \leq \ldots \leq i_k \leq n} \sum_{\sigma \in S_k} 1 \cdot \ldots \cdot \widehat{i_1} \cdot \ldots \cdot \widehat{i_k} \cdot \ldots \cdot n \cdot i_{\sigma(1)} \cdot \ldots \cdot i_{\sigma(k)}.\]

En revanche, $\pi_k^\top$ requiert la notion de battage, présentée au \autoref{chap:chap0}. Ainsi,
\[\pi_k^\top = \sum_{1 \leq i_1 \leq \ldots \leq i_k \leq n} \sum_{\sigma \in S_k} \left( 1 \cdot \ldots \cdot \widehat{i_1} \cdot \ldots \cdot \widehat{i_k} \cdot \ldots \cdot n\right) \shuffle \left(i_{\sigma(1)} \cdot \ldots \cdot i_{\sigma(k)}\right), \]
et $\pi_k^\top$ est aussi une somme de permutations.

Comme le produit d'éléments de l'algèbre du groupe est aussi dans l'algèbre, les opérateurs de mélange $\allnuk$ sont bien dans l'algèbre.

\subsection{Modules de permutations}\label{ssec:modules_perm}

Le \textit{contenu} d'un mot (parfois aussi appelé son \textit{évaluation}) est une liste \mbox{$\lambda = (\lambda_1, \ldots, \lambda_n)$}, où $\lambda_i$ est le nombre d'occurrences de la lettre $i$ dans ce mot. 
\begin{ex}
	Le contenu du mot $123321121$ est $(4,3,2)$, parce qu'il contient quatre occurrences de la lettre $1$, trois occurrences de la lettre $2$ et que $2$ apparaît deux fois. En suivant la définition donnée plus haut, on trouve $(4,3,2,0,0,0,0,0,0)$. Toutefois, on omet généralement les $0$ à la fin du tuple.
\end{ex}

Il est assez pratique de considérer $\lambda$ comme un partage (comme on le faisait à la \autoref{ssec:vp_mots}). Pour que ça soit possible, il faut réordonner les valeurs de $\lambda$. C'est possible de le faire en réordonnant l'alphabet en suivant l'injection présentée à la \autopageref{par:alph_contenu}.
Enfin, pour considérer $\lambda$ comme un partage, on retire les $0$ qui apparaissent à la fin; c'est ce qui est fait dans l'exemple précédent.

On considère le sous-espace vectoriel $M^\lambda$\index{Modules!$M^\lambda$} engendré par tous les mots de contenu $\lambda$. Pour l'action à droite du groupe symétrique, $M^\lambda$ est un module; bien sûr, changer les positions des lettres du mot n'altère pas le contenu de ce mot. On appelle $M^\lambda$ le \textit{module de permutations} associé à $\lambda$ et une décomposition de l'algèbre des mots de longueur $n$ est la suivante~:
\[ \C A^n \cong \bigoplus_{c\ \vDash n} M^{\text{partage}(c)}, \]
où $\text{partage}(c)$ est le partage obtenu de la composition $c$ en ordonnant ses parts en ordre décroissant.

\subsubsection{Modules de permutation et modules de Specht}
Contrairement aux modules de Specht (présentés à la \autoref{ssec:defn_Specht}), les modules de permutations ne sont pas simples. Ils contiennent différentes copies des modules de Specht, et nous présentons leur décomposition ici. Cela est particulièrement utile pour comprendre les théorèmes de la \autoref{ssec:vp_mots}, mais aussi pour réaliser les preuves dans le présent chapitre.

\begin{thm}[Règle de Young pour les modules de permutation]\label{thm:regle_de_young_mod_perm}
	Les modules de permutations se décomposent en modules simples de la façon suivante~:
	\[M^\mu \cong \bigoplus_{\lambda \unrhd \mu} m_{\lambda\mu} S^\lambda,\]
	où $m_{\lambda\mu}$ est le nombre de tableaux semi-standards de forme $\lambda$ et de contenu $\mu$. Pour la définition de tableaux semi-standards, voir à la \autopageref{par:SSYT}.
\end{thm}

\begin{ex}
	On décompose le module de permutations $M^{\smalldia(2,2,1)}$ de la façon suivante.
	\begin{itemize}
		\item On énumère les partages qui dominent $(2,2,1)$~: il s'agit de $(5)$, $(4,1)$, $(3,2)$, $(3,1,1)$ et $(2,2,1)$.
		\item Pour chacun de ces partages $\lambda$, on énumère les tableaux semi-standards de forme $\lambda$ qui ont pour contenu $(2,2,1)$. Ceux-ci sont illustrés au \autoref{tab:SSYT221}.
		\item On trouve, à l'aide du \autoref{tab:SSYT221}, que 
		\[ M^{\smalldia(2,2,1)} \cong S^{\smalldia(5)} \oplus2\  S^{\smalldia(4,1)} \oplus2\ S^{\smalldia(3,2)} \oplus S^{\smalldia(3,1,1)} \oplus S^{\smalldia(2,2,1)}.\]
	\end{itemize}
\begin{table}
	\caption{Tableaux semi-standards de contenu $(2,2,1)$.}\label{tab:SSYT221}
	\[
	\renewcommand{\arraystretch}{2.5}
	\Yvcentermath1
	\begin{array}{|c|c|c|c|c|}
	\hline
	\dia(5) & \dia(4,1) & \dia(3,2) & \dia(3,1,1) & \dia(2,2,1)\\
	\hline
	\tab(11223) & \tab(1123,2) \quad \tab(1122,3) & \tab(112,23) \quad \tab(113,22)  & \tab(112,2,3) & \tab(11,22,3)\\
	\hline
	\end{array}\]
\end{table}
\end{ex}
	 
La règle de Young pour les modules de permutations est en fait l'ingrédient-clé de la preuve du \autoref{thm:val_pr_mots}, qui liste toutes les valeurs propres d'un opérateur de mélange $\Nu{k}$ lorsque la collection d'objets peut admettre des répétitions. La règle de Young pour les modules de permutations est l'élément qui dit quelles valeurs propres de $\Nu{k}$ sont aussi des valeurs lorsque la collection a des objets répétés, en plus de permettre de calculer la multiplicité de chacune des valeurs propres.

\subsection{Opérateurs sur l'algèbre des mots}

On aborde maintenant les collections triées d'objets comme des mots. Pour ce faire, on réutilise le vocabulaire de la combinatoire des mots, présenté au \autoref{chap:chap0}. Les objets deviennent alors des lettres, et la collection triée est vue comme un mot. 
En général, on utilise l'alphabet latin $\{a,b,\ldots\}$ pour représenter les lettres. Par exemple, $a$ pourrait aussi représenter $6$ ou $\alpha$ si nous modifions l'alphabet.

\subsubsection{Insertion et suppression de lettres}
On définit, sur l'algèbre des mots deux opérateurs linéaires, correspondant respectivement à l'insertion d'une lettre dans les mots et au retrait d'une lettre.

\paragraph{Insertion d'une lettre~: $\sh_a$}\index{Opérateurs!$\sh_a$}
Nous avons défini au \autoref{chap:chap0} le battage de deux mots $u$ et $v$ en disant que c'était la combinaison linéaire formelle de tous les mots qu'on peut obtenir en intercalant des lettres de $u$ et de $v$, tout en s'assurant que $u$ et $v$ soient des sous-mots disjoints dans le résultat. Lorsqu'un des mots est composé d'une seule lettre. Ainsi, on définit la fonction $\sh_a(w)$ comme $a\shuffle w$.

Soit $\vec{e_a}$\index{$\vec{e_a}$} le vecteur qui est de la même longueur que la taille de l'alphabet et qui ne contient que des $0$, à l'exception d'un $1$ à la position de la lettre $a$. Explicitement, $\sh_a$ est une application $M^\lambda \to M^{\lambda+\vec{e_a}}$ qui est définie par
\[ \sh_a (w) = \sum_{i=1}^{n+1} w_1 \cdot \ldots \cdot w_{i-1} \cdot a \cdot w_i  \cdot \ldots \cdot w_n. \]
\begin{ex}
	Lorsque $w = abaa$,
	\begin{align*}
	\sh_a(abaa) &= 2 \cdot aabaa + 3 \cdot abaaa,\\
	\sh_b(abaa) &= babaa + 2 \cdot abbaa + ababa + abaab.
	\end{align*}
\end{ex}
\paragraph{Retrait d'une lettre~: $\del_a$}\index{Opérateurs!$\del_a$} De la même façon qu'on peut exprimer l'insertion d'une occurrence de la lettre, on peut exprimer l'opération contraire, soit le retrait d'une occurrence de la lettre $a$. On note ainsi l'opérateur de retrait \mbox{$\del_a : M^{\lambda} \to M^{\lambda-\vec{e_a}}$}, défini explicitement comme
\[ \del_a(w)  = \sum_{1 \leq i \leq n} \delta_{a, w_i} w_1 \ldots w_{i-1}w_{i+1} \ldots w_n. \]

\begin{ex}
	Lorsque $w = abaa$,
	\begin{align*}
	\del_a(abaa) &= baa + 2\cdot aba\\
	\del_b(abaa) &= aaa
	\end{align*}
\end{ex}

Comme expliqué dans \cite{DS}, les opérateurs $\del_a$ et $\sh_a$ sont mutuellement adjoints. En effet, considérons le produit scalaire $\langle\cdot ,\cdot\rangle$\index{$\langle\cdot ,\cdot\rangle$ (mots)} associé à la base des mots. Ainsi, $\langle u, w \rangle$ vaut $1$ lorsque $u$ et $w$ sont égaux, et $0$ sinon. Alors,
\[ \langle \sh_a(u), v\rangle  = \langle u, \del_{a}(v) \rangle. \]

\begin{prop}[Équations (11) et (12) de \cite{DS}]\label{prop:commutShDel}
	Les opérateurs d'insertion commutent entre eux, et il en est de même pour les opérateurs de suppression~:
	\begin{align*}
	\sh_a\circ \sh_b &= \sh_b \circ \sh_a,\\
	\del_a \circ \del_b &= \del_b \circ \del_a.
	\end{align*}
\end{prop}
Cependant, les opérateurs d'insertion  et de suppression ne commutent pas entre eux, tel qu'expliqué au \autoref{lem:lem35}.

\subsubsection{L'opérateur de remplacement $\theta_{a,b}$}\index{Opérateurs!$\theta_{a,b}$}

Enfin, on définit un opérateur pour remplacer une lettre dans un mot par une autre. Ainsi, on note $\theta_{a,b}$ l'opérateur visant à remplacer, de toutes les façons possibles, une occurrence de $a$ par une occurrence de $b$.
Formellement,
\[ \theta_{a,b} (w) = \sum_{i=1}^n \delta_{w_i,a}\ w_1\cdot \ldots\cdot w_{i-1}\cdot b \cdot w_{i+1} \cdot \ldots \cdot w_n.  \]

\begin{lem}\label{lem:lem11_thet}
	L'opérateur de remplacement, $\theta_{a,b}$, est un morphisme de $S_n$-modules à droite. \\
	En particulier, $\theta$ commute avec les actions à droite de $S_n$.
\end{lem}
\begin{proof}
	Pour prouver le \autoref{lem:lem11_thet}, on montre d'abord le deuxième énoncé, c'est-à-dire que $\theta$ commute avec les actions de $S_n$. Comme l'opérateur $\theta$ est une application linéaire, il est automatiquement un morphisme de $S_n$-modules s'il commute avec les actions de $S_n$.
	
	Il faut donc montrer que pour tout mot $w \in \langle A\rangle$ et pour toute permutation $\sigma \in S_n$, 
	\[(\theta_{a,b} (w))\cdot \sigma = \theta_{a,b} (w \cdot \sigma).\]
	Comme les mots forment une base de l'algèbre des mots $\A$, montrer cet énoncé démontre également la commutativité sur toute l'algèbre.
	Or, 
	\begin{align*}
	(\theta_{a,b} (w)) \cdot \sigma = &\ \left(\sum_{\substack{1\leq i \leq n\\ w_i = a}} w_1 \ldots w_{i-1} \cdot b\cdot  w_{i+1} \ldots w_n \right) \cdot \sigma \\
	= &\ \sum_{\substack{1 \leq i \leq n\\ w_{\sigma(i)}= a}} w_{\sigma(1)} \ldots w_{\sigma(i-1)} \cdot b\cdot w_{\sigma(i+1)} \ldots w_{\sigma(n)}\\
	= &\ \theta_{a,b}(w \cdot \sigma),
	\end{align*}
	achevant la démonstration du \autoref{lem:lem11_thet}.
\end{proof}

\subsection{Opérateurs de mélange $\allnuk$ sur l'algèbre des mots}
\begin{prop}\label{prop:defn_op}\index{Opérateurs!$\nu_k$}
	À l'aide des opérateurs $\sh$ et $\del$, on peut donner cette nouvelle définition des opérateurs $\allnuk$ dans l'algèbre des mots~:
	\[ \Nu{k} = \sum_{1\leq a_1 \leq \ldots \leq a_k\leq n} \frac{1}{\prodnuk} \shdelUnak.  \]
\end{prop}

\begin{rmq}
	La proposition plus haut apparaissait déjà (sans preuve) dans l'article de Dieker et Saliola \cite[remarque 37]{DS}, mais une erreur dans la formule faisait que l'égalité ne tenait que pour le cas du mélange doublement aléatoire (lorsque $k=1$.)
\end{rmq}

\begin{proof}
	En utilisant la \autoref{defn:nu_k_cartes} de $\Nu{k}$, on sait que $\Nu{k}(w)$ retire $k$ lettres (pas nécessairement distinctes) à $w$, puis réinsère les mêmes lettres dans n'importe quel ordre.

	Supposons que les lettres retirées soient $a_1, \ldots, a_k$. L'ordre dans lequel ces lettres sont retirées n'importe pas, puisque les opérateurs $\del_{a_i}$ et $\del_{a_j}$ commutent (\autoref{prop:commutShDel}).
	Sans perte de généralité, on peut supposer que $a_1 \leq \ldots \leq a_k$.
	Le retrait puis la réinsertion des lettres correspond à $\shdelUnak$.
	
	Pour décrire $\Nu{k}$, on doit considérer la somme sur tous les ensembles de lettres $\{a_1, \ldots, a_k\}$ possibles. Toutefois, dans l'expression 
	\[\smashoperator{\sum_{1\leq a_1 \leq \ldots \leq a_k\leq n}}\  \shdelUnak, \]
	le terme $\shdelUnak$ peut apparaître plus d'une fois~: lorsque l'ensemble $\{a_1, \ldots, a_k\}$ compte $m$ occurrences d'une certaine lettre \mbox{$a_{i} = \ldots = a_{i+m}$}, il y a $m!$ façons d'ordonner les lettres dans la suite d'inégalités. En faisant le produit pour chacune des lettres de l'alphabet, on trouve qu'il y a \mbox{$\prodnuk$} façons possibles d'ordonner $\{a_1, \ldots, a_k\}$.
	
	En éliminant les doublons, on trouve que l'opérateur $\Nu{k}$ s'écrit aussi \[ \sum_{1\leq a_1 \leq \ldots \leq a_k\leq n} \frac{1}{\prodnuk} \shdelUnak.  \qedhere \]
\end{proof}

\subsection{Identités sur l'algèbre des mots}
La prochaine sous-section est dédiée à prouver le \autoref{thm:thm1}, énoncé plus loin. On y présente différentes propriétés des opérateurs $\Nu{k}$, ainsi que leurs relations avec d'autres opérateurs sur l'algèbre des mots. Notamment, on étudie ce qui se produit lorsqu'on insère ou retire une lettre du mot.

Les identités qui suivent permettent en quelque sorte de quantifier \textit{à quel point} les différents opérateurs sur les mots ne commutent pas. S'ils commutaient, alors le côté droit de ces équations serait nul.

\begin{lem}\label{lem:lem35}
	Sur l'algèbre des mots de longueur $n$, les identités suivantes sont vérifiées~:
	\begin{enumerate}[label=\roman*.]
		\item $\del_b \circ \sh_a - \sh_a\circ\del_b = \theta_{b,a} + (n+1)\delta_{a,b}\ \Id $ \label{lem35i},
		\item $\theta_{b,c} \circ \sh_a - \sh_a\circ\theta_{b,c} = \delta_{a,b}\sh_c $\label{regleDeNadia},
	\end{enumerate}
	où $\delta_{a,b}$ vaut $1$ si $a=b$ et $0$, sinon.
\end{lem}

\begin{proof} La première de ces deux équations est déjà présente dans le lemme 36 de \cite{DS}.
	La seconde partie se démontre en distinguant les différents cas. Pour nous aider, on peut rappeler une autre partie du lemme 36 de \cite{DS}~:
	\[ \theta_{a,c} \circ \sh_a - \sh_a \circ \theta_{a,c} = \sh_c. \]
	Dans l'identité qu'on souhaite démontrer, cela correspond au cas où $a = b$.
	
	Si $a\neq b$, alors on doit montrer que $\sh_a$ et $\theta_{b,c}$ commutent. Or, si $a$ et $b$ sont distinctes, l'insertion de la lettre $a$ ne change rien aux façons de modifier un $b$ en $c$, et réciproquement.
\end{proof}

Ce dernier lemme est grandement utilisé dans les preuves tout au long de ce chapitre. Il est notamment utilisé dans la preuve du prochain lemme.

Celui-ci est plutôt technique et donne une façon plus compliquée de décrire l'action de $\sh_a \circ \Nu{k-1}$. Cependant, il semble que ce détour soit nécessaire pour la preuve du \autoref{thm:thm1}.
\begin{lem}\label{lem:lem2}
	Sur l'algèbre des mots de longueur $n$,
	\begin{align*}
		\sh_a \circ \Nu{k-1} = \sum_{i=1}^k \sum_{1 \leq a_1 \leq \ldots \leq a_k\leq n}& \delta_{a, a_i} \frac{1}{\prodnuk} \sh_a \circ \sh_{a_1} \circ \ldots \circ \nonumber\\ &\widehat{\sh_{a_{i}}}\circ \ldots\circ\sh_{a_k}\circ \del_{a_1} \circ\ldots \widehat{\del_{a_{i}}}\circ \ldots \circ \del_{a_k}.
	\end{align*}
\end{lem}

Ici, notons que $\sh_{a_i}$ et $\del_{a_{i}}$ sont retirées du côté droit de l'équation afin d'obtenir $\nu_{k-1}$.
\begin{proof}
On sait déjà, par la définition de $\Nu{k-1}$ (\autoref{prop:defn_op}) que 
\begin{align}
\sh_a \circ \Nu{k-1} = \sh_a \sum_{1 \leq b_1 \leq \ldots \leq b_{k-1}\leq n}& \frac{1}{\prod_{j\in [n]}\#\{l \in [k-1] \mid b_l = j \}!}\nonumber\\ & \sh_{b_1} \circ \ldots \circ\sh_{b_{k-1}}\circ \del_{b_1} \circ\ldots \circ \del_{b_{k-1}}.\label{eq:pr_lem2}
\end{align}

Pour chacun des termes de la somme à la droite de l'équation, on fixe le multi-ensemble (ou ensemble admettant les répétitions) $\{b_1, \ldots, b_{k-1}\}$. On note alors 
\[ \{a, b_1, \ldots, b_{k-1}\} = \{a_1, \ldots, a_k\}. \]
Supposons qu'il y ait $q$ occurrences de $a$ dans $\{b_1, \ldots, b_{k-1}\}$. Alors,
\[\sum_{i=1}^k \delta_{a, a_i} = q+1 \textrm{ et } \prodnuk = (q+1)\prod_{j\in [n]}\#\{l \in [k-1] \mid b_l = j\}!. \]

Comme la lettre $a$ est fixée, l'équation \eqref{eq:pr_lem2} vaut aussi
\begin{align*}
\sh_a \sum_{1 \leq a_1 \leq \ldots \leq a_k\leq n}& \frac{\sum_{i=1}^k\delta_{a, a_i}}{\prodnuk} \sh_{a_1} \circ \ldots \circ \\ &\widehat{\sh_{a_i}}\circ \ldots\circ\sh_{a_k}\circ \del_{a_1} \circ\ldots \widehat{\del_{a_{i}}}\circ \ldots \circ \del_{a_k}.
\qedhere\end{align*}
\end{proof}

L'équation qui suit nous donne une nouvelle façon d'écrire $\Nu{k}$ en fonction de $\Nu{k-1}$, ce qui nous permettra de réaliser des preuves par induction faisant intervenir $\Nu{k}$.
\begin{thm}\label{thm:lem8}
	Sur l'algèbre des mots,
	\begin{align*}
		\sum_{a = 1}^{n}\sh_a \circ\ \Nu{k-1}\circ \del_a = k\ \Nu{k}.
	\end{align*}
\end{thm}

\begin{proof}
	La preuve se fait essentiellement à l'aide du \autoref{lem:lem2} et de quelques égalités~:
	{\allowdisplaybreaks
	\begin{align*}
	\sum_{a = 1}^n \sh_a \circ& \Nu{k-1} \circ \del_a
	\\
	\eqexp{\autoref{lem:lem2}}{& \sum_{a = 1}^n \sum_{i=1}^k \sum_{1 \leq a_1 \leq \ldots \leq a_k \leq n} \delta_{a, a_i} \frac{1}{\prodnuk} \sh_a \circ \sh_{a_1} \circ \ldots \circ\\
		& \widehat{\sh_{a_{i}}}\circ\ldots\circ\sh_{a_k}\circ \del_{a_1} \circ\ldots \widehat{\del_{a_{i}}}\circ \ldots \circ \del_{a_k}\circ \del_a }\\
	\eqexp{\autoref{prop:commutShDel}}{&\quad  \sum_{a = 1}^n \sum_{i=1}^k \sum_{1 \leq a_1 \leq \ldots \leq a_k \leq n} \delta_{a, a_i} \frac{1}{\prodnuk} \sh_{a_1} \circ \ldots \circ\\
		& \sh_{a_{i-1}} \circ \sh_a \circ\sh_{a_{i+1}}\circ \ldots \circ\sh_{a_k}\circ \del_{a_1} \circ\ldots \del_{a_{i-1}}\circ \del_{a} \circ \del_{a_{i+1}}\circ \ldots \circ \del_{a_k}}\\
	\eqexp{}{& \sum_{a = 1}^n \sum_{1 \leq a_1 \leq \ldots\leq a_k \leq n} \frac{\#\{l \in [k] \mid a_l = a \}}{\prodnuk} \sh_{a_1} \circ \ldots \circ
		\sh_{a_k}\circ \del_{a_1} \circ\ldots \circ \del_{a_k}}\\
	\eqexp{}{& \sum_{1 \leq a_1 \leq \ldots \leq a_k \leq n} \frac{\sum_{a = 1}^n  \#\{l \in [k] \mid a_l = a \}}{\prodnuk} \sh_{a_1} \circ \ldots \circ
		\sh_{a_k}\circ \del_{a_1} \circ\ldots \circ \del_{a_k}}\\
	\eqexp{}{& \sum_{1 \leq a_1 \leq \ldots \leq a_k \leq n} \frac{k}{\prodnuk} \sh_{a_1} \circ \ldots \circ
		\sh_{a_k}\circ \del_{a_1} \circ\ldots \circ \del_{a_k}}\\
	\eqexp{\autoref{prop:defn_op}}{&\quad k\ \Nu{k}},
	\end{align*}}
 	où l'avant-dernière égalité vient du fait que $\sum_{a = 1}^n  \#\{l \in [k] \mid a_l = a\}$ compte le nombre total d'éléments dans $\{a_1, \ldots, a_k\}$.
\end{proof}

On est maintenant en mesure d'énoncer le théorème principal de cette section. C'est une règle essentielle pour procéder à l'induction de $\Nu{k-1}$ à $\Nu{k}$.
\begin{thm}\label{thm:thm1}
	Sur les mots de longueur $n$,
	\[ \Nu{k}\circ \sh_a - \sh_a \circ\Nu{k} = (n+2-k) \sh_a \circ \Nu{k-1} + \sum_{1\leq b \leq n} \sh_b \circ \theta_{b,a}\circ\Nu{k-1}. \]
\end{thm}

De façon grossière, cette équation décrit la différence entre insérer d'abord une lettre dans notre mot de longueur $n$, puis mélanger avec $\Nu{k}$, ou faire ces deux opérations dans l'ordre inverse. On pourra ainsi comparer l'action de $\Nu{k}$ sur un mot de taille $n$ et sur un mot de taille $n+1$.

Dans la suite d'égalités, des parties de certaines équations ont été encadrées; ceci indique des parties qui changeront à la ligne suivante.
\begin{proof}
	La preuve se fait par induction sur $k$. Le cas de base est celui du mélange doublement aléatoire, $\nu_1$, et correspond au théorème 38 de \cite{DS}~:
	\[ \Nu{1} \circ \sh_a - \sh_a\circ\Nu{1} = (n+1)\sh_a + \sum_{1\leq b \leq n} \sh_b \circ \theta_{b,a}. \]
	Supposons maintenant que, pour une certaine valeur $k \geq 2$, et sur des mots de longueur $n-1$,
	\[ \Nu{k-1}\circ \sh_a - \sh_a \circ\Nu{k-1} = (n+2-k) \sh_a \circ \Nu{k-2} + \sum_{1\leq b \leq n-1} \sh_b \circ \theta_{b,a}\circ\Nu{k-2}. \]
	C'est notre hypothèse d'induction.
	
	Alors, sur des mots de longueur $n$,
	\allowdisplaybreaks
	\begin{align*}
		\Nu{k}\circ \sh_a - &\sh_a \circ \Nu{k}\\
		\eqexp{\autoref{thm:lem8}}{&  \Nu{k}\circ \sh_a - \frac{1}{k} \left( \bbox{\sh_a \sum_{1\leq b \leq n} \sh_b} \circ \Nu{k-1}\circ \del_b \right)}\\
		\eqexp{\autoref{prop:commutShDel}}{&\Nu{k}\circ \sh_a - \frac{1}{k} \left( \sum_{1\leq b \leq n} \sh_b  \circ \bbox{\sh_a \circ \Nu{k-1}}\circ \del_b \right)} \\
		\eqexp{\text{Hyp. induction}}{ &\Nu{k}\circ \sh_a - \frac{1}{k} \bigg( \sum_{1\leq b \leq n} \sh_b  \big(\Nu{k-1} \circ \sh_a -(n+2-k) \sh_a\circ\Nu{k-2}}\\
		& - \sum_{1\leq c \leq n-1} \sh_c \circ \theta_{c,a} \circ\Nu{k-2} \big)\circ \del_b \bigg)\\
		\eqexp{\autoref{prop:commutShDel}}{}&\Nu{k}\circ \sh_a - \frac{1}{k} \bigg( (\sum_{1\leq b \leq n} \sh_b  \circ \Nu{k-1} \circ \bbox{\sh_a \circ \del_b}) - (n+2-k)\\& \cdot  \sh_a (\bbox{\sum_{1\leq b \leq n} \sh_b\circ \Nu{k-2}\circ\del_b})
		 - \sum_{1\leq b,c \leq n} \sh_b\circ \sh_c \circ \theta_{c,a} \circ \Nu{k-2} \circ \del_b \bigg)\\
		\eqexp{$\stackunder{\text{\autoref{lem:lem35}}}{\text{\autoref{thm:lem8}}}$}{ } & \Nu{k}\circ \sh_a - \frac{1}{k} \bigg( (\sum_{1\leq b \leq n} \bbox{\sh_b  \circ \Nu{k-1} \circ (\del_b}\circ \sh_a - \theta_{b,a}))-(n+1) \sh_a \circ \Nu{k-1}\\
		& - (n+2-k)(k-1) \sh_a \circ \Nu{k-1}- \sum_{1\leq b,c \leq n} \sh_b\circ \sh_c \circ \theta_{c,a} \circ \Nu{k-2} \circ \del_b \bigg) \\
		\eqexp{\autoref{thm:lem8}}{ & \Nu{k}\circ \sh_a - \frac{1}{k} \bigg( k\ \Nu{k} \circ \sh_a-\sum_{1\leq b \leq n} \sh_b  \circ \Nu{k-1} \circ \theta_{b,a} -(n+1) \sh_a \circ \Nu{k-1}}\\& - (n+2-k)(k-1) \sh_a \circ \Nu{k-1}- \sum_{1\leq b,c \leq n} \bbox{\sh_b\circ \sh_c} \circ \theta_{c,a} \circ \Nu{k-2} \circ \del_b \bigg) \\
		 \eqexp{\autoref{prop:commutShDel}}{ &\ \frac{1}{k} \bigg( (\sum_{1\leq b \leq n} \sh_b  \circ \bbox{\Nu{k-1} \circ \theta_{b,a}}) +\left((n+1) + (n+2-k)(k-1)\right) \sh_a \circ \Nu{k-1}}\\&  +  \sum_{1\leq b,c \leq n} \sh_c\circ \sh_b \circ \theta_{c,a} \circ \Nu{k-2} \circ \del_b \bigg) \\
		 \eqexp{\autoref{lem:lem11_thet}}{ & \left(n+3-k -\frac{1}{k}\right) \sh_a \circ \Nu{k-1} + \frac{1}{k} \bigg( (\sum_{1\leq b \leq n} \sh_b \circ \theta_{b,a} \circ \Nu{k-1} )}\\&  +  \sum_{1\leq c \leq n} \sh_c \sum_{1\leq b \leq n} \bbox{\sh_b \circ \theta_{c,a}} \circ \Nu{k-2} \circ \del_b \bigg) \\
		 \eqexp{\autoref{lem:lem35}}{ & \left(n+3-k -\frac{1}{k}\right) \sh_a \circ \Nu{k-1} + \frac{1}{k} \bigg( (\sum_{1\leq b \leq n} \sh_b \circ \theta_{b,a} \circ \Nu{k-1} )}\\&  +  \sum_{1\leq c \leq n} \bbox{\sh_c} \sum_{1\leq b \leq n}(\theta_{c,a} \circ \sh_b \circ \Nu{k-2} \circ \del_b - \delta_{b,c} \bbox{\sh_a} \circ \Nu{k-2} \circ \del_b)\bigg) \\
		 \eqexp{\autoref{prop:commutShDel}}{&\quad  \left(n+3-k -\frac{1}{k}\right) \sh_a \circ \Nu{k-1} + \frac{1}{k} \bigg( (\sum_{1\leq b \leq n} \sh_b \circ \theta_{b,a} \circ \Nu{k-1} )}\\&  +  \sum_{1\leq c \leq n} \sh_c \circ \theta_{c,a}\bbox{\sum_{1\leq b \leq n} \sh_b \circ \Nu{k-2} \circ \del_b} - \bbox{\sum_{1\leq b \leq n} \sh_a \sh_b \circ \Nu{k-2} \circ \del_b}\bigg) \\
		 \eqexp{\autoref{thm:lem8}}{ &\quad \left(n+3-k -\frac{1}{k}\right) \sh_a \circ \Nu{k-1} + \frac{1}{k} \bigg( (\sum_{1\leq b \leq n} \sh_b \circ \theta_{b,a} \circ \Nu{k-1} )}\\&  +  (k-1)\sum_{1\leq c \leq n} \sh_c \circ \theta_{c,a}\circ \Nu{k-1} - (k-1) \sh_a \circ \Nu{k-1}\bigg) \\
		  = \qquad & \left(n+2-k\right) \sh_a \circ \Nu{k-1} +  \sum_{1\leq b \leq n} \sh_b \circ \theta_{b,a} \circ \Nu{k-1}. \qedhere
	\end{align*}
\end{proof}

\section{De l'algèbre des mots aux modules de Specht}\label{sec:restr_Specht}
\begin{figure}
	\begin{center}
		\begin{tikzcd}[ampersand replacement=\&, column sep=12em, row 
			sep=5.5em]
			M^{\lambda+\Box} \arrow[r, "\textrm{Projection}", two 
			heads] \& S^{\lambda+\Box}\\
			M^\lambda \arrow[u, "\text{Théorème \ref{thm:thm1}}"] \& S^\lambda \arrow [l, "\textrm{Inclusion}", hook'] \arrow[ul, "\text{Théorème \ref{thm:thm2}}"] \arrow[u, swap, "\text{Théorème \ref{thm:thm3}}"]
		\end{tikzcd}
		\caption[Schéma détaillé de la preuve du \autoref{thm:main}.]{Schéma de la preuve du \autoref{thm:main}. Rappelons que nous cherchons à comprendre comment les valeurs propres de $S^\lambda$ et de $S^{\lambda + \Box}$ sont liées. Ici, les numéros des théorèmes concernés ont été ajoutés. C'est le \autoref{thm:thm3} qui est utilisé pour la preuve du \autoref{thm:main}.}
	\end{center}
\end{figure}
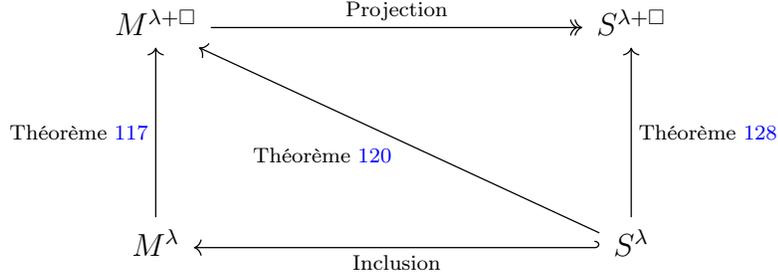

\paragraph{Modules de Specht}
Les modules de Specht sont les modules simples du groupe symétrique et sont définis à la \autoref{ssec:defn_Specht}. Comme c'est décrit à la \autoref{fig:recap_vp_tableaux}, ils sont particulièrement intéressants parce qu'ils forment le lien entre les valeurs propres et les tableaux standards.

Les deux prochains lemmes décrivent le comportement des opérateurs $\allnuk$ et $\theta_{b,a}$ lorsqu'ils sont restreints aux modules de Specht.
\begin{lem}\label{lem:lem10_nu_k_Slambda->Slambda}
	Pour toute valeur de $k$ et pour tout partage $\lambda$, l'image de $\nu_k|_{S^\lambda} $ est comprise dans $S^\lambda$.
\end{lem}

\begin{proof}
	Les modules de Specht $S^\lambda$ sont des modules, et, pour toute action $\rho$ du groupe symétrique, l'image de $\rho|_{S^\lambda}$ est comprise dans $S^{\lambda}$. Pour prouver le théorème, on se rappelle que les opérateurs de mélange sont des actions du groupe symétrique (\autoref{ssec_contenant_par:melanges_alg_du_groupe}).\qedhere
\end{proof}

\begin{lem}[Lemme 43 de \cite{DS}]\label{lem:lem43}
	Soit $\lambda$ un partage.  Si $a<b$, alors $\theta_{b,a}|_{S^\lambda} = 0$.
\end{lem}

\begin{proof}[Idée de preuve]
	L'argument-clé pour la preuve du lemme ci-dessous est la proposition 2.4.5 de \cite{sagan}. Celle-ci dit que si un homomorphisme de $S^\lambda$ vers $M^{\lambda'}$ est non-nul, alors $\lambda \unrhd \lambda'$. Or, le partage $\lambda$ est dominé par $\lambda'= \lambda+ \vec{e_a} - \vec{e_b}$ lorsque $a< b$, et $\theta_{b,a}|_{S^\lambda}$ doit forcément être nul.
\end{proof}

On a maintenant tous les outils en main pour prouver le théorème suivant~:

\begin{thm}\label{thm:thm2}
	Soit $\lambda \vdash n$. Alors, pour toute lettre $a$,
	\[ \left(\nu_k \circ \sh_a - \sh_a \circ \nu_k\right)|_{S^\lambda} = (n+2-k)\left(\sh_a \circ \Nu{k-1}\right)|_{S^\lambda} + \sum_{1\leq b \leq a} \left(\sh_b \circ \theta_{b,a}\circ \Nu{k-1}\right)|_{S^\lambda}.  \]
\end{thm}

\begin{proof}
	On procède par suite d'égalités~:
	\begin{align*}
		\big(\Nu{k}\circ \sh_a - \sh_a \circ& \Nu{k}\big)|_{S^\lambda}\\ \eqexp{\autoref{thm:thm1}}{& (n+2-k)\left(\sh_a \circ \Nu{k-1}\right)|_{S^\lambda} + \left(\sum_{1\leq b \leq n} \sh_b \circ \theta_{b,a}\right)\circ \Nu{k-1}|_{S^\lambda}}\\
		\eqexp{$\stackunder{\text{\autoref{lem:lem10_nu_k_Slambda->Slambda}}}{\text{\autoref{lem:lem43}}}$}{& (n+2-k)\left(\sh_a \circ \Nu{k-1}\right)|_{S^\lambda} + \left(\sum_{1\leq b \leq a} \sh_b \circ \theta_{b,a}\right)\circ \Nu{k-1}|_{S^\lambda}}.
	\end{align*}
	Dans la dernière égalité, le \autoref{lem:lem10_nu_k_Slambda->Slambda} est nécessaire pour utiliser le \autoref{lem:lem43} et, ainsi, affirmer que la somme est nulle lorsque $a < b$.
\end{proof}

Enfin, un fait intéressant concernant l'insertion d'une lettre dans un élément d'un module de Specht est rapporté ici et sera utilisé plus tard.

\begin{prop}[Lemme 46 de \cite{DS}]\label{prop:lem45DS}
	L'image de $(\sh_a)|_{S^\lambda}$ est contenue dans un sous-module de $M^{\lambda+\vec{e_a}}$ isomorphe à $\bigoplus_{\mu = \lambda+\vec{e_r},\  r\leq a} S^\mu$.
\end{prop}

\subsection{Projecteurs isotypiques}\label{ssec:proj_iso}

Considérons une algèbre de groupe $\C G$ pouvant se décomposer en modules simples, notés $\{S^{i}\}_{i}$, et notons $S^{[i]} = \bigoplus m_i S^i$,\index{Modules!$S^{[i]}$} où $m_i$ est le nombre de copies de $S^i$ dans $\C G$. On suppose que $S^i$ et $S^j$ ne sont pas isomorphes si $i\neq j$. Ainsi, \[\C G \cong c_1S^{[1]} \oplus c_2S^{[2]} \oplus \ldots \oplus c_lS^{[l]}.\]

Un module est dit \textit{isotypique} s'il est une somme directe de sous-modules simples isomorphes, et $S^{[i]}$, tel que défini plus haut, est la $i$-ième composante isotypique.

Il est connu que tout module admet une unique décomposition maximale en sous-modules isotypiques. Dans le cas de l'algèbre du groupe symétrique, cette décomposition se fait sur les partages~:
\[ \CSn \cong \bigoplus_{\lambda\vdash n} (\CSn)^{[\lambda]} \cong \bigoplus_{\lambda\vdash n} (S^{\lambda})^{\oplus f^\lambda}, \]
où $f^\lambda$ désigne le nombre de tableaux standards de forme $\lambda$ et correspond à la multiplicité de $S^\lambda$ dans $\CSn$ (d'après la règle de Young, \autoref{thm:regle_de_Young}).

\begin{defn}
	Soit $G$ un groupe fini et soit $S$ un $G$-module simple. Soit $\chi^S$ le caractère de $S$. Alors, on note $p_S$ l'élément de $\C G$
	\[ p_S = \frac{\dim(S)}{|G|} \sum_{g\in G} \overline{\chi^S(g)} g, \]
	où $\overline{\chi}$ désigne le conjugué complexe de $\chi$.\index{$p_S$}
	La multiplication à droite par cet élément constitue la fonction sur l'algèbre du groupe $G$
	\begin{align*}
		\isoproj_S : \C G &\to \C G^{[S]}\\
					 w & \mapsto w\cdot p_S.\index{Opérateurs!$\isoproj$}
	\end{align*}
	
	Dans le cas du groupe symétrique, il sera démontré prochainement que $\isoproj_\lambda$ est un projecteur isotypique sur le module simple $S^\lambda$. C'est également le cas pour les groupes en général et la preuve se fait de la même façon.
	\end{defn}

\begin{prop}[Propriétés du projecteur isotypique $\isoproj$]\label{prop:proprietes_isoproj}
	La fonction $\isoproj_\mu$ satisfait les énoncés suivants~:
	\begin{enumerate}
		\item elle commute avec les actions du groupe symétrique.
		\item c'est un morphisme de $S_n$-modules.
		\item c'est un projecteur isotypique de $\CSn$ vers la composante isotypique associée au module simple $S^\mu$.
	\end{enumerate}
\end{prop}

Avant de pouvoir donner une démonstration de la proposition, il y a deux notions qu'il est essentiel de définir. Il s'agit du produit de deux caractères et des relations d'orthogonalité. Ces deux concepts bien connus sont l'objet des prochains paragraphes.
\paragraph{Produit de caractères}
\begin{defn}\label{defn:prod_caracteres}
	Soit $G$ un groupe et soit $\varphi$ et $\psi$ deux fonctions $G\to\C$ (des caractères, par exemple). Leur produit interne est donné par l'équation
	\[ \langle\varphi, \psi \rangle = \frac{1}{|G|}\sum_{g \in G} \overline{\varphi(g)}\ \psi(g).\]\index{$\langle\cdot, \cdot\rangle$ (caractères)}
\end{defn}

Ce produit est utile pour comprendre les caractères irréductibles, c'est-à-dire ceux qui forment une base de tous les caractères. Les caractères irréductibles sont associés aux modules simples, ce qui nous donne la proposition suivante~:
\begin{prop}[Première relation d'orthogonalité]\label{prop:1e_relation_orthogonalite}
	Si $U$ et $V$ sont des modules simples, les caractères associés à $U$ et à $V$ sont orthonormaux, c'est-à-dire que
	\[ \langle \chi^U, \chi^V \rangle = \delta_{U,V} = \left\{ \begin{array}{ll}
	1 & \text{si $U = V$}\\
	0 & \text{si $U \neq V$}
	\end{array}\right.. \]
\end{prop}

\paragraph{De retour aux projecteurs isoptypiques}
On a maintenant les outils nécessaires pour démontrer la \autoref{prop:proprietes_isoproj}.
\begin{proof}[Démonstration de la \autoref{prop:proprietes_isoproj}]
	\begin{enumerate}
		\item \textit{Elle commute avec les actions du groupe symétrique.}
		On montre d'abord que $\isoproj_\mu$ est dans le centre de $\CSn$. Pour ce faire, choisissons $\tau\in S_n$. Étant donné que le caractère est une fonction centrale, $\chi^\mu(\sigma) =  \chi^\mu(\tau\sigma\tau^{-1})$, et 
		\begin{align}
			p_\mu \cdot \tau &= \left(\frac{f^\mu}{n!} \sum_{\sigma \in S_n} \overline{\chi^\mu(\sigma)} \sigma \right)\tau \nonumber\\
			&= \frac{f^\mu}{n!} \sum_{\tau \sigma \tau^{-1}\in S_n} \overline{\chi^\mu(\sigma)} \tau \sigma\nonumber\\
			&=  \tau \left(\frac{f^\mu}{n!} \sum_{\sigma\in S_n} \overline{\chi^\mu(\sigma)} \sigma\right)\nonumber\\
		    &= \tau\cdot p_\mu.\label{arg:isoprojZSn}
		\end{align}
		Ainsi, l'action à droite de toute permutation $\sigma$ commute avec $\isoproj_\mu$~:
		\[ \isoproj_\mu(v \cdot \sigma ) =(v \cdot \sigma)  \cdot p_\mu  \overset{\text{\tiny \normalfont \'Eq.\eqref{arg:isoprojZSn}}}{=\ }(v \cdot p_\mu)\cdot \sigma  = \isoproj_\mu(v) \cdot \sigma. \]
		\item \textit{C'est un morphisme de $S_n$-modules.}
		On a vu que $\isoproj_\mu$ commute avec les actions du groupe symétrique. De plus, 
		\[\isoproj_\mu(u + v) = (u+v) \cdot p_\mu  = u \cdot p_\mu + v \cdot p_\mu  = \isoproj_\mu(u) + \isoproj_\mu(v) \]
		et $\isoproj_\mu$ est un morphisme de $S_n$-modules.
		\item \textit{C'est un projecteur isotypique de $\CSn$ vers la composante isotypique associée au module simple $S^\mu$.}
		On fait cette preuve en deux étapes~: d'abord, on montre que, pour tout module de Specht $S^\lambda$, $\isoproj_\mu|_{S^\lambda}$ est nul si $\lambda \neq \mu$ et l'identité sinon. Ensuite, on montre comment on étend cette propriété à $(\isoproj_\mu)|_{\CSn}$.
		\begin{enumerate}
			\item On a montré plus tôt que $\isoproj_{\mu}$ est un morphisme de $S_n$-modules. Comme $S^\mu$ et $S^\lambda$ sont des modules simples, on peut utiliser le lemme de Schur pour montrer que $\isoproj_{\mu}|_{S^\lambda}$ est un multiple de l'identité~: $\isoproj_{\mu}|_{S^\lambda} = c\cdot \Id$.
			
			Pour connaître la valeur de $c$, on calcule la trace de $\isoproj_{\mu}|_{S^\lambda}$. Notons qu'une notion essentielle pour ce calcul est le produit des caractères, présenté à la \autoref{defn:prod_caracteres}. Ainsi,
			\allowdisplaybreaks
			\begin{align*}
			\tr(\isoproj_{\mu}|_{S^\lambda}) =& \tr\left( \underbrace{\frac{f^\mu}{n!} \sum_{\sigma \in S_n} \overline{\chi^\mu(\sigma)} \sigma}_{p_\mu \cdot \Id|_{S^\lambda}}|_{S^\lambda}\right)\\
			= & \frac{f^\mu}{n!} \sum_{\sigma \in S_n} \overline{\chi^\mu(\sigma)} \tr(\sigma|_{S^\lambda})\\
			= & \frac{f^\mu}{n!} \sum_{\sigma \in S_n} \overline{\chi^\mu(\sigma)} \chi^\lambda(\sigma)\\
			= & f^\mu \langle\chi^\mu, \chi^\lambda\rangle\\
			= & \left\{\begin{array}{ll}
				f^\mu  & \text{si $S^\lambda \cong S^\mu$}\\
				0 & \text{sinon}
				\end{array} \right.,
			\end{align*}
			par la première relation d'orthogonalité.
			
			Ainsi, puisque la trace de $\isoproj_{\mu}|_{S^\lambda}$ est égale à la dimension de $S^\lambda$ ou nulle, on conclut que, $\isoproj_{\mu}|_{S^\lambda}$ est l'identité lorsque $\lambda = \mu$ et nul sinon.
			\item Considérons la décomposition de $\CSn$ en modules simples~:
			\[ \CSn \cong \bigoplus_{\lambda \vdash n} S^{[\lambda]} \cong \bigoplus_{\lambda \vdash n}  (S^\lambda)^{\oplus f^\lambda}. \]
			Par ce qui précède, $\isoproj_{\mu}(S^\lambda)$ est non-nul seulement si $\mu = \lambda$, auquel cas il agit comme l'identité.
			Ainsi,
			\[\isoproj_{\mu}(\CSn) \cong (S^\mu)^{\oplus f^\mu} \cong S^{[\mu]}, \]
			ce qui permet de démontrer que $\isoproj_{\mu}$ agit véritablement comme le projecteur isotypique.\qedhere
		\end{enumerate}
	\end{enumerate}
\end{proof}

Les propriétés des projecteurs isotypiques sont notamment utiles pour évaluer avec quels opérateurs ils commutent (comme au \autoref{lem:lem11_commut}) et avec lesquels ils ne commutent pas (\autoref{lem:lem13}).
\begin{lem}\label{lem:lem11_commut}
	Le projecteur isotypique $\isoproj_\mu$ commute avec l'opérateur de remplacement $\theta_{a,b}$~:
	\[ \theta_{a,b}(\isoproj_\mu (v)) = \isoproj_\mu(\theta_{a,b}(v)). \]
\end{lem}

\begin{proof}
	Comme $\theta_{a,b}$ est un morphisme de $S_n$-modules (\autoref{lem:lem11_thet}), il commute avec les éléments de l'algèbre du groupe. Or, $p_\mu$ est un élément de l'algèbre du groupe.
	\[ \theta_{a,b}(\isoproj_\mu (v)) = \theta_{a,b}(v \cdot p_\mu) =  \theta_{a,b}(v)  \cdot p_\mu  = \isoproj_\mu(\theta_{a,b}(v)).\qedhere \]
\end{proof}

Nous avons introduit la notion de projecteur isotypique pour avoir une façon de revenir vers les modules simples même en ajoutant une lettre. Ainsi, pour simplifier la notation, nous reprenons la notation de \cite{DS} et introduisons l'opérateur $\projlift$, qui est la composition de l'insertion puis de la projection isotypique~:
\begin{align*}
	\projlift_a^{\mu} : M^{\lambda} & \to S^{\mu}\\
	v & \mapsto \isoproj_{\mu}(\sh_a(v))\index{Opérateurs!$\projlift$}
\end{align*}
lorsque $\mu = \lambda+\vec{e_a}$ et $v \in M^\lambda$.

\begin{lem}\label{lem:lem13}
	Soit $\lambda = (\lambda_1, \ldots, \lambda_l) \vdash n$, $a \in [l+1]$ et $\mu = \lambda + \vec{e_r}$, avec $r \in [a]$.  
	Alors,
	 \begin{align*}
		 \big(\nu_k \circ &\projlift_a^\mu - \projlift_a^\mu \circ \nu_k \big)|_{S^\lambda} =\\
		 &(n+3-k+\lambda_a - a) \left(\projlift_a^\mu \circ \Nu{k-1}\right)|_{S^\lambda} + \left(\sum_{r \leq b < a} \theta_{b,a} \circ \projlift_b^\mu \circ \Nu{k-1}\right)|_{S^\lambda}.
	 \end{align*}
\end{lem}

\begin{proof}
	Cette preuve se fait encore une fois par une suite d'égalités. Mentionnons toutefois l'argument suivant, qui est utile pour la preuve~: \textit{Si $\mu = \lambda+\vec{e_r}$ et $b<r$~:}
	\begin{equation}\label{eq:lem12}
		\projlift_b^\mu|_{S^\lambda} = 0.
	\end{equation}
	
	Pour démontrer cette affirmation, on se sert de deux résultats auxiliaires, la  \autoref{prop:lem45DS} et le \autoref{lem:lem43}, qui nous indiquent dans quels espaces les modules sont envoyés par la suite d'opérateurs. Graphiquement~:
	\[ S^\lambda \overset{\sh_b}{\longrightarrow} \bigoplus_{r\leq b}S^{\lambda+\vec{e_r}} \overset{\isoproj_{\mu}}{\longrightarrow} \left(\sum_{r\leq b} \delta_{\mu, \lambda+\vec{e_r}}\right) S^\mu.\] 
	Or, si $r>b$, comme dans l'hypothèse, l'image de $\projlift_b^{\lambda+\vec{e_r}}$ à partir de $S^\lambda$ est nulle.
	
	De plus, $\theta_{a,a}(w)$ est un multiple de l'identité. Le résultat est toujours le mot lui-même, et le nombre de façons de l'obtenir est le nombre d'occurrences de $a$ dans $w$. En particulier, $\theta_{a,a}|_{M^\lambda} = \lambda_a \cdot \Id$.

	Enfin,
	\allowdisplaybreaks
	\begin{align*}
		\big(\nu_k \circ &\projlift_a^\mu - \projlift_a^\mu \circ \nu_k \big)|_{S^\lambda} \\
		& = \left(\nu_k \circ \isoproj_\mu \circ \sh_a\right)|_{S^\lambda} - \left(\isoproj_\mu \circ \sh_a \circ \nu_k \right)|_{S^\lambda}\\
		& \eqexp{\autoref{lem:lem11_commut}}{ \isoproj_\mu\left(\nu_k \circ \sh_a - \sh_a \circ \nu_k\right)|_{S^\lambda}}\\
		& \eqexp{\autoref{thm:thm2}}{}\isoproj_{\mu} \left((n+2-k)(\sh_a \circ \Nu{k-1})|_{S^\lambda} + \sum_{1\leq b \leq a} \bbox{\sh_b\circ\theta_{b,a}}\circ\Nu{k-1}|_{S^\lambda}\right)\\
		& \eqexp{\autoref{lem:lem35}}{} (n+2-k) \projlift_a^\mu \circ \Nu{k-1}|_{S^\lambda} + \sum_{1\leq b \leq a} \bbox{\isoproj_\mu  \circ\theta_{b,a}}\circ \sh_b \circ\Nu{k-1}|_{S^\lambda}\\
		& \qquad \qquad - \sum_{1\leq b \leq a} \underbrace{\isoproj_\mu  \circ \sh_a}_{\projlift_a^\mu} \circ\Nu{k-1}|_{S^\lambda}\\
		& \eqexp{\autoref{lem:lem11_commut}}{} (n+2-k-a) \projlift_a^\mu \circ \Nu{k-1}|_{S^\lambda} + \sum_{\bbox{\scriptstyle 1\leq b \leq a}}  \theta_{b,a}\circ \projlift_b^\mu \circ\Nu{k-1}|_{S^\lambda}\\
		& \eqexp{\text{équation }\eqref{eq:lem12}}{} (n+2-k-a) \projlift_a^\mu \circ \Nu{k-1}|_{S^\lambda} + \sum_{r\leq b \leq a}  \theta_{b,a}\circ \projlift_b^\mu \circ\Nu{k-1}|_{S^\lambda}\\
		& =  (n+2-k-a) \projlift_a^\mu \circ \Nu{k-1}|_{S^\lambda} + \underbrace{\theta_{a,a}}_{(\lambda_a+1) \Id}\circ \projlift_a^\mu \circ\Nu{k-1}|_{S^\lambda}\\
		& \qquad\qquad + \sum_{r\leq b < a}  \theta_{b,a}\circ \projlift_b^\mu \circ\Nu{k-1}|_{S^\lambda}\\
		& = (n+3-k+\lambda_a-a) \projlift_a^\mu \circ \Nu{k-1}|_{S^\lambda} + \sum_{r\leq b < a}  \theta_{b,a}\circ \projlift_b^\mu \circ\Nu{k-1}|_{S^\lambda}. \qedhere
	\end{align*}
\end{proof}

On peut maintenant décrire le comportement de nos opérateurs sur les modules de Specht~: la prochaine équation vient simplifier celle du \autoref{lem:lem13}.
\begin{thm}\label{thm:thm3}
	Soit $\lambda = (\lambda_1, \ldots, \lambda_l) \vdash n$, $a \in [l+1]$ et $\mu = \lambda + \vec{e_r}$, avec $r \in [a]$. Alors, 
	\[\left(\nu_k \circ \projlift_a^\mu - \projlift_a^\mu \circ \nu_k \right)|_{S^\lambda} = (n+2-k + (\lambda_r + 1- r)) \left(\projlift_a^\mu \circ \Nu{k-1}\right)|_{S^\lambda}.\]
\end{thm}
	
Une grosse partie du travail nécessaire à la preuve de ce théorème a déjà été réalisé à ce point. En fait, le \autoref{lem:lem13} est l'ingrédient-clé de la preuve qui suit.
\begin{proof}
	Grâce au \autoref{lem:lem13}, on sait que 
	\begin{align*}
		\big(\nu_k \circ &\projlift_a^\mu - \projlift_a^\mu \circ \nu_k \big)|_{S^\lambda} =\\
		&(n+3-k+\lambda_a - a) \left(\projlift_a^\mu \circ \Nu{k-1}\right)|_{S^\lambda} + \left(\sum_{r \leq b < a} \theta_{b,a} \circ \projlift_b^\mu \circ \Nu{k-1}\right)|_{S^\lambda}.
	\end{align*}
	On distingue deux cas~:\\
	\underline{Si $r = a$~:} Alors le deuxième terme du côté droit de l'équation est nul (la somme est vide), et l'équation est simplement
	\begin{equation}
	\left(\nu_k \circ \projlift_r^\mu - \projlift_r^\mu \circ \nu_k \right)|_{S^\lambda} = (n+3-k+\lambda_r - r) \left(\projlift_r^\mu \circ \Nu{k-1}\right)|_{S^\lambda} \label{eq:Slambda->Slambda,a=r}.
	\end{equation}
	
	\underline{Sinon, $r < a$~:} On se sert de l'équation \eqref{eq:Slambda->Slambda,a=r}. On multiplie à gauche chaque côté de cette équation par $\theta_{r,a}$. On a alors~:
	\[\begin{aligned}
	(n&+3-k+\lambda_r-r) \theta_{r,a} \circ \projlift_r^\mu \circ \Nu{k-1}|_{S^\lambda}\\
	&\eqexp{équation \eqref{eq:Slambda->Slambda,a=r}}{}\bbox{\theta_{r,a} \circ \nu_k \circ \projlift_r^\mu}|_{S^\lambda} - \bbox{\theta_{r,a} \circ \projlift_r^\mu} \circ \nu_k |_{S^\lambda}\\
	& \eqexp{$\stackunder{\text{\autoref{lem:lem11_thet}}}{\text{\autoref{lem:lem11_commut}}}$}{} 	\nu_k \circ \isoproj_\mu \circ \bbox{\theta_{r,a} \circ \sh_r}|_{S^\lambda} - \isoproj_\mu \circ \bbox{\theta_{r,a} \circ \sh_a} \circ \nu_k |_{S^\lambda}\\
	&\eqexp{\autoref{lem:lem35}}{}	\nu_k \circ \isoproj_\mu \circ \left(\sh_r \circ \theta_{r,a}  + \sh_a\right)|_{S^\lambda} - \isoproj_\mu \circ \left(\sh_r \circ \theta_{r,a}  + \sh_a\right) \circ \nu_k|_{S^\lambda}\\
	& = \qquad\big( \nu_k \circ \projlift_r^\mu \circ \theta_{r,a} + \nu_k \circ \projlift_a^\mu\\
	&\qquad\qquad- \projlift_r^\mu \circ \bbox{\theta_{r,a} \circ \nu_k} - \projlift_a^\mu \circ \nu_k\big)|_{S^\lambda}\\
	& \eqexp{\autoref{lem:lem11_thet}}{} \left( \nu_k \circ \projlift_a^\mu - \projlift_a^\mu \circ \nu_k\right)|_{S^\lambda}\\
	&\qquad\qquad + \left(\nu_k \circ \projlift_r^\mu  - \projlift_r^\mu \circ \nu_k\right) \circ \theta_{r,a}|_{S^\lambda} \\
	&\eqexp{équation \eqref{eq:Slambda->Slambda,a=r}}{} \left( \nu_k \circ \projlift_a^\mu - \projlift_a^\mu \circ \nu_k\right)|_{S^\lambda}\\
	& \qquad\qquad + (n+3-k+\lambda_r-r)\projlift_r^\mu  \circ \Nu{k-1} \circ \theta_{r,a}|_{S^\lambda}. 
	\end{aligned}\]
	En comparant la première et la dernière ligne de cette suite d'équations, on trouve
	\allowdisplaybreaks
	\begin{align*}
	\big( \nu_k \circ& \projlift_a^\mu - \projlift_a^\mu \circ \nu_k\big)|_{S^\lambda} \\
	&= (n+3-k+\lambda_r-r) \left( \bbox{\theta_{r,a}\circ \projlift_r^\mu} \circ \Nu{k-1} - \projlift_r^\mu  \circ \Nu{k-1} \circ \theta_{r,a}\right)|_{S^\lambda}\\
	& \eqexp{\autoref{lem:lem11_commut}}{} (n+3-k+\lambda_r-r) \isoproj_\mu \left( \theta_{r,a} \circ \sh_r \circ \Nu{k-1} - \sh_r  \circ \bbox{\Nu{k-1} \circ \theta_{r,a}}\right)|_{S^\lambda}\\
	& \eqexp{\autoref{lem:lem11_thet}}{}  (n+3-k+\lambda_r-r) \isoproj_\mu \left( \bbox{\theta_{r,a} \circ \sh_r - \sh_r  \circ \theta_{r,a}}\right) \circ \Nu{k-1}|_{S^\lambda}\\
	& \eqexp{\autoref{lem:lem35}}{}  (n+3-k+\lambda_r-r) \isoproj_\mu \circ \sh_a \circ \Nu{k-1}|_{S^\lambda}\\
	&  = \qquad  (n+3-k+\lambda_r-r) \projlift_a^\mu  \circ \Nu{k-1}|_{S^\lambda},
	\end{align*}
	ce qui achève la démonstration du \autoref{thm:thm3}.	
\end{proof}

On peut déduire du \autoref{thm:thm3} le corollaire suivant sur les valeurs propres de $\Nu{k}$.
\begin{corl}\label{corl:corl4}
	Soit $\lambda\vdash n$ et soit $v \in S^\lambda$ un vecteur propre à la fois pour $\Nu{k}$ et $\Nu{k-1}$, associé aux valeurs propres $v_k$ et $v_{k-1}$, respectivement. Alors, soit
	\begin{itemize}
		\item $\projlift_a^{\lambda+\vec{e_r}}(v) = 0$
		\item[ou] 
		\item $\begin{aligned}
		\Nu{k} \circ \projlift_a^{\lambda+\vec{e_r}}(v)& = \projlift_a^{\lambda+\vec{e_r}}\left( v_k + (n+2-k+(\lambda_r+1-r)) v_{k-1}\right) (v)\\ & = \left( v_k + (n+2-k+(\lambda_r+1-r)) v_{k-1}\right) \projlift_a^{\lambda+\vec{e_r}}(v).
		\end{aligned}$
	\end{itemize}
\end{corl}
L'hypothèse faite ici qu'un vecteur soit à la fois un vecteur propre pour $\nu_k$ et $\Nu{k-1}$ peut sembler contraignante. Toutefois, on a démontré que les opérateurs $\allnuk$ admettent une base commune de vecteurs propres (\autoref{thm:vec_prop_communs_allnuk}). L'existence de vecteurs propres à la fois pour $\Nu{k}$ et $\Nu{k-1}$ est non seulement garantie par ce théorème; on sait en plus que, dans une certaine base, tous les vecteurs propres sont ainsi. Cela nous garantit que ce processus nous permettra de trouver \textit{toutes} les valeurs propres des opérateurs $\Nu{k}$. 

\section{Théorème principal}\label{sec:vp_preuve}
On souhaite maintenant calculer explicitement toutes les valeurs propres de $\Nu{k}$. On aura ainsi démontré le \autoref{thm:main}. Pour ce faire, on procède comme suit~:
\begin{enumerate}
	\item On décrit, à l'aide d'un résultat de \cite{RSW}, le noyau de chacun des opérateurs $\Nu{k}$ (\autoref{ssec:vect_pr_noyau}).
	\item On présente une procédure pour obtenir tous les autres vecteurs propres, à partir de ceux du noyau (\autoref{ssec:vect_pr_pas_noyau}).
	\item Pour ces vecteurs propres, on peut calculer la valeur propre associée à l'aide du \autoref{corl:corl4} (\autoref{ssec:vp_pas_noyau}).
\end{enumerate}

\subsection{Vecteurs propres du noyau}\label{ssec:vect_pr_noyau}
Le noyau des opérateurs $\allnuk$ est connu depuis \cite{RSW}, où il est l'objet du théorème VI.10.5. Ce que Reiner, Saliola et Welker ont démontré, c'est que l'opérateur $\Delta$ de Schützenberger connecte les tableaux de taille $n$ et type $j$  à ceux de taille $n-1$ et type $j-1$. Ils ont de plus démontré que les espaces propres des opérateurs $\allnuk$ sur les collections de $n$ et $n-1$ objets étaient reliés de la même façon. Les trois auteurs expliquent la construction du noyau en fonction des tableaux d'un certain type~:
\begin{thm}[Théorème VI.10.5 de \cite{RSW}]\label{thm:noyau_nu_k}
	Le noyau de $\Nu{k}$ est
\begin{equation*}
\ker(\Nu{k}) \cong \bigoplus_{\substack{\text{tableau  standard $t$ ,}\\ \type(t) < k}} S^{\mathrm{forme}(t)}. 
\end{equation*}
\end{thm}

Ce théorème est pratique parce qu'il nous permet non seulement de donner la dimension du noyau, mais aussi parce qu'il est la base de l'induction que nous faisons pour décrire les autres espaces propres. Malheureusement, le dernier théorème ne donne pas une base du noyau; si nous avions une telle base, ce qui est présenté dans les prochains paragraphes permettrait de décrire une base de tous les espaces propres.

\subsection{Vecteurs propres qui ne sont pas dans le noyau}\label{ssec:vect_pr_pas_noyau}
On étudie d'abord la structure de l'image de $\Nu{k}|_{S^\lambda}$. On sait déjà que $\Nu{k}$ est un élément de l'algèbre de groupe, d'où $\im(\Nu{k}|_{S^\lambda}) \subseteq S^\lambda$. On peut toutefois raffiner ce résultat.
\begin{prop}\label{prop:im_nu_k_specht}
	L'image de $\Nu{k}$ sur $S^\lambda$ est comprise dans la somme suivante~:
	\begin{align*}
	\im(\Nu{k}|_{S^\lambda}) \ &\subseteq \sum_{1 \leq a_1 \leq \ldots \leq a_k\leq n}  \im\left(\isoproj_\lambda \circ \sh_{a_k} \circ \ldots \circ \sh_{a_1}|_{S^{\lambda-\vec{e_{a_1}}-\ldots-\vec{e_{a_k}}}}\right)\\
	& \subseteq\ \sum_{1 \leq a \leq n} \im\left(\projlift_a^\lambda |_{S^{\lambda-\vec{e_a}}}\right).
	\end{align*}
\end{prop}

\begin{proof}
	On fait d'abord quelques observations sur $\nu_k$ qui sont utiles pour démontrer la proposition.
	\begin{enumerate}
		\item Soit $\pi_k$ l'opérateur de \bhr\ tel que $\nu_k$ est un multiple de $\pi_k^\top\circ\pi_k$ (pour plus de détails, voir la \autoref{ssec:1efamille_bhr}). Alors, l'image de $\nu_k$ est comprise dans l'image de $\pi_k^\top$. En effet, si un vecteur $\vec{v}$ est dans l'image de $\nu_k$, il existe $\vec{u}$ tel que $\vec{v} = \nu_k(\vec{u}) = \pi_k^\top (\pi_k(\vec{u}))$.
		\item Soit $w$ un mot de longueur $n$, et notons $v = w_1\ldots w_{n-k}$. Soit $u$ le mot trié de longueur $k$ formé des lettres $w_{n-k+1},\ldots, w_n$ avec répétitions, c'est-à-dire le mot formé des lettres de $w$ qui ne sont pas dans $v$, placées en ordre croissant. La définition de $\pi_k$ donnée à la \autoref{ssec:1efamille_bhr} nous indique qu'on peut écrire 
		\[\pi_k^\top (w) = \sum_{1\leq u_1 \leq \ldots \leq u_k\leq n} \sh_{u_k} \circ \ldots \circ \sh_{u_1}(v).\]
		\item De plus, $v$ est le mot obtenu de $w$ en retirant les $k$ dernières lettres. Dans \cite{DS}, l'opérateur $\proj_a(w)$ est défini ainsi~: si $a$ est la dernière lettre de $w$, alors $\proj_a(w) = w_1 \ldots w_{n-1}$; sinon, $\proj_a(w) = 0$.\footnote{On devine, par le nom de l'opérateur, que $\proj_a|_{S^\lambda}$ est une projection sur le sous-module $S^{\lambda-\vec{e_a}}$.} Dans cet article, il est démontré (dans la preuve de la proposition 53) que l'image de $\proj_a|_{S^\lambda}$ est comprise dans $S^{\lambda-\vec{e_a}}$.\\
		Remarquons que $v$ peut s'écrire comme $v = \proj_{w_{n-k+1}} \circ \ldots \circ \proj_{w_n} (w)$.
	\end{enumerate}
	
	En combinant les deuxième et troisième observations, on trouve que
	\[ \im(\pi_k^\top|_{S^\lambda}) \subseteq  \sum_{1\leq u_1 \leq \ldots \leq u_k\leq n} \sh_{u_k} \circ \ldots \circ \sh_{u_1}|_{S^{\lambda-\vec{e_{u_1}} - \ldots -\vec{e_{u_k}}}}.  \]
	
	Rappelons que $\pi_k$ est un élément de l'algèbre du groupe. Ainsi, son image, lorsque restreint à un module, est dans ce module.
	D'où~:
	\begin{align*}
		\im(\nu_k|_{S^\lambda})
		& \subseteq \im(\pi_k^\top |_{S^\lambda})\\
		& = \im\left(\isoproj_\lambda \circ\  \pi_k^\top \right)|_{S^\lambda}\\
		& \subseteq \sum_{1\leq u_1 \leq \ldots \leq u_k\leq n} \im \left(\isoproj_\lambda \circ \sh_{u_k} \circ \ldots \circ \sh_{u_1}\right)|_{S^{\lambda-\vec{e_{u_1}} - \ldots -\vec{e_{u_k}}}},
	\end{align*}
	ce qui correspond à la première inclusion à démontrer.
	
	De plus, la proposition 55 de \cite{DS} indique que, lorsque \mbox{$u_1 \leq u_2 \leq \ldots \leq u_k$} représentent des rangées dans lesquelles se trouvent des boîtes de $\lambda$, alors
	\begin{align*}
		\isoproj_\lambda \circ  \sh_{u_k}& \circ \ldots \circ \sh_{u_1}|_{S^{\lambda-\vec{e_{u_1}}-\ldots-\vec{e_{u_k}}}}\\ 
		& = \left(\isoproj_\lambda \circ  \sh_{u_k}\right) \circ \left(\isoproj_{\lambda-\vec{e_{u_k}}} \circ  \sh_{u_{k-1}}\right)  \\
		&\quad\circ \ldots \circ\left(\isoproj_{\lambda-\vec{e_{u_2}}-\ldots-\vec{e_{u_k}}}\circ \sh_{u_1}|_{S^{\lambda-\vec{e_{u_1}}-\ldots-\vec{e_{u_k}}}}\right)
	\end{align*}
	Ainsi,
	\begin{align*}
		\im(\nu_k|_{S^\lambda}) & \subseteq \sum_{1\leq u_1 \leq \ldots \leq u_k\leq n} \im \left(\isoproj_\lambda \sh_{u_k} \circ \ldots \circ \sh_{u_1}\big|_{S^{\lambda-\vec{e_{u_1}} - \ldots -\vec{e_{u_k}}}}\right)\\
		& \subseteq \sum_{1\leq u_1 \leq \ldots \leq u_k\leq n} \im \big(\isoproj_\lambda \circ \sh_{u_k} \circ \isoproj_{\lambda-\vec{e_{u_k}}}\circ \sh_{u_{k-1}}\\
		&\qquad\qquad\qquad\qquad \circ\ldots \circ \sh_{u_1}\big|_{S^{\lambda-\vec{e_{u_1}} - \ldots -\vec{e_{u_k}}}}\big)\\
		 & \subseteq \sum_{1\leq u_k\leq n} \im \big(\underbrace{\isoproj_\lambda \circ \sh_{u_k}}_{\projlift_{u_k}^\lambda} \big|_{S^{\lambda-\vec{e_{u_k}}}}\big).\qedhere
	\end{align*}
\end{proof}

D'après le résultat de la dernière proposition, tout vecteur de $S^\lambda$ qui est dans l'image de $\nu_k$ peut s'écrire comme
\[ \sum_{1 \leq a \leq n}  \projlift_a^\lambda(\vec{v_a}), \]
où $\vec{v_a}$ est un vecteur de $S^{\lambda-\vec{e_a}}$. En particulier, si $\vec{v} \in S^\lambda$ est un vecteur propre de $\nu_k$ et de $\nu_{k-1}$ qui n'est pas dans le noyau, alors il existe des vecteurs propres $\vec{v_1} \in S^{\lambda-\vec{e_1}}$, $\ldots$, $\vec{v_n} \in S^{\lambda-\vec{e_n}}$ de $\nu_k$ et de $\Nu{k-1}$ tels que 
\[\vec{v} =  \sum_{1 \leq a \leq n} \delta_a \projlift_a^\lambda(\vec{v_a}), \]
où $\delta_a$ peut valoir $0$ ou $1$.
Autrement dit, les vecteurs propres de $\Nu{k}$ qui ne sont pas dans le noyau sont tous construits à partir de vecteurs propres de $\Nu{k}$ pour un mot dont la longueur est inférieure (de $1$) et avec l'opérateur $\projlift$.

De la même façon, pour tout $\vec{v_a} \in S^{\lambda - \vec{e_a}}$ qui est un vecteur propre pour les opérateurs $\Nu{k}$ et $\Nu{k-1}$, $\projlift_a^\lambda(\vec{v_a})$ est un vecteur propre de $\Nu{k}$ (\autoref{corl:corl4}). De plus, le corollaire donne deux choix pour la valeur propre associée à $\projlift_a^\lambda(\vec{v_a})$~: elle est soit nulle, soit entièrement déterminée par les valeurs propres de $\vec{v_a}$, par la lettre ($a$) qui est ajoutée au mot et par $k$ (les détails sont dans l'énoncé du \autoref{corl:corl4}).

\subsection{Valeurs propres}\label{ssec:vp_pas_noyau}

Les sous-sections \ref{ssec:vect_pr_noyau} et \ref{ssec:vect_pr_pas_noyau} décrivent certaines propriétés des vecteurs propres. Avec celles-ci, on a suffisamment d'informations pour donner toutes les valeurs propres des opérateurs $\allnuk$.

On rappelle de la \autoref{fig:recap_vp_tableaux} que les valeurs propres distinctes des opérateurs $\allnuk$ correspondent aux tableaux standards de taille $n$. Pour un tableau $t$, la multiplicité de la valeur propre associée est la dimension du module de Specht $S^{\mathrm{forme}(t)}$, soit le nombre de tableaux standards qui ont la même forme que $t$.

\begin{tcolorbox}
	\centering
	Tableaux standards $\longrightarrow$ Valeurs propres 
\end{tcolorbox} 

Le sujet de la \autoref{ssec:vect_pr_noyau} est le noyau des opérateurs $\allnuk$. Grâce à cette section, on connaît la multiplicité de la valeur propre $0$ et on sait combien de copies de chaque module de Specht sont dans le noyau~:
\begin{tcolorbox}
	\centering
	Tableaux standards de type inférieur à $k-1$ $\rightarrow$ Valeur propre nulle pour $\Nu{k}$
\end{tcolorbox}

Les deux résultats plus haut nous disent que les valeurs propres non-nulles de $\Nu{k}$ sont données par les tableaux standards de type supérieur ou égal à $k$.

On a vu plus tôt que les vecteurs propres de $\nu_k$ sur le module de Specht $S^\lambda$ qui ne sont pas dans le noyau s'écrivent tous comme la somme $\sum_{1 \leq a \leq n} \delta_a\projlift_a(\vec{v_a})$, où $\vec{v_a}$ est un vecteur propre situé dans $S^{\lambda-\vec{e_a}}$ et $\delta_a$ vaut soit $0$ ou $1$. On peut donc dire que $\projlift_a^\lambda$ permet d'obtenir des vecteurs propres de $S^\lambda$ à partir de ceux de $S^{\lambda-\vec{e_a}}$. Or, on sait d'après la règle de branchement (\autoref{ssec:regle_branchement}) que 
\[ S^\lambda \downarrow_{S_{n-1}} \cong \bigoplus_{1\leq a \leq n} S^{\lambda-\vec{e_a}}. \]
Donc, pour associer valeurs propres et tableaux standards, 
\begin{itemize}
	\item on sait que, si la valeur propre est nulle, on peut associer un tableau de type (strictement) inférieur à $k$;
	\item sinon, une valeur propre de $S^\lambda$ doit être associée à un tableau $t$ de forme $\lambda \vdash n$. De plus, il existe une certaine valeur $a$ telle que $\lambda - \vec{e_a}$ est un partage et $t$ est lié à un tableau $t'$ de forme $\lambda - \vec{e_a}$ de la façon suivante~: la valeur propre de $\Nu{k}$ associée à $t$ est 
	\[v_k(t) = v_k(t') + (n+1-k+\lambda_a - a)v_{k-1}(t').\]
	Cette dernière condition découle du \autoref{corl:corl4}.
	\item la procédure qui associe à $t$ exactement un tel tableau est l'opérateur $\Delta$ de Schützenberger, présenté à la \autopageref{sssec:op_delta_schutz}. Ainsi, à $t$, on associe un tableau $t' = \Delta(t)$, et la rangée dans laquelle se trouve l'unique boîte de $t$ qui n'apparaît pas dans $t'$ est notée $a$; on peut donc calculer la valeur propre associée à $t$. Réciproquement, à partir de $t'$, on peut construire un tableau de taille $n$ pour chaque rangée dans laquelle on peut ajouter une boîte  et obtenir un diagramme. En spécifiant qu'on veut ajouter une boîte dans la rangée $a$, on retrouve le tableau $t$ (le détail de la procédure inverse de $\Delta$ se trouve à la \autoref{rmq:inv_schutz_delta}).
	
	Comme le type de $\Delta(t)$ est strictement inférieur à celui de $t$, on sait qu'en un nombre fini d'étapes, on atteint la valeur propre $0$.
\end{itemize}
Ceci achève de démontrer le \autoref{thm:main}.

\begin{tcolorbox}
	\centering
	  Tableaux standards de type au moins $k$ $\rightarrow$ Valeurs propres non-nulles
\end{tcolorbox}

\section{Commutativité des opérateurs}\label{sec:commutativite}
Le but de cette section est de donner une nouvelle preuve à un résultat de Victor Reiner, Franco Saliola et Volkmar Welker, affirmant que les opérateurs $\allnuk$ commutent entre eux (Théorème I.1.1 de \cite{RSW}). Dans le même article, ils invitent le lectorat à chercher une preuve alternative de leur résultat (problème VI.1.4 du même article). En effet, leur preuve se fait par récurrence et compare différents cas, ce qui rend la lecture plutôt technique. Sans prétendre que la preuve donnée ici soit nécessairement plus éclairante, elle présente une façon alternative de comprendre le résultat.

\subsection{Une autre identité sur l'algèbre des mots}
Afin de démontrer la commutativité des opérateurs $\allnuk$, nous avons besoin d'un résultat concernant le retrait d'une lettre dans un mot et son interaction avec les opérateurs $\allnuk$. C'est un outil très similaire au \autoref{thm:thm1}, qui comparait l'effet de l'ajout d'une lettre avant ou après avoir mélangé, et les techniques de preuve utilisées sont les mêmes.
\begin{thm}\label{thm:thm5}
	La propriété suivante concernant l'opérateur de suppression est vraie sur l'algèbre des mots de longueur $n$~:
	\[ \del_a \circ \nu_k - \nu_k \circ \del_a = (n+1-k) \Nu{k-1} \circ \del_a + \sum_{1\leq b \leq n} \theta_{a,b} \circ  \Nu{k-1} \circ \del_b. \]
\end{thm}
\begin{proof}
	La preuve se fait par induction sur $k$. Le cas de base est celui du mélange doublement aléatoire, soit lorsque $k=1$. Alors,
	\begin{align*}
		\del_a \circ \nu_1 -  \nu_1 \circ \del_a & \eqexp{\autoref{prop:defn_op}}{} \del_a \circ \sum_{1 \leq b\leq n} \sh_b \circ \del_b - \sum_{1 \leq b\leq n} \sh_b \circ \del_b \circ \del_a\\
		& \eqexp{\autoref{lem:lem35}}{} \sum_{1 \leq b\leq n} \left( \del_{a} \circ \sh_b \circ \del_b - (\del_a \circ \sh_b - \theta_{a,b} - n \delta_{a,b}) \del_b \right)\\
		& \qquad =  \sum_{1 \leq b\leq n} \theta_{a,b} \circ \del_b + n \del_a.
	\end{align*}
	On pose l'hypothèse d'induction suivante~:  pour une certaine valeur de $k \geq 0$, et sur les mots de longueur $n-1$,
	\[\del_a \circ \nu_k - \nu_k \circ \del_a = (n-k) \Nu{k-1} \circ \del_a + \sum_{1\leq b \leq n-1} \theta_{a,b} \circ  \Nu{k-1} \circ \del_b.\]
	
	L'étape d'induction va comme suit~: 
	\allowdisplaybreaks
	\begin{align*}
		\del_a& \circ \nu_{k+1} - \nu_{k+1} \circ \del_a\\
		& \eqexp{\autoref{thm:lem8}}{} \del_a \circ \nu_{k+1}- \frac{1}{k+1} \sum_{1\leq b \leq n} \left(\sh_b \circ \Nu{k} \circ \del_b\right) \circ \del_a\\
		& \eqexp{\autoref{prop:commutShDel}}{} \del_a \circ \nu_{k+1}- \frac{1}{k+1} \sum_{1\leq b \leq n} \sh_b \circ \Nu{k} \circ \del_a \circ \del_b\\
		&\eqexp{Hyp. induction}{} \del_a \circ \nu_{k+1}- \frac{1}{k+1} \sum_{1\leq b \leq n} \sh_b \circ \bigg(\del_a \circ \Nu{k} - (n-k) \Nu{k-1} \circ \del_a\\
		& \qquad \qquad - \sum_{1\leq c \leq n-1} \theta_{a,c} \circ \Nu{k-1} \circ \del_c \bigg) \circ \del_b\\	
		& \eqexp{\autoref{prop:commutShDel}}{} \del_a \circ \nu_{k+1}- \frac{1}{k+1} \bigg(\sum_{1\leq b \leq n} \bbox{\sh_b \circ \del_a }\circ \Nu{k} \circ  \del_b\\
		& \qquad \qquad  - (n-k) \underbrace{\sum_{1\leq b \leq n} \sh_b \circ \Nu{k-1} \circ \del_b}_{k \cdot \Nu{k}} \circ  \del_a - \sum_{1\leq b,c \leq n} \bbox{\sh_b \circ \theta_{a,c}} \circ \Nu{k-1} \circ \del_c  \circ \del_b\bigg)\\	
		& \eqexp{\autoref{lem:lem35}}{} \del_a \circ \nu_{k+1}- \frac{1}{k+1} \bigg(\sum_{1\leq b \leq n} (\del_a \circ \sh_b -\theta_{a,b} - n \delta_{a,b}) \circ \Nu{k} \circ  \del_b\\
		& \qquad \qquad  - (n-k) k \cdot \Nu{k} \circ  \del_a - \sum_{1\leq b,c \leq n} \left(\theta_{a,c} \circ \sh_b  - \delta_{a,b} \sh_c\right) \circ \Nu{k-1} \circ \bbox{\del_c  \circ \del_b}\bigg)\\	
		& \eqexp{\autoref{prop:commutShDel}}{} \del_a \circ \nu_{k+1}- \frac{1}{k+1} \bigg(\del_{a} \circ \underbrace{\sum_{1\leq b \leq n} \sh_b\circ \Nu{k} \circ  \del_b}_{(k+1)\cdot \Nu{k+1}} - \sum_{1\leq b \leq n}\theta_{a,b} \circ \Nu{k} \circ  \del_b - n\  \Nu{k} \circ  \del_a\\
		& \qquad \qquad  - (n-k) k \cdot \Nu{k} \circ  \del_a - \sum_{1\leq c \leq n}\theta_{a,c} \circ \underbrace{\left( \sum_{1\leq b \leq n} \sh_b \circ \Nu{k-1} \circ \del_b \right)}_{k \cdot \Nu{k}}  \circ \del_c\\
		& \qquad \qquad  +\underbrace{ \sum_{1\leq c \leq n} \sh_c \circ \Nu{k-1} \circ \del_c}_{k\cdot \Nu{k}}  \circ \del_a \bigg)\\
		& \eqexp{\autoref{thm:lem8}}{} \frac{1}{k+1} \bigg( \sum_{1\leq b \leq n}\theta_{a,b} \circ \Nu{k} \circ  \del_b + (n+(n-k)k) \cdot \Nu{k} \circ  \del_a\\
		& \qquad \qquad  + k \sum_{1\leq c \leq n}\theta_{a,c} \circ \Nu{k}  \circ \del_c - k\cdot \Nu{k}  \circ \del_a \bigg)\\
		& =  \qquad \frac{1}{k+1} \left((n-k)(k+1) \Nu{k}\circ\del_a + (k+1)   \sum_{1\leq b \leq n}\theta_{a,b} \circ \Nu{k} \circ  \del_b \right)\\
		& = \qquad (n-k) \Nu{k} \circ \del_a +  \sum_{1\leq b \leq n}\theta_{a,b} \circ \Nu{k} \circ  \del_b,
	\end{align*}
	ce qui termine la preuve.
\end{proof}

\subsection{Nouvelle preuve de commutativité}
Une des conséquences des résultats démontrés jusqu'ici est que les opérateurs $\allnuk$ commutent entre eux. Ceci était déjà connu \cite[Théorème I.1.1]{RSW}, mais le résultat apparaît de façon plus directe ici.

\paragraph{Théorème \ref{thm:commutativite_nu}.}
Les opérateurs $\allnuk$ commutent deux-à-deux.

La preuve qui suit n'est pas la seule qui existe; \cite{RSW} en contient également une. Toutefois, leur preuve se termine en mettant au défi son lectorat de trouver une autre preuve, plus directe. Avec les théorèmes précédents, il est plus direct de montrer la commutativité. On le fait ici par induction.

\begin{proof}[Preuve du \autoref{thm:commutativite_nu}.]
	Il suffit de montrer que, pour toutes valeurs de $j$ et $k$, $\Nu{j}\circ \Nu{k} = \Nu{k} \circ \Nu{j}$.
	On fait la preuve en fixant $k$ et par induction sur $j$.
	
	\underline{Cas de base : $j=1$}\\
	Remarquons d'abord que  la \autoref{prop:defn_op} et le \autoref{thm:lem8} impliquent que
    \[ \nu_1 \circ \nu_k - (k+1) \Nu{k+1}  = \sum_{1 \leq a \leq n} \left(\sh_a \circ \del_a \circ \Nu{k} - \sh_a \circ \nu_k \circ  \del_a\right).  \]
	Ainsi, sur les mots de longueur $n$, \allowdisplaybreaks
	\begin{align*}
		\nu_1 \circ \nu_k
		& \eqexp{}{} (k+1) \Nu{k+1} + \sum_{1 \leq a \leq n} \sh_a \left(\bbox{\del_a \circ  \Nu{k} - \nu_k \circ  \del_a}\right)\\
		& \eqexp{\autoref{thm:thm5}}{} (k+1) \Nu{k+1} \\
		& \qquad \qquad+ \sum_{1 \leq a \leq n} \sh_a \circ \left( (n+1-k) \Nu{k-1} \circ \del_a + \sum_{1\leq b \leq n} \theta_{a,b} \circ  \Nu{k-1} \circ \del_b \right)\\
		& = \qquad (k+1) \Nu{k+1} +(n+1-k)  \underbrace{\sum_{1 \leq a \leq n} \sh_a \circ \Nu{k-1} \circ \del_a}_{k\cdot \Nu{k}} \\
		&\qquad \qquad + \sum_{1\leq a,b \leq n} \sh_a \circ \theta_{a,b} \circ  \Nu{k-1} \circ \del_b \\
		& \eqexp{\autoref{thm:lem8}}{} (k+1) \Nu{k+1} +(n+1-k)k\cdot \Nu{k} + \sum_{1\leq a,b \leq n} \sh_a \circ \theta_{a,b} \circ  \Nu{k-1} \circ \del_b.
	\end{align*}
	
	De plus,\footnote{Dans la deuxième égalité, $\nu_k\circ \sh_a$ agit sur les mots de longueur $n-1$ et le \autoref{thm:thm1} est ajusté en conséquences.}
	\allowdisplaybreaks
	\begin{align*}
		\Nu{k} \circ \Nu{1} & \eqexp{\autoref{prop:defn_op}}{} \Nu{k} \circ \sum_{1 \leq a \leq n} \sh_a \circ \del_a\\
		& \eqexp{\autoref{thm:thm1}}{}  \sum_{1 \leq a \leq n} \left( \sh_a \circ \Nu{k} + (n+1-k) \sh_a \circ \Nu{k-1} + \sum_{1 \leq b\leq n} \sh_b \circ \theta_{b,a} \circ \Nu{k-1} \right) \circ \del_a\\
		& = \qquad \underbrace{\sum_{1 \leq a \leq n} \sh_a \circ \Nu{k} \circ \del_a}_{(k+1)\Nu{k+1}} + (n+1-k) \underbrace{\sum_{1 \leq a \leq n} \sh_a \circ \Nu{k-1}\circ \del_a}_{k\cdot \Nu{k}}\\
		&\qquad \qquad + \sum_{1 \leq a,b\leq n} \sh_b \circ \theta_{b,a} \circ \Nu{k-1}  \circ \del_a\\
		& \eqexp{\autoref{thm:lem8}}{} (k+1)\Nu{k+1} + (n+1-k) k\cdot \Nu{k} + \sum_{1 \leq a,b\leq n} \sh_b \circ \theta_{b,a} \circ \Nu{k-1}  \circ \del_a,
	\end{align*}
	
	Ceci permet de démontrer l'égalité de $\Nu{1}\circ \Nu{k} = \Nu{k} \circ \Nu{1}$.
	
	\underline{Hypothèse d'induction~:}\\
	Supposons que, pour toute valeur $i \leq j$, $\nu_i\circ  \nu_k = \nu_k \circ \Nu{i}$.
	
	\underline{Étape d'induction~:}\\ 
	Il reste à montrer que $\Nu{j+1} \circ \Nu{k} = \Nu{k} \circ \Nu{j+1}$~:
	\allowdisplaybreaks
	\begin{align*}
		\Nu{j+1} \circ \Nu{k} & \eqexp{\autoref{thm:lem8}}{} \frac{1}{j+1} \bigg( \sum_{1 \leq a \leq n} \sh_a \circ \Nu{j} \circ \bbox{\del_a  \bigg)\circ \Nu{k}}\\
		& \eqexp{\autoref{thm:thm5}}{} \frac{1}{j+1} \sum_{1 \leq a \leq n} \sh_a \circ \bbox{\Nu{j} \circ \bigg( \Nu{k}} \circ \del_{a}  + (n+1-k) \bbox{\Nu{k-1}} \circ \del_a\\
		&\qquad\qquad + \sum_{1 \leq b\leq n} \bbox{\theta_{a,b} \circ \Nu{k-1}}\circ \del_b\bigg)\\
		& \eqexp{$\stackunder{\text{Hyp. induction et}}{\text{\autoref{lem:lem11_thet}}}$}{} \frac{1}{j+1} \bigg( \sum_{1 \leq a \leq n} \bbox{\sh_a \circ \Nu{k}} \circ \Nu{j} \circ \del_{a}\\
		&\qquad\qquad + (n+1-k) \sum_{1\leq a\leq n} \sh_a \circ \Nu{k-1} \circ \Nu{j} \circ \del_a\\
		& \qquad \qquad + \sum_{1 \leq a,b\leq n} \sh_a \circ \theta_{a,b} \circ \Nu{k-1}\circ \Nu{j}\circ \del_b\bigg)\\
		& \eqexp{\autoref{thm:thm1}}{}  \frac{1}{j+1} \bigg( \sum_{1 \leq a \leq n} \bigg(\Nu{k} \circ \sh_a - (n+1-k) \sh_a \circ \Nu{k-1}\\
		&\qquad\qquad - \sum_{1\leq b \leq n} \sh_b \circ \theta_{b,a} \circ \Nu{k-1}\bigg) \circ \Nu{j} \circ \del_{a} \\
		& \qquad \qquad + (n+1-k) \sum_{1\leq a\leq n} \sh_a \circ \Nu{k-1} \circ \Nu{j} \circ \del_a\\
		&\qquad\qquad + \sum_{1 \leq a,b\leq n} \sh_a \circ \theta_{a,b} \circ \Nu{k-1}\circ \Nu{j}\circ \del_b\bigg)\\
		& \eqexp{}{}  \frac{1}{j+1} \Nu{k} \circ \sum_{1 \leq a \leq n} \sh_a \circ \Nu{j} \circ \del_{a} \\
		& \eqexp{\autoref{thm:lem8}}{}\nu_k \circ \Nu{j+1}.
	\end{align*}
	On conclut ainsi que les opérateurs $\allnuk$ commutent entre eux.
\end{proof}

\chapter{Une deuxième famille d'opérateurs qui commutent}\label{chap:2efamille}

Jusqu'à présent, on a défini les opérateurs de mélange symétrisés comme une généralisation du mélange doublement aléatoire dans laquelle on retirait plus d'un objet à chaque itération. On peut cependant aller beaucoup plus loin dans la généralisation~: dans leur travail précurseur sur les opérateurs de mélange symétrisés, Victor Reiner, Franco Saliola et Volkmar Welker ont défini ceux-ci pour des groupes de réflexion. Nous nous en tiendrons ici au cas du groupe symétrique, mais nous examinerons en particulier un ensemble d'opérateurs qui ont des propriétés similaires à ceux de la première famille. Nommément, ils commutent entre eux et ont des valeurs propres entières. Cela les distingue des autres opérateurs définis dans \cite{RSW} pour le cas du groupe symétrique.

On décrit spécifiquement les opérateurs $\allgammak$, ceux qui nous intéressent, dans la section qui suit. Pour avoir tout le matériel nécessaire pour comprendre la multitude de définitions équivalentes, on s'attarde brièvement à définir les autres opérateurs de mélange sur le groupe symétrique (nommés plus loin $\nu_\lambda$) à la \autoref{ssec:gammak_bhr}. Ceux-ci ne sont cependant pas l'objet de ce chapitre et nous les délaissons rapidement. Après avoir décrit les matrices des opérateurs $\allgammak$ et leurs propriétés, nous nous attardons à leurs valeurs propres et présentons des conjectures pour certaines d'entre elles.
\section{Définitions équivalentes}
On présente quatre définitions des opérateurs $\allgammak$. La première s'exprime en termes de mélange d'une collection d'éléments distincts; comme plus tôt, les objets représentés peuvent être des cartes d'un paquet de cartes à jouer. Ensuite, nous rappelons que nous étudions les versions symétrisées de certains opérateurs de \bhr\ et détaillons lesquels. Puis, nous nous intéressons à l'interprétation en termes de suites croissantes de tous les opérateurs $\{\nu_\lambda\}_{\lambda \vdash n}$, en spécifiant ce que cela veut dire pour la famille $\allgammak$. Enfin, nous examinons les différentes propriétés des matrices associées à $\allgammak$.
\subsection[Comme personne ne mélangerait]{Comme \sout{une} personne (ne) mélangerait}\label{ssec:gammak_intuitif}\index{Opérateurs!$\gamma_k$}
Tout comme les opérateurs $\allnuk$, les opérateurs $\allgammak$ sont d'abord des opérateurs de mélange. On peut ainsi les expliquer à quelqu'un qui souhaiterait mélanger une collection d'objets. Pour prendre l'exemple le plus simple, considérons un paquet de cartes. Pour le mélanger, on le divise en $k$ plus petits paquets de deux cartes et $n-2k$ paquets d'une seule carte. Ensuite, on combine de nouveau tous ces paquets de la façon suivante~: on bat le premier paquet avec le deuxième, puis le paquet résultant avec le troisième, et ainsi de suite jusqu'à ce que tous les paquets soient fusionnés. Pour battre deux paquets, on intercale leurs cartes en préservant l'ordre des cartes dans chacun des deux plus petits paquets, comme illustré sur la \autoref{fig:melange_americain}.
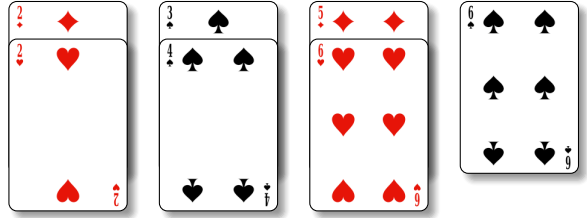
\begin{figure}
	\centering
	\begin{tikzpicture}
	\node (A) at (0,0) {\carda};
	\node (B) at (0,-0.5) {\cardb};
	\node (C) at (2,0) {\cardc};
	\node (D) at (2,-0.5) {\cardd};
	\node (E) at (4,0) {\carde};
	\node (F) at (4,-0.5) {\cardf};
	\node (G) at (6,0) {\cardg};
	\end{tikzpicture}
	\caption{Paquets d'une et de deux cartes pour l'exécution du mélange $\gammak{3}$.}
\end{figure}
L'opérateur $\gammak{k}$ attribue à la permutation $\sigma$ la combinaison linéaire formée de toutes les permutations $\tau$ obtenues de cette façon, auxquelles on affecte, comme coefficient, le nombre de façons de les obtenir.

\begin{ex}
	Appliquer $\gammak{2}$ revient à séparer le paquet de quatre cartes en deux paquets de deux cartes, puis à les battre. Ainsi, $\gammak{2}(\mathsf{1234}) = 3 \cdot \mathsf{1234} + 2 \cdot \mathsf{1243} + 2 \cdot \mathsf{1324} +  \mathsf{1423} + \mathsf{1342} + 2  \cdot \mathsf{2134} + 2  \cdot \mathsf{2143} + \mathsf{2314} + \mathsf{2413} + \textsf{3124} + \textsf{3142} + \textsf{3412}$.
	Par exemple, pour obtenir le résultat $\textsf{2143}$, on peut soit séparer $\textsf{1}$ et $\textsf{3}$ de $\textsf{2}$ et $\textsf{4}$, ou encore placer $\textsf{1}$ et $\textsf{4}$ ensemble, isolés de $\textsf{2}$ et $\textsf{3}$. C'est pourquoi on retrouve $2$ comme coefficient devant $\textsf{2143}$.
\end{ex}

\subsection{Opérateurs de BHR symétrisés}\label{ssec:gammak_bhr}
Les opérateurs de la famille $\allgammak$ sont également des multiples des versions symétrisées de certains opérateurs de \bhr. De la même façon qu'on pouvait obtenir $\Nu{k}$  en séparant les cartes pour obtenir un paquet de taille $n-k$ et des paquets de taille $1$, puis en battant les cartes, on peut séparer notre paquet en $k$ paquets de taille $2$ et $n-2k$ paquets de taille 1, puis les recombiner pour obtenir $\gammak{k}$. Or, en séparant puis recombinant les paquets, on applique en fait un opérateur de $\bhr$ et son opérateur dual~:
\begin{itemize}
	\item l'opérateur de \bhr\ $\pi_\lambda$\index{Opérateurs!$\pi_\lambda$}, défini comme le produit par la somme de compositions $\sum_{c \text{ de forme }\lambda} c$, divise le paquet en sous-ensembles de taille $\lambda_1, \ldots, \lambda_{|\lambda|}$, puis les empile dans cet ordre.
	\item son opérateur dual (dont la matrice est $\pi_\lambda^\top$) divise le paquet de la façon suivante~: il retire les $\lambda_1$ premières cartes et en fait un nouveau paquet, il fait de même avec les $\lambda_2$ cartes suivantes, et ainsi de suite jusqu'à la fin du paquet. Ensuite, les paquets sont recombinés par battage.
\end{itemize}

Ainsi, l'opérateur $\gammak{k}$ est \textit{presque} $\pi_{(2^k,1^{n-2k})}^\top \circ \pi_{(2^k,1^{n-2k})}$. 

\begin{rmq}
	Afin d'alléger la notation, les partages avec $k$ parts de taille $2$ et tels que toutes les autres parts sont des singletons sont notés $(2^k,1^{n-2k})$. Bien qu'il puisse y avoir une certaine ambiguïté avec le partage à deux parts, une de taille $2^k$ et une taille $1$, le contexte fera que nous parlerons toujours de partage avec des parts de taille $1$ ou $2$, ce qui est d'ailleurs un partage de $n$.
\end{rmq}

\begin{defn}
	L'opérateur $\gamma_k$ est un multiple de la version symétrisée d'un opérateur de \bhr. Plus précisément,
	\[ \gammak{k} = \frac{1}{k!(n-2k)!}\pi_{(2^k, 1^{n-2k})}^\top \circ \pi_{(2^k,1^{n-2k})}, \]
	où $\pi_{(2^k,1^{n-2k})}$ est la multiplication par la somme des compositions de forme $(2^k,1^{n-2k})$.\index{Opérateurs!$\gamma_k$}
\end{defn}

\begin{proof}[Preuve de l'équivalence des définitions]
	On a déjà remarqué que $\gammak{k}$ est un multiple de $\pi_{(2^k, 1^{n-2k})}^\top \circ \pi_{(2^k,1^{n-2k})}$. La présence d'un dénominateur s'explique par le fait que l'ordre dans lequel les paquets de même taille sont séparés pour être recombinés n'importe pas dans $\gammak{k}$. Il est cependant pris en compte dans le produit de compositions. En effet, deux compositions sont différentes si les blocs sont les mêmes, mais que des blocs de même taille ont été permutés.
	
	Pour calculer le dénominateur, il suffit de comparer la somme des coefficients dans $\gammak{k}(\Id)$ tel que défini à la \autoref{ssec:gammak_intuitif} et celle dans $\pi_{(2^k, 1^{n-2k})}^\top \circ \pi_{(2^k, 1^{n-2k})}(\Id)$. Remarquons que, dans ce dernier cas, la somme des coefficients est le carré de celle dans $\pi_{(2^k, 1^{n-2k})}(\Id)$.
	
	La somme des coefficients dans $\gammak{k}(\Id)$ est le nombre de façons de faire ce qui suit~:
	\begin{enumerate}
		\item Faire $k$ paquets de taille $2$, et $n-2k$ paquets de taille $1$. Comme on ne souhaite pas tenir compte de l'ordre des paquets identiques, il y a $\frac{\binom{n}{(2^k, 1^{n-2k})}}{k! (n-2k!)}$ façons de faire ceci.
		\item Battre successivement les paquets. Pour chaque paquet de $l$ cartes battu avec un paquet de $m$ cartes, on trouve $\binom{l+m}{l}$ recombinaisons possibles. Ainsi, il y a 
		\[ \binom{4}{2} \cdot \binom{6}{2}  \cdot \cdots \cdot \binom{2k}{2} \cdot (2k+1)\cdot \ldots \cdot n = \frac{n!}{2^k} = \binom{n}{(2^k, 1^{n-2k})}  \]
		façons de recombiner les paquets.
		\item Pour obtenir la somme des coefficients, on multiplie le nombre de résultats possibles pour chacune des deux étapes~:
		\[ \frac{\binom{n}{(2^k, 1^{n-2k})}^2}{k!(n-2k)!}. \]
	\end{enumerate}
	En revanche, la somme des coefficients dans $\pi_k(\Id)$ correspond au nombre de termes dans $\sum_{c\text{ de forme }(2^k,1^{n-2k})} c$, et donc au nombre de compositions ensemblistes de cette forme. Il y en a $\binom{n}{(2^k,1^{n-2k})}$.
	On peut conclure que la somme des coefficients dans $\pi_{(2^k,1^{n-2k})}^\top \circ \pi_{(2^k,1^{n-2k})}(\Id)$ est $\binom{n}{(2^k,1^{n-2k})}^2$. C'est $k! (n-2k)!$ fois plus que dans la définition de $\gammak{k}$.
\end{proof}

\subsection{De multiples suites croissantes}\label{ssec:suites_croissantes_gamma_k}
De la même façon qu'on peut écrire les opérateurs $\allnuk$ en termes de suites croissantes, on peut le faire pour les opérateurs $\allgammak$. De manière plus générale, on pourrait le faire pour des opérateurs de mélange symétrisés associés à un partage.

\paragraph{Les opérateurs $\Nu{\lambda}$ sur le groupe symétrique}\index{Opérateurs!$\nu_\lambda$}
À partir de la définition des opérateurs comme la version symétrisée d'opérateurs de \bhr, on peut généraliser la définition des opérateurs $\nu_k$ en termes de suites croissantes. Ceci a été fait par Victor Reiner, Franco Saliola et Volkmar Welker; leur travail est toutefois beaucoup plus général, puisqu'il porte sur un grand nombre de groupes de réflexions, et pas seulement sur le groupe symétrique. En se restreignant au groupe symétrique, on trouve la définition suivante~:
\begin{defn}\label{defn:entrees_mat_v_lambda}
	Soit $\lambda$ un partage de $n$ et soit $\pi_\lambda$ l'opérateur de \bhr\ correspondant à la multiplication par la somme de compositions $\sum_{c \text{ de forme } \lambda} c$. 
	Alors, l'entrée $(\sigma, \tau)$ de la matrice $\nu_\lambda = \frac{\pi_\lambda\pi_\lambda^\top}{\prod_{j\in [n]} \#\{i \in [|\lambda|] \ |\ \lambda_i = j\}!}$ est le nombre de façons de partitionner la permutation $\sigma^{-1}\tau$ en suites croissantes de longueurs $\lambda_1, \lambda_2, \ldots, \lambda_{|\lambda|}$.
\end{defn}
Tout au long de la preuve qui suit, nous utilisons $y(\lambda)$ pour désigner le nombre de permutations des blocs qui ne changent pas la taille des blocs (dans l'ordre), c'est-à-dire $y(\lambda)=\prod_{j\in [n]} \#\{i \in [|\lambda|] \ |\ \lambda_i = j\}!$.

\begin{proof}[Preuve de l'équivalence des définitions]
	Par la définition de $\pi_\lambda$, l'entrée $(\sigma, \tau)$ de la matrice est égale au nombre de compositions $c$ de forme $\lambda$ telles que $\sigma \ast c = \tau$. Rappelons que $c$ agit à droite par le produit de composition sur les permutations.
	Remarquons que, pour un certain partage $\lambda$, il y a $y(\lambda)$ fois plus de compositions ensemblistes de forme $\lambda$ que de façons de partitionner les nombres de $1$ à $n$ en ensembles de taille $\lambda_1$ à $\lambda_{|\lambda|}$. En effet, ce ratio est le nombre de façons de permuter les blocs de la composition ensembliste.
	
	Pour simplifier la notation, nous notons $\pi$ pour désigner $\pi_\lambda$. 
	Ainsi,
	\allowdisplaybreaks
	\begin{align*}
	[\nu_\lambda]_{\sigma,\tau} & = \frac{[\pi\pi^\top]_{\sigma,\tau}}{y(\lambda)}\\
	& = \sum_{\omega \in S_n} \frac{\pi_{\sigma,\omega} \pi_{\omega,\tau}^\top}{y(\lambda)}\\
	& = \sum_{\omega \in S_n} \frac{\pi_{\sigma,\omega} \pi_{\tau,\omega}}{y(\lambda)}\\
	& = \sum_{\omega \in S_n} \frac{\#\{c \text{ de forme }\lambda \mid \sigma \ast c = \omega = \tau \ast c \}}{y(\lambda)}\\
	& = \frac{\#\{c \text{ de forme }\lambda \mid \sigma\ast c = \tau\ast c \}}{y(\lambda)}\\
	& = \frac{\#\{c \text{ de forme }\lambda \mid \Id \ast c = \sigma^{-1}\tau \ast c \}}{y(\lambda)}\\
	& = \#\left\{\begin{array}{l|l}
		i_{1,1} < i_{1,2} < \ldots < i_{1,\lambda_1} & \tau^{-1}\sigma(i_{1,1}) < \ldots < \tau^{-1}\sigma(i_{1,\lambda_1})\\
		\qquad \vdots & \qquad \vdots\\
		i_{|\lambda|,1} < i_{|\lambda|,2} < \ldots < i_{|\lambda|,\lambda_{|\lambda|}} & \tau^{-1}\sigma(i_{|\lambda|,1}) < \ldots < \tau^{-1}\sigma(i_{|\lambda|,\lambda_1})\\
	\end{array}\right\}\\
	& = \#\left\{\begin{array}{l}\text{ Façons de partitionner la permutation $\sigma^{-1}\tau$ en suites}\\
	\text{ croissantes de longueurs }\lambda_1, \lambda_2, \ldots, \lambda_{|\lambda|}	
	\end{array}\right\}.
	\end{align*}
	Notons que l'avant-dernière égalité est déduite de la remarque faite précédemment sur le nombre de façons de permuter les blocs de même taille et du fait que les nombres dans un même bloc doivent être présentés dans le même ordre dans $\sigma^{-1}\tau$ et dans la permutation identité.
\end{proof}

Pour simplifier les notations, on note le nombre de façons de partitionner la permutation $\sigma$ en suites croissantes de longueurs $\lambda_1, \lambda_2, \ldots, \lambda_{|\lambda|}$ par $\noninv_\lambda(\sigma)$.\index{$\noninv_\lambda$} Ainsi, $\noninv_{i}(\sigma)$, tel que présenté à la \autoref{defn:nu_k_noninv_perm}, correspond aussi à $\noninv_{(i, 1, 1, \ldots, 1)}(\sigma)$. Cela nous permet de définir les opérateurs $\{\nu_\lambda\}_{\lambda \vdash n}$.

\begin{defn}
	Soit $\sigma$ une permutation de $n$. L'opérateur de mélange symétrisé $\nu_\lambda$ peut être défini en termes de suites croissantes de la façon suivante~:
	\[ v_\lambda(\sigma) = \sum_{\tau \in S_n} \noninv_\lambda(\tau) \sigma \cdot \tau. \]
\end{defn}
Cette définition des opérateurs $\{\nu_\lambda\}_{\lambda \vdash n}$ est notamment compatible avec les opérateurs $\allnuk$ et peut être élargie aux mots de la même façon~:
\begin{defn}
	Soit $w$ un mot de longueur $n$. L'opérateur de mélange symétrisé $\nu_\lambda$ agit à droite sur $w$ ainsi~:
	\[ \nu_\lambda(w) = \sum_{\tau \in S_n} \noninv_\lambda(\tau) w \cdot \tau. \]
\end{defn}

On est alors en mesure de définir les opérateurs $\allgammak$~:
\begin{defn}\label{defn:op_gamma_k}
	L'opérateur de mélange symétrisé $\gammak{k}$ est associé au partage $(2^k,1^{n-2k})$ de $n$~:
	\[\gammak{k}=\nu_{(2^k, 1^{n-2k})}. \]
\end{defn}

\subsection{Matrices}\label{ssec:gammak_matrices}
Comme expliqué à la \autoref{defn:entrees_mat_v_lambda}, la valeur de la matrice $\nu_\lambda$ à la position $(\sigma, \tau)$ est le nombre de façons de partitionner $\sigma^{-1}\tau$ en suites croissantes de taille $\lambda_1, \ldots, \lambda_{|\lambda|}$. En particulier, pour $\gamma_k = \nu_{(2^k, 1^{n-2k})}$, on doit partitionner $\sigma^{-1}\tau$ en $k$ suites croissantes de taille $2$ et en $n-2k$ singletons. 

Cette définition est donnée de sorte que, si $\vec{v}$ est le vecteur formé de toutes les permutations (dans le même ordre que les lignes de la matrice), alors la colonne $\sigma$ de $\vec{v} \cdot M_\lambda$ vaut $\nu_\lambda(\sigma)$.

Grâce à la structure de ces matrices, on peut déduire certaines de leurs propriétés.
\paragraph{Matrices entières}
Les entrées de la matrice étant le nombre de partitions des nombres de $1$ à $n$, les valeurs de la matrices sont évidemment entières (et positives).
\paragraph{Matrices symétriques}
La matrice de $\gammak{k}$ est un multiple de la matrice $\pi\pi^\top$, lorsque $\pi$ est un opérateur de \bhr\ (décrit à la \autoref{ssec:gammak_bhr}). Or,
\mbox{$(\pi\pi^\top)^\top = \pi\pi^\top$}, et c'est donc une matrice symétrique.
\paragraph{Matrices semi-définies positives} De la même façon que les opérateurs de mélange symétrisés de la première famille correspondent à des matrices semi-définies positives, les opérateurs de la deuxième famille le sont également. En effet, l'argument-clé pour montrer qu'elles sont semi-définies positives est que c'est le cas pour les matrices qui s'écrivent sous la forme $\pi\pi^\top$. On a vu, à la \autoref{ssec:gammak_bhr}, que c'était le cas pour les opérateurs $\allgammak$. Pour un rappel de la preuve complète, voir à la \autopageref{par:nuk_semi_def_pos}.

\paragraph{Matrices stochastiques, après renormalisation}
Les matrices des opérateurs $\allgammak$ sont construites de façon à ce que toutes leurs entrées soient positives et que la somme d'une ligne soit la même pour chaque ligne. Ce sont donc des multiples de matrices stochastiques.

On rappelle que la valeur à la position $(\sigma,\tau)$ dans la matrice associée à $\gammak{k}$ correspond au nombre de façons d'obtenir $\tau$ à partir de $\sigma$. La somme de chaque ligne est donc le nombre d'événements pouvant se produire, avec multiplicités.  Or, on a vu à la \autoref{ssec:gammak_bhr} qu'il y avait $\frac{(n!)^2}{(2^k)^2 \cdot k! \cdot (n-2k)!}$ suites d'événements possibles; ce nombre est donc égal à la somme des valeurs sur une ligne.


\paragraph{Famille de matrices qui commutent}La propriété la plus intéressante des matrices $\allgammak$ est qu'elles commutent entre elles. Cette propriété se distingue des autres propriétés énoncées jusqu'ici qui sont vraies en général pour les opérateurs $\nu_\lambda$. 
Prouvé dans \cite{RSW} (Théorème I.4.3) en utilisant plusieurs notions qui ne sont pas définies ici, ce résultat a d'ailleurs mené à la conjecture suivante~:
\begin{cj}[Conjecture I.4.4 de \cite{RSW}]
	Soit $\lambda$ et $\mu$ deux différents partages de $n$, qui sont distincts de $(n)$ et de $(1^n)$. Les opérateurs $\nu_\lambda$ et $\nu_\mu$ commutent si et seulement s'ils appartiennent tous les deux à la première famille d'opérateurs de mélange symétrisés, $\allnuk$, ou s'ils appartiennent tous les deux à la deuxième famille d'opérateurs, $\allgammak$.
\end{cj}

\section{Valeurs propres}
Le travail de Victor Reiner, Franco Saliola et Volkmar Welker a aussi permis de comprendre les valeurs propres des opérateurs $\allgammak$. Le résultat principal à cet égard est que celles-ci sont des entiers (Théorème I.4.3). Cela est particulièrement intéressant, puisque les auteurs de ce théorème ont aussi conjecturé que seuls les opérateurs de mélange symétrisés qui appartiennent aux familles $\allnuk$ et $\allgammak$ ont des valeurs propres entières. Ils ont de plus donné une formule pour calculer explicitement toutes les valeurs propres de $\gammak{k}$ en utilisant les caractères irréductibles du groupe symétrique.
\subsection{Avec les caractères}\label{ssec:gammak_caracteres}
On connaît déjà une façon de décrire toutes les valeurs propres des matrices $\gamma_k$. Celle-ci est donnée par le théorème 5.2 de \cite{RSW}~:
\begin{thm}[Théorème 5.2 de \cite{RSW}]\label{thm:vp_gammak_rsw}
	On peut décomposer $\mathbb{R}S_n$ en une somme directe de modules orthogonaux~:
	\[ \mathbb{R}S_n = K \bigoplus_{\lambda \vdash n} U^\lambda, \]
	où $K$ est dans le noyau de tous les opérateurs $\allgammak$ et $U^\lambda$ est compris dans un espace propre de chacun des opérateurs $\allgammak$. Pour une certaine valeur de $k$, l'unique valeur propre de $\gammak{k}|_{U^\lambda}$ est notamment déterminée par le caractère irréductible $\chi^\lambda$ associé à ce partage. Plus précisément, cette valeur propre est
	\[c_k^\lambda = \sum_{\sigma \in S_n} \noninv_{(2^k, 1^{n-2k})}(\sigma) \cdot \chi^\lambda(\sigma).\index{$c_k^\lambda$}\]
\end{thm}

Ce théorème nous permet entre autres de dire que les valeurs propres des opérateurs $\allgammak$ sont des entiers. Rappelons, à ce titre, que les caractères du groupe symétrique sont des entiers (voir le \autoref{corl:caracteres_Sn_entiers} pour plus de détails).

Ce qui est plus difficile à voir dans le dernier théorème est que  
les valeurs propres sont toutes positives~: ce sont en effet les valeurs propres d'un matrice semi-définie positive (comme démontré à la \autoref{ssec:gammak_matrices}). Nous avons vu (\autoref{ssec:positivite}) que ces matrices ont des valeurs propres positives.

Par contre, déduire la positivité des valeurs propres à partir du \autoref{thm:vp_gammak_rsw} n'est pas évident, puisque les caractères du groupe symétrique peuvent aussi être négatifs. C'est pourquoi on aimerait avoir une autre façon d'exprimer $c_k^\lambda$. Victor Reiner, Franco Saliola et Volkmar Welker souhaitaient déjà, dans \cite{RSW}, qu'il existe une expression explicite des valeurs propres qui mettrait en évidence la positivité (remarque V.2.3). La \autoref{ssec:obs_vp_gammak} explique quelques conjectures à cet effet.

\subsection{Formules récursives pour les valeurs propres non-nulles}
\label{ssec:obs_vp_gammak}

Pour la première famille d'opérateurs de mélange symétrisés, $\allnuk$, on a trouvé une façon d'exprimer les valeurs propres de l'opérateur $\nu_k$ sur $S_n$ en fonction de $\nu_{k-1}$ et de $\nu_k$ sur $S_{n-1}$. Ceci nous permettait donc de  calculer rapidement toutes les valeurs propres de tous les opérateurs de la famille, en plus de démontrer que toutes les valeurs propres sont des entiers positifs.

On s'est donc demandé s'il pouvait exister une façon similaire d'exprimer les valeurs propres des opérateurs $\allgammak$. Expérimentalement, nous pensons avoir trouvé une expression récursive pour certaines valeurs propres. Les observations les plus intéressantes ont été faites sur les valeurs propres lorsqu'on se restreint aux modules de Specht $S^{(n-m,1^m)}$. Ce sont les modules indexés par les diagrammes en forme d'équerre et qui sont associés à certaines des représentations les plus connues (triviale, standard et signe). Pour ces modules, les valeurs $c_k^\lambda$ sont affichées aux tableaux \ref{tab:vp_gammak_triv} à \ref{tab:vp_gammak_m3} pour $n$ entre $4$ et $8$. Certaines relations entre les valeurs propres sont mises en évidence dans ces tableaux. Dans ceux-ci, les lignes représentent des modules et, sur une même colonne, on retrouve différentes valeurs propres d'un même opérateur. Des relations ont lieu autant entre les valeurs propres de différents opérateurs sur un même module (ce sont les flèches horizontales) et entre les valeurs propres du même opérateur, sur différents modules. Elles sont données par les flèches verticales.

\begin{table}[H]
	\caption[Valeurs propres des opérateurs $\gamma_k$ sur les modules $S^{(n)}$.]{Valeurs propres des opérateurs $\gamma_k$ sur les modules $S^{(n)}$, associés à la représentation triviale, pour $n$ entre $2$ et $8$. Les flèches verticales arrivant à la ligne $n$ sont étiquetées par $\frac{n^2}{n-2k}$, et les flèches horizontales arrivant dans la colonne $k$, par des produits d'un nombre triangulaire et de $\frac{1}{2k}$.}\label{tab:vp_gammak_triv}
	\[\begin{tikzcd}
		k & 0 & 1 & 2 & 3 & 4 \\[-15pt]
		S^{(2)} & 2 \arrow[r, "\times\frac{1}{2}", color=magenta] \arrow[d, dotted, color=blue,"\times3"]  & 1 \arrow[d, dotted, color=blue, "\times9"] \\
		S^{(3)} & 6 \arrow[r, "\times\frac{3}{2}", color=magenta] \arrow[d, dotted, color=blue,"\times4"]  & 9 \arrow[d, dotted, color=blue, "\times8"] \\
		S^{(4)} & 24 \arrow[r, "\times\frac{6}{2}", color=magenta] \arrow[d, dotted, color=blue,"\times5"]  & 72 \arrow[r, color=magenta, "\times\frac{1}{4}"] \arrow[d, dotted, color=blue, "\times\frac{25}{3}"]  & 18 \arrow[d, dotted, color=blue,"\times25"] \\
		S^{(5)} & 120 \arrow[r, color=magenta, "\times\frac{10}{2}"] \arrow[d, dotted, color=blue, "\times6"] &   600 \arrow[r, color=magenta, "\times\frac{3}{4}"] \arrow[d, dotted,color=blue,"\times9"] & 450 \arrow[d, dotted,color=blue,"\times18"] &   \\
		S^{(6)} & 720 \arrow[r, color=magenta, "\times\frac{15}{2}"]\arrow[d, dotted,color=blue,"\times7"]  & 5400 \arrow[r, color=magenta,"\times\frac{6}{4}"]\arrow[d, dotted, color=blue, "\times\frac{49}{5}"] & 8100 \arrow[r, color=magenta,"\times\frac{1}{6}"]\arrow[d, dotted,color=blue,"\times\frac{49}{3}"] & 1350\arrow[d, dotted,color=blue,"\times49"]\\
		S^{(7)} & 5040 \arrow[r, color=magenta, "\times\frac{21}{2}"] \arrow[d, dotted,color=blue,"\times8"] & 52920\arrow[r, color=magenta, "\times\frac{10}{4}"] \arrow[d, dotted,color=blue,"\times\frac{32}{3}"] & 132300 \arrow[r, color=magenta,"\times\frac{3}{6}"] \arrow[d, dotted,color=blue,"\times16"]& 66150 \arrow[d, dotted,color=blue,"\times32"]\\
		S^{(8)} & 40320 \arrow[r, color=magenta, "\times\frac{28}{2}"] & 564480 \arrow[r, color=magenta, "\times\frac{15}{4}"] & 2116800\arrow[r, color=magenta, "\times\frac{6}{6}"] &  2116800 \arrow[r, color=magenta, "\times\frac{1}{8}"] & 264600
	\end{tikzcd}\]
\end{table}

\begin{table}[H]
	\caption[Valeurs propres des opérateurs $\gamma_k$ sur les modules $S^{(n-1,1)}$.]{Valeurs propres des opérateurs $\gamma_k$ sur les modules $S^{(n-1,1)}$, associés à la représentation standard, pour $n$ entre $4$ et $8$. Les flèches verticales vers la ligne $n$ sont étiquetées par $\frac{(n-2)(n+1)}{n-2k}$, et les flèches horizontales arrivant dans la colonne $k$, par des produits d'un nombre triangulaire et de $\frac{1}{2(k-1)}$.}\label{tab:vp_gammak_reg}
	\[\begin{tikzcd}
	k &  1 & 2 & 3 & 4 \\[-15pt]
	S^{(3,1)} & 20 \arrow[r, "\times\frac{1}{2}", color=magenta] \arrow[d, dotted, color=blue,"\times6"]  & 10 \arrow[d, dotted, color=blue, "\times18"] \\
	S^{(4,1)} & 120 \arrow[r, color=magenta, "\times\frac{3}{2}"] \arrow[d, dotted, color=blue, "\times7"] & 180 \arrow[d, dotted,color=blue,"\times14"]\\
	S^{(5,1)} & 840 \arrow[r, color=magenta, "\times\frac{6}{2}"]\arrow[d, dotted,color=blue,"\times8"]  & 2520 \arrow[r, color=magenta,"\times\frac{1}{4}"]\arrow[d, dotted, color=blue, "\times\frac{40}{3}"] & 630 \arrow[d, dotted,color=blue,"\times40"]\\
	S^{(6,1)} & 6720 \arrow[r, color=magenta, "\times\frac{10}{2}"]\arrow[d, dotted,color=blue,"\times9"]  & 33600 \arrow[r, color=magenta,"\times\frac{3}{4}"]\arrow[d, dotted, color=blue, "\times\frac{27}{2}"] & 25200 \arrow[d, dotted,color=blue,"\times27"]\\
	S^{(7,1)} & 60480 \arrow[r, color=magenta, "\times\frac{15}{2}"] & 453600\arrow[r, color=magenta, "\times\frac{6}{4}"] & 680400 \arrow[r, color=magenta,"\times\frac{1}{6}"] & 113400
	\end{tikzcd}\]
\end{table}

\begin{table}[H]
	\caption[Valeurs propres des opérateurs $\gamma_k$ sur les modules $S^{(n-2,1,1)}$.]{Valeurs propres des opérateurs $\gamma_k$ sur les modules $S^{(n-2,1,1)}$ pour $n$ entre $4$ et $8$. Les flèches verticales arrivant à la ligne $n$ sont étiquetées par $\frac{n(n-2)}{n-2k}$, et les flèches horizontales arrivant dans la colonne $k$, par des produits d'un nombre triangulaire et de $\frac{1}{2(k-1)}$.}
	\[\begin{tikzcd}
	k & 1 & 2 & 3 & 4 \\[-15pt]
	S^{(2,1,1)} & 4 \arrow[r, "\times\frac{1}{2}", color=magenta] \arrow[d, dotted, color=blue,"\times5"]  & 2 \arrow[d, dotted, color=blue, "\times15"]\\
	S^{(3,1,1)} & 20 \arrow[r, color=magenta, "\times\frac{3}{2}"] \arrow[d, dotted, color=blue, "\times6"] &   30 \arrow[d, dotted,color=blue,"\times12"]  \\
	S^{(4,1,1)} & 120 \arrow[r, color=magenta, "\times\frac{6}{2}"]\arrow[d, dotted,color=blue,"\times7"]  & 360 \arrow[r, color=magenta,"\times\frac{1}{4}"]\arrow[d, dotted, color=blue, "\times\frac{35}{3}"] & 90 \arrow[d, dotted,color=blue,"\times35"]\\
	S^{(5,1,1)} & 840 \arrow[r, color=magenta, "\times\frac{10}{2}"]\arrow[d, dotted,color=blue,"\times8"]  & 4200 \arrow[r, color=magenta,"\times\frac{3}{4}"]\arrow[d, dotted, color=blue, "\times12"] & 3150 \arrow[d, dotted,color=blue,"\times24"]\\
	S^{(6,1,1)} & 6720 \arrow[r, color=magenta, "\times\frac{15}{2}"] & 50400\arrow[r, color=magenta, "\times\frac{6}{4}"] & 75600 \arrow[r, color=magenta,"\times\frac{1}{6}"] & 12600
	\end{tikzcd}\]
\end{table}

\begin{table}[H]
	\caption[Valeurs propres des opérateurs $\gamma_k$ sur les modules $S^{(n-3,1,1,1)}$.]{Valeurs propres des opérateurs $\gamma_k$ sur les modules $S^{(n-3,1,1,1)}$ pour $n$ entre $4$ et $8$. Les flèches verticales arrivant à la ligne $n$ sont étiquetées par $\frac{(n-4)(n+1)}{n-2k}$, et les flèches horizontales arrivant dans la colonne $k$, par des produits d'un nombre triangulaire et de $\frac{1}{2(k-2)}$.}\label{tab:vp_gammak_m3}
	\[\begin{tikzcd}
		k & 1 & 2 & 3 & 4 \\[-15pt]
		S^{(1,1,1,1)} & 0 & 2 \arrow[d, dotted, color=blue, "\times6"]\\
		S^{(2,1,1,1)} & 0 & 12 \arrow[d, dotted,color=blue,"\times7"]  \\
		S^{(3,1,1,1)} & 0 & 84 \arrow[r, color=magenta,"\times\frac{1}{2}"]\arrow[d, dotted, color=blue, "\times8"] & 42 \arrow[d, dotted,color=blue,"\times24"]\\
		S^{(4,1,1,1)} & 0 & 672 \arrow[r, color=magenta,"\times\frac{3}{2}"]\arrow[d, dotted, color=blue, "\times9"] & 1008 \arrow[d, dotted,color=blue,"\times18"]\\
		S^{(5,1,1,1)} & 0 & 6048 \arrow[r, color=magenta, "\times\frac{6}{2}"] & 18144 \arrow[r, color=magenta,"\times\frac{1}{4}"] & 4536
	\end{tikzcd}\]
\end{table}

\subsubsection{Sur les modules de Specht associés aux diagrammes en forme d'équerre}

Les tableaux \ref{tab:vp_gammak_triv} à \ref{tab:vp_gammak_m3} montrent l'unique valeur propre non-nulle (lorsqu'elle existe) de $\gammak{k}|_{S^{(n-m,1^m)}}$ pour $n$ entre $4$ et $8$. Les seuls tableaux représentés sont pour $m$ de $0$ à $3$. Au-delà de ce nombre, la quantité de données facilement calculables dans le tableau est trop faible pour vraiment voir des relations entre elles. Rappelons que les données deviennent rapidement difficiles à obtenir lorsque $n$ croît.

Après avoir regardé ces tableaux, il semble que certaines soient multiples des autres. En faisant le ratio des valeurs sur les mêmes colonnes d'un même tableau, nous trouvons les valeurs indiquées à côté des flèches verticales. Ceci mène à la conjecture suivante~:

\begin{cj}\label{cj:vp_gammak_rec}
	La  valeur propre non-nulle $c_k^{(n-m, 1^m)}$ de l'opérateur $\gammak{k}|_{S^{(n-m, 1^m)}}$ satisfait l'équation~:
	\[\frac{c_k^{(n-m, 1^m)}}{c_k^{(n-m-1, 1^m)}} = \left\{\begin{array}{ll}
	\frac{(n-m)n}{n-2k}= \frac{\lambda_1 n}{n-2k} & \text{si $m$ est pair;}\\
	\frac{(n-m-1)(n+1)}{n-2k}= \frac{(\lambda_1-1)(n+1)}{n-2k} & \text{si $m$ est impair.}
	\end{array}\right.  \]
	
	
\end{cj}

Notons toutefois que cette dernière formule n'est pas tout à fait élégante.
On aimerait pouvoir expliquer les divergences entre les cas pairs et impairs pour $m$. Or, comme les valeurs apparaissent dans des tableaux différents pour chaque $m$, dans notre démarche, nous arrivons difficilement à voir pourquoi c'est le cas. C'est toutefois quelque chose qui devrait être exploré.
Enfin, même si la \autoref{cj:vp_gammak_rec} était prouvée, il resterait à trouver au moins une valeur propre pour $\gammak{k}$ sur chaque module $S^{(n-m,1^m)}$.

On aimerait donc aussi connaître des liens entre les valeurs propres des différents opérateurs $\allgammak$ sur les collections de même taille (les différents éléments d'une même ligne des tableaux \ref{tab:vp_gammak_triv} à \ref{tab:vp_gammak_m3}). C'est d'ailleurs l'objet de la prochaine observation~:
\begin{cj}\label{cj:gammak_vp_horiz}
	En regardant les valeurs dans les tableaux \ref{tab:vp_gammak_triv} à \ref{tab:vp_gammak_m3}, il semble que les valeurs sur une même ligne soient des multiples l'une de l'autre. Plus spécifiquement, si la colonne indexée par $k$ est la $i$-ième à contenir des valeurs propres non-nulles sur $S^{(n-m,1^m)}$, alors
	\[ \frac{c_{k+1}^{(n-m,1^m)}}{c_{k}^{(n-m,1^m)}} = \frac{\binom{n-2k}{2}}{2i} \]
	pour chacun des cas présents dans les tableaux \ref{tab:vp_gammak_triv} à \ref{tab:vp_gammak_m3}. 
\end{cj}

\begin{rmq}
	La \autoref{cj:gammak_vp_horiz}, contrairement à la \autoref{cj:vp_gammak_rec}, semble vraie sur tous les modules simples, et pas seulement sur $S^{(n-m,1^m)}$.
\end{rmq}
\subsubsection{Sur les modules de Specht $S^{(n)}$}
Le cas pour lequel nous savons vérifier la \autoref{cj:vp_gammak_rec} pour l'instant est celui du module $S^{(n)}$. La raison pour cela est que ce module simple correspond à la représentation triviale; le caractère associé est donc la fonction constante $1$. En utilisant le \autoref{thm:vp_gammak_rsw}, on trouve~:
\begin{align*}
c^{(n)}_k = & \sum_{w \in S_n} \noninv_{(2^k,  1^{n-2k})}(w) \underbrace{\chi^{(n)}(w)}_{=1}\\
= &\sum_{w \in S_n} \noninv_{(2^k, 1^{n-2k})}(w).
\end{align*}

Cette dernière valeur est égale à la somme des éléments sur la première ligne de la matrice associée à $\gamma_k$. On a démontré, à la \autoref{ssec:gammak_bhr}, que cette valeur était 
$\frac{\binom{n}{2^k, 1^{n-2k}}^2}{k!(n-2k)!}$.

Cela concorde avec la \autoref{cj:vp_gammak_rec}. On disait alors que 
\[c_k^{(n)} = \frac{n^2}{n-2k} \cdot c_k^{(n-1)}. \]
C'est bien ce que l'on obtient ici~:
\[\frac{c_k^{(n)}}{c_k^{(n-1)}} = \frac{\binom{n}{(2^k, 1^{n-2k})}^2}{(n-2k)! k!} \cdot \frac{(n-1-2k)!k!}{\binom{n-1}{(2^k, 1^{n-1-2k})}^2} = \frac{n^2}{n-2k}.\]

\subsection{L'opérateur $\gamma_1$ et $\nu_{n-2}$}
Enfin, notons qu'il existe un opérateur qui est commun aux deux familles. Il s'agit de $\gammak{1} = \nu_{n-2}$. Cela nous permet, par exemple, de justifier que, sur $S^{(n-3,1^3)}$, il n'y ait aucune valeur propre non-nulle pour cet opérateur (dans le tableau \ref{tab:vp_gammak_m3}). En effet, il n'y a aucun tableau de type au moins $n-2$ qui soit de forme $(n-3,1^3)$, et la théorie sur les opérateurs de la première famille implique que toute copie de $S^{(n-3,1^3)}$ est dans le noyau de $\Nu{n-2}$, et donc de $\gamma_1$.

De plus, si nous étions capables de prouver la \autoref{cj:gammak_vp_horiz} et de mieux expliquer les liens entre $c_k^{(n-m,1^m)}$ et $c_{k+1}^{(n-m,1^m)}$, l'égalité de $\gamma_1$ et de $\Nu{n-2}$ pourrait être très intéressante. En effet, les valeurs propres de $\nu_{n-2}$ sont toutes données par le \autoref{thm:main}. En plus de nous permettre de calculer les valeurs propres de $\gamma_1$, mieux comprendre les liens entre les valeurs propres des différents opérateurs pourrait nous permettre de décrire toutes les valeurs propres des $\allgammak$ sur les espaces qui ne sont pas compris dans le noyau de $\gamma_1$.

\begin{conclusion}
	Cette thèse porte essentiellement sur deux familles d'opérateurs de mélange symétrisés qui commutent entre eux. Pour ces deux familles, on sait calculer les valeurs propres (par des moyens différents)~: dans le premier cas, c'est le résultat du \autoref{thm:main}, alors que, dans le second, on sait le faire depuis \cite{RSW}. Bien qu'on connaisse des algorithmes pour calculer ces valeurs, plusieurs questions restent irrésolues. Certaines sont listées ici.
	
	\subsection*{La première famille~: le mélange doublement aléatoire et les mélanges analogues}
	La procédure pour calculer les valeurs propres est récursive~: on obtient celles de $\nu_k$ à partir de celles de $\nu_{k-1}$.  C'est intéressant, notamment parce que ça nous permet de démontrer que toutes les valeurs propres sont entières, mais aussi parce que ça donne un algorithme rapide pour les calculer. Cependant, certaines choses ne sont pas faciles à faire avec une telle façon d'énumérer les valeurs propres. Par exemple, faire la somme de toutes les valeurs propres (pour calculer le temps de mélange) s'avère très ardu. Ceci suggère le problème suivant~:
	\begin{pb}
		Trouver une formule close pour les valeurs propres des opérateurs 
		$\{\nu_k\}_{k\in \N}$.
	\end{pb} 
	Ceci n'est probablement pas facile. À défaut de trouver une formule close, on pourrait se satisfaire de bornes suffisamment proches des valeurs propres. Ceci pourrait être suffisant pour calculer le temps de mélange. Une conjecture sur une borne supérieure pour les valeurs propres est donnée à la \autopageref{cj:bornes_v_k}. Il faudrait toutefois vérifier que c'est suffisamment proche pour rendre les calculs intéressants. Si c'est le cas, il est possible que le temps de mélange puisse être calculé.
	
%
	\begin{pb}
		Calculer le temps de mélange des opérateurs $\{\nu_k\}_{k\in\N}$.
	\end{pb}

	Bien que les valeurs propres aient été toutes trouvées, nous ne savons pas décrire complètement les espaces propres. À partir d'une base du noyau, nous savons décrire les autres espaces propres (par l'opérateur $\projlift$, décrit au \autoref{chap:preuves}). Toutefois, connaître une base du noyau est suffisant pour décrire tous les vecteurs propres des opérateurs $\allnuk$.
	\begin{pb}
		Donner une base du noyau des opérateurs $\allnuk$. Plus généralement, décrire une base des espaces propres des opérateurs.
	\end{pb}
	
	\subsection*{L'autre famille d'opérateurs}
	En plus de la famille d'opérateurs $\allnuk$, Reiner, Saliola  et Welker ont identifié une autre famille d'opérateurs $\{\gamma_k\}_{k\in\N}$ qui commutent et dont les valeurs propres sont entières \cite{RSW}. Ces dernières peuvent être explicitement calculées par une formule faisant intervenir les caractères  du groupe symétrique. Toutefois, cette formule ne permet pas, par exemple, d'expliquer pourquoi les valeurs propres sont toujours positives. C'est pourquoi, au \autoref{chap:2efamille}, plusieurs conjectures sont proposées pour expliquer les relations entre ces valeurs propres et, notamment, mettre en évidence le fait qu'elles sont positives.
	\begin{pb}\label{pb:formule_close}
		Trouver une explication combinatoire aux valeurs propres des 
		opérateurs $\allgammak$ et résoudre les conjectures du \autoref{chap:2efamille}.
	\end{pb}

	\subsection*{Plus d'opérateurs de mélanges symétrisés}
	Les opérateurs $\allnuk$ et $\allgammak$ sont des opérateurs de mélange symétrisés, mais ils sont loin d'être les seuls. Toutefois, si nous nous intéressons à ces deux familles, c'est qu'elles sont telles que les opérateurs d'une même famille commutent, ce qui n'est pas le cas des opérateurs de mélange symétrisés en général (ils sont décrits à la \autoref{ssec:suites_croissantes_gamma_k} et dans \cite{RSW}). La conjecture I.4.4 de \cite{RSW} détaille l'intuition de ses auteurs~: elle dit que si deux opérateurs de mélange symétrisés commutent, ils sont soit tous les deux dans la famille $\allnuk$ ou tous les deux dans la famille $\allgammak$, sauf quelques cas triviaux. De plus, elle dit que si un opérateur n'a que des valeurs propres entières, il doit être dans une de ces deux familles. Pourtant, les techniques, pour démontrer la commutativité et pour démontrer que les valeurs propres sont entières sont très différentes d'une famille à l'autre, et il est difficile de comprendre pourquoi la conjecture devrait être vraie.
	\begin{pb}
		Comprendre \textit{pourquoi} les valeurs propres des autres opérateurs de mélange symétrisés (ceux qui ne sont pas dans $\allnuk$, ni dans $\allgammak$) ne sont pas entières.
	\end{pb}
	\begin{pb}
		Donner une explication combinatoire de la commutativité des opérateurs $\allnuk$ et des opérateurs $\allgammak$.
	\end{pb}
	
	\subsection*{Sur d'autres structures}
	Les opérateurs décrits plus haut agissent sur le groupe symétrique. C'est utile de considérer ce cas parce que la théorie des modules y est bien développée. Il serait intéressant de s'attarder à ce qui se produit sur d'autres structures algébriques pour généraliser le problème. Je pense notamment aux monoïdes dont la théorie de la représentation et les relations avec les chaînes de Markov ont été développées dans les dernières années \cite{ASST}, mais aussi aux extensions linéaires de posets, sur lesquelles le mélange doublement aléatoire a déjà été étudié \cite{AST}.
	
	Même parmi les groupes, on pourrait avoir des surprises. Par exemple, si on remplace le groupe symétrique par le groupe hyperoctaédral, les opérateurs de mélange ont aussi une interprétation en termes de mélange de cartes~: en plus de mélanger les cartes, on peut choisir de les retourner, individuellement, ou pas.
	
	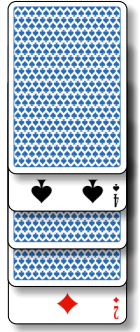
\begin{figure}
		\centering
		\begin{tikzpicture}
		\node (A) at (0,0) {\carda};
		\node (B) at (0,0.5) {\versocarte};
		\node (C) at (0,1) {\versocarte};
		\node (D) at (0,1.5) {\cardd};
		\node (E) at (0,2) {\versocarte};
		\end{tikzpicture}
		\caption{Un paquet avec des cartes à l'endroit et à l'envers représente un élément du groupe hyperoctaédral.}
	\end{figure}
	Il serait intéressant de savoir lesquels, parmi les résultats de cette thèse, sont toujours vrais sur cette structure. Par exemple, il serait intéressant de savoir si des opérateurs similaires à ceux des première et deuxième familles ont toujours des valeurs propres entières. Les calculs préliminaires suggèrent que cela est vrai.
	
\end{conclusion}

\appendix
\chapter{Théorèmes du chapitre 4}\label{chap:annexe_thm}

\paragraph{\Cref{thm:regle_de_young_mod_perm} ou règle de Young pour les modules de permutation} 
	Les modules de permutations se décomposent en modules simples de la façon suivante~:
	\[M^\mu \cong \bigoplus_{\lambda \unrhd \mu} m_{\lambda\mu} S^\lambda,\]
	où $m_{\lambda\mu}$ est le nombre de tableaux semi-standards de forme $\lambda$ et de contenu $\mu$. Pour la définition de tableaux semi-standards, voir à la \autopageref{par:SSYT}.

\paragraph{\Cref{prop:commutShDel}}
	Les opérateurs d'insertion commutent entre eux, et il en est de même pour les opérateurs de suppression~:
	\begin{align*}
	\sh_a\circ \sh_b &= \sh_b \circ \sh_a,\\
	\del_a \circ \del_b &= \del_b \circ \del_a.
	\end{align*}

\paragraph{\Cref{lem:lem11_thet}}
	L'opérateur de remplacement, $\theta_{a,b}$, est un morphisme de $S_n$-modules à droite. \\
	En particulier, $\theta$ commute avec les actions à droite de $S_n$.

\paragraph{\Cref{prop:defn_op}}
	À l'aide des opérateurs $\sh$ et $\del$, on peut donner cette nouvelle définition des opérateurs $\allnuk$ dans l'algèbre des mots~:
	\[ \Nu{k} = \sum_{1\leq a_1 \leq \ldots \leq a_k\leq n} \frac{1}{\prodnuk} \shdelUnak.  \]

\paragraph{\Cref{lem:lem35}}
	Sur l'algèbre des mots de longueur $n$, les identités suivantes sont vérifiées~:
	\begin{enumerate}[label=\roman*.]
		\item $\del_b \circ \sh_a - \sh_a\circ\del_b = \theta_{b,a} + (n+1)\delta_{a,b}\ \Id $ \label{lem35i},
		\item $\theta_{b,c} \circ \sh_a - \sh_a\circ\theta_{b,c} = \delta_{a,b}\sh_c $\label{regleDeNadia},
	\end{enumerate}
	où $\delta_{a,b}$ vaut $1$ si $a=b$ et $0$, sinon.

\paragraph{\Cref{lem:lem2}}
	Sur l'algèbre des mots de longueur $n$,
	\begin{align*}
		\sh_a \circ \Nu{k-1} = \sum_{i=1}^k \sum_{1 \leq a_1 \leq \ldots \leq a_k\leq n}& \delta_{a, a_i} \frac{1}{\prodnuk} \sh_a \circ \sh_{a_1} \circ \ldots \circ \nonumber\\ &\widehat{\sh_{a_{i}}}\circ \ldots\circ\sh_{a_k}\circ \del_{a_1} \circ\ldots \widehat{\del_{a_{i}}}\circ \ldots \circ \del_{a_k}.
	\end{align*}

\paragraph{\Cref{thm:lem8}}
	Sur l'algèbre des mots,
	\begin{align*}
		\sum_{a = 1}^{n}\sh_a \circ\ \Nu{k-1}\circ \del_a = k\ \Nu{k}.
	\end{align*}

\paragraph{\Cref{thm:thm1}}
	Sur les mots de longueur $n$,
	\[ \Nu{k}\circ \sh_a - \sh_a \circ\Nu{k} = (n+2-k) \sh_a \circ \Nu{k-1} + \sum_{1\leq b \leq n} \sh_b \circ \theta_{b,a}\circ\Nu{k-1}. \]

\paragraph{\Cref{lem:lem10_nu_k_Slambda->Slambda}}
	Pour toute valeur de $k$ et pour tout partage $\lambda$, l'image de $\nu_k|_{S^\lambda} $ est comprise dans $S^\lambda$.

\paragraph{\Cref{lem:lem43}}
	Soit $\lambda$ un partage.  Si $a<b$, alors $\theta_{b,a}|_{S^\lambda} = 0$.

\paragraph{\Cref{thm:thm2}}
	Soit $\lambda \vdash n$. Alors,
	\[ \left(\nu_k \circ \sh_a - \sh_a \circ \nu_k\right)|_{S^\lambda} = (n+2-k)\left(\sh_a \circ \Nu{k-1}\right)|_{S^\lambda} + \sum_{1\leq b \leq a} \left(\sh_b \circ \theta_{b,a}\circ \Nu{k-1}\right)|_{S^\lambda}.  \]

\paragraph{\Cref{prop:lem45DS}}
	L'image de $(\sh_a)|_{S^\lambda}$ est contenue dans un sous-module de $M^{\lambda+\vec{e_a}}$ isomorphe à $\bigoplus_{\mu = \lambda+\vec{e_r},\  r\leq a} S^\mu$.

\paragraph{\Cref{prop:proprietes_isoproj}}
	La fonction $\isoproj_\mu$ satisfait les énoncés suivants~:
	\begin{enumerate}
		\item elle commute avec les actions du groupe symétrique.
		\item c'est un morphisme de $S_n$-modules.
		\item c'est un projecteur isotypique de $\CSn$ vers la composante isotypique associée au module simple $S^\mu$.
	\end{enumerate}

\paragraph{\Cref{prop:1e_relation_orthogonalite}}
	Si $U$ et $V$ sont des modules simples, les caractères associés à $U$ et à $V$ sont orthonormaux, c'est-à-dire que
	\[ \langle \chi^U, \chi^V \rangle = \delta_{U,V} = \left\{ \begin{array}{ll}
	1 & \text{si $U = V$}\\
	0 & \text{si $U \neq V$}
	\end{array}\right.. \]

\paragraph{\Cref{lem:lem11_commut}}
	Le projecteur isotypique $\isoproj_\mu$ commute avec l'opérateur de remplacement $\theta_{a,b}$~:
	\[ \theta_{a,b}(\isoproj_\mu (v)) = \isoproj_\mu(\theta_{a,b}(v)). \]

\paragraph{\Cref{lem:lem13}}
	Soit $\lambda = (\lambda_1, \ldots, \lambda_l) \vdash n$, $a \in [l+1]$ et $\mu = \lambda + \vec{e_r}$, avec $r \in [a]$.  
	Alors,
	 \begin{align*}
		 \big(\nu_k \circ &\projlift_a^\mu - \projlift_a^\mu \circ \nu_k \big)|_{S^\lambda} =\\
		 &(n+3-k+\lambda_a - a) \left(\projlift_a^\mu \circ \Nu{k-1}\right)|_{S^\lambda} + \left(\sum_{r \leq b < a} \theta_{b,a} \circ \projlift_b^\mu \circ \Nu{k-1}\right)|_{S^\lambda}.
	 \end{align*}

\paragraph{\Cref{thm:thm3}}
	Soit $\lambda = (\lambda_1, \ldots, \lambda_l) \vdash n$, $a \in [l+1]$ et $\mu = \lambda + \vec{e_r}$, avec $r \in [a]$. Alors, 
	\[\left(\nu_k \circ \projlift_a^\mu - \projlift_a^\mu \circ \nu_k \right)|_{S^\lambda} = (n+2-k + (\lambda_r + 1- r)) \left(\projlift_a^\mu \circ \Nu{k-1}\right)|_{S^\lambda}.\]

\paragraph{\Cref{corl:corl4}}
	Soit $\lambda\vdash n$ et soit $v \in S^\lambda$ un vecteur propre à la fois pour $\Nu{k}$ et $\Nu{k-1}$, associé aux valeurs propres $v_k$ et $v_{k-1}$, respectivement. Alors, soit
	\begin{itemize}
		\item $\projlift_a^{\lambda+\vec{e_r}}(v) = 0$
		\item[ou] 
		\item $\begin{aligned}
		\Nu{k} \circ \projlift_a^{\lambda+\vec{e_r}}(v)& = \projlift_a^{\lambda+\vec{e_r}}\left( v_k + (n+2-k+(\lambda_r+1-r)) v_{k-1}\right) (v)\\ & = \left( v_k + (n+2-k+(\lambda_r+1-r)) v_{k-1}\right) \projlift_a^{\lambda+\vec{e_r}}(v).
		\end{aligned}$
	\end{itemize}

\paragraph{\Cref{thm:noyau_nu_k}}
Le noyau de $\Nu{k}$ est
\begin{equation*}
\ker(\Nu{k}) \cong \bigoplus_{\substack{\text{tableau  standard $t$ ,}\\ \type(t) < k}} S^{\text{forme}(t)}. 
\end{equation*}

\paragraph{\Cref{prop:im_nu_k_specht}}
	L'image de $\Nu{k}$ sur $S^\lambda$ est comprise dans la somme suivante~:
	\begin{align*}
	\im(\Nu{k}|_{S^\lambda}) \ &\subseteq \sum_{1 \leq a_1 \leq \ldots \leq a_k\leq n}  \im\left(\isoproj_\lambda \circ \sh_{a_k} \circ \ldots \circ \sh_{a_1}|_{S^{\lambda-\vec{e_{a_1}}-\ldots-\vec{e_{a_k}}}}\right)\\
	& \subseteq\ \sum_{1 \leq a \leq n} \im\left(\projlift_a^\lambda |_{S^{\lambda-\vec{e_a}}}\right).
	\end{align*}

\paragraph{\Cref{thm:thm5}}
	La propriété suivante concernant l'opérateur de suppression est vraie sur l'algèbre des mots de longueur $n$~:
	\[ \del_a \circ \nu_k - \nu_k \circ \del_a = (n+1-k) \Nu{k-1} \circ \del_a + \sum_{1\leq b \leq n} \theta_{a,b} \circ  \Nu{k-1} \circ \del_b. \]




\bibliographystyle{apalike-uqam}
\bibliography{bibliodoc}

\renewcommand\indexname{Index des notations}
\addtocontents{toc}{\protect\vspace{1.5ex}}
\printindex
\end{document}